\documentclass[10pt,a4paper]{article}

\usepackage{ifthen}	
\usepackage{amssymb}
\usepackage{amstext}
\usepackage{amsfonts}
\usepackage{amsmath}
\usepackage{amscd}
\usepackage{latexsym}
\usepackage{epsfig}
\usepackage{color}
\usepackage{epic}
\usepackage{pstricks,pst-node,pst-text,pst-tree}
\usepackage[all]{xy}
\usepackage{amsthm}
\usepackage{eucal,mathrsfs}
\usepackage[dvips, ps2pdf, draft]{hyperref}         

\usepackage{psfrag}

\newcommand{\numberset}[1]{\ensuremath{\mathbb{#1}}}    
\newcommand{\C}{\numberset{C}}  
\newcommand{\N}{\numberset{N}}  
\newcommand{\R}{\numberset{R}}  
\newcommand{\Z}{\numberset{Z}}  
\newcommand{\PP}{\numberset{P}}  

\newcommand{\zbar}{\ensuremath{\overline{z}}}

\newcommand{\half}{\ensuremath{\frac{1}{2}}}
\newcommand{\inner}[2]{ \langle {#1}, {#2} \rangle}
\newcommand{\znor}{Z_{\text{nor}}}
\newcommand{\zbnor}{\bar{Z}_{\text{nor}}}
\newcommand{\gnor}{\Gamma_{\text{nor}}}

\theoremstyle{definition}

\newtheorem{thm}{Theorem}[section]
\newtheorem{prop}[thm]{Proposition}
\newtheorem{lem}[thm]{Lemma}
\newtheorem{cor}[thm]{Corollary}
\newtheorem{rem}[thm]{Remark}
\newtheorem{ex}[thm]{Example}
\newtheorem{defi}[thm]{Definition}

\newtheorem{ass}[thm]{Assumption}

\DeclareMathOperator{\im}{Im} 
\DeclareMathOperator{\re}{Re}
 
\DeclareMathOperator{\Arg}{Arg}

\DeclareMathOperator{\Crit}{Crit}
\DeclareMathOperator{\Log}{Log}

\DeclareMathOperator{\id}{Id}
\DeclareMathOperator{\Gl}{Gl}
\DeclareMathOperator{\spn}{span}
\DeclareMathOperator{\Hom}{Hom}
\DeclareMathOperator{\I}{Id}
\DeclareMathOperator{\inv}{inv}
\DeclareMathOperator{\aff}{Aff}

\newboolean{mynotes}\setboolean{mynotes}{false}

\newcommand{\mycomments}[1]{
           \ifthenelse{\boolean{mynotes}}
                      {#1}{}
           }
	   
\newcommand{\marginlabel}[1]
{\mbox{}\marginpar[\raggedleft $\longrightarrow$ \\ \tiny\sf #1]{\raggedright $\longleftarrow$\\ \tiny\sf #1}}

\newcommand{\todo}[1]{\mycomments{\marginlabel{\tiny\sf{#1}}}}

\addtocounter{secnumdepth}{-2}

\author{R. Casta\~no-Bernard and D. Matessi}
\title{Lagrangian $3$-torus fibrations}

\begin{document}
\maketitle

\begin{abstract}
We give a method to construct singular Lagrangian 3-torus fibrations over certain a priori given integral affine manifolds with singularities, which we call \textit{simple}. The main result of this article is the proof that M. Gross' topological Calabi-Yau compactifications \cite{TMS} can be made into symplectic compactifications. As an example, we obtain a pair of compact symplectic 6-manifolds together with Lagrangian fibrations whose underlying affine structures are dual. The symplectic manifolds obtained are homeomorphic to a smooth quintic Calabi-Yau 3-fold and its mirror.
\end{abstract}

\tableofcontents

\section{Introduction.}

A map $f:X\rightarrow B$ from a smooth symplectic manifold onto a smooth manifold is a Lagrangian fibration if the regular locus of fibres has half the dimension of $X$ and the symplectic form restricts to zero there. The fibration is allowed to have singular fibres. In fact, interesting examples should in general include singular fibres. If the fibration map is smooth and proper, it is a well-known fact that the non-singular fibres are tori. Furthermore, away from the discriminant locus parametrizing the singular fibres, the base has the structure of an integral affine manifold. In other words, $B$ has an atlas whose change of coordinates are integral affine linear transformations.

\medskip 
Lagrangian fibrations lie at the crossroads of integrable systems,  toric symplectic geometry and more recently, Mirror Symmetry. For all three subjects, important issues are: the global topology of the fibration, the singularities of the fibres, the regularity of the fibration map and the affine structures induced on the base. In the recent years, integral affine geometry started to play a remarkably important role in Mirror Symmetry. The first evidence of this is given by Hitchin \cite{Hitchin}, who  observed that the SYZ duality \cite{SYZ} can be interpreted as a Legendre transform between integral affine manifolds. Later, Kontsevich and Soibelman \cite{Kontsevich-Soibelman} and Gross and Wilson \cite{G-Wilson2} proposed a landmark conjecture which, roughly speaking, says:

\begin{itemize}
\item[(1)] Degenerating families of Calabi-Yau manifolds approaching large complex structure limits should collapse down to a singular integral affine $S^n$.
\item[(2)] Mirror families should be (re)constructed starting from the affine manifolds in (1).  
\end{itemize} 
 
The first part of this conjecture is referred to as the Gromov-Hausdorff collapse, while the second part is usually called the reconstruction problem \cite{Gross_SYZ_rev}. We know that the Gromov-Hausdorff collapse does happen in dimension two \cite{G-Wilson2}. More recently, Gross and Siebert \cite{G-Siebert, G-Siebert2003} develop a program to reconstruct the ``complex side'' of the mirror using Logarithmic geometry. Kontsevich and Soibelman \cite{K-S-archim} approach the complex reconstruction problem using non-Archimedean analytic spaces. The final explanation of Mirror Symmetry is likely to emerge from the work deriving from these two main streams.

\medskip
On the ``symplectic side'' of the mirror, there is an analogous reconstruction problem. This paper is motivated by the following question. Can we construct symplectic manifolds starting from integral affine manifolds with singularities and obtain total spaces homeomorphic to mirror pairs of Calabi-Yau manifolds?

\medskip
To answer this question we take Gross' Topological Mirror Symmetry \cite{TMS} as a starting point. Gross developed a method to construct topological $T^3$ fibrations of 6-manifolds. This method consists, roughly, on the compactification of certain $T^3$ bundles by means of gluing suitable singular fibres. The discriminant locus in this case is a 3-valent graph with vertices labeled positive or negative. There are three types of singular fibres: \emph{generic fibres}, \emph{positive fibres} and \emph{negative fibres}, mapping to either points on the edges, or positive or negative vertices of the graph, respectively. The names are given according to the Euler characteristic of the fibres which can be 0, $+1$ or $-1$ respectively\footnote{Gross uses a different convention: $(2,2)$, $(1,2)$ and $(2,1)$, for generic, positive and negative fibres, respectively}. Gross' compactification produces a class of fibrations that can be dualized. As an example of this construction, Gross obtained a pair of smooth manifolds with dual topological $T^3$ fibrations, the first one being homeomorphic to the quintic 3-fold and the second one homeomorphic to a mirror of the quintic. 

\medskip
The main result of this paper is the proof that a compactification similar to that of Gross can be carried out in the symplectic category. The basic idea is the following. We start with an integral affine manifold with singularities $(B,\Delta,\mathscr A)$ with 3-valent graph singular locus $\Delta$. The affine structure on $B_0=B-\Delta$ induces a family of maximal lattices $\Lambda\subseteq T^\ast B_0$, together with a symplectic manifold $X(B_0)$ and an exact sequence

\[
0\rightarrow\Lambda\rightarrow T^\ast B_0\rightarrow X(B_0)\rightarrow 0.
\]

This gives us a Lagrangian $T^n$ bundle $f_0:X(B_0)\rightarrow B_0$. When $\mathscr A$ is \emph{simple} (cf. Definition \ref{def:simple}), $X(B_0)$ can be compactified to a topological 6-manifold $X(B)$ using Gross method. To define a symplectic structure on $X(B)$, in other words, to achieve a symplectic compactification of $X(B_0)$, one needs \emph{Lagrangian} models of generic, positive and negative singular fibres. The first two models have already been studied by the first author \cite{RCB1}. The construction of a Lagrangian negative model is much more delicate. An important part of this article is devoted to the construction of Lagrangian fibrations of negative type.

\medskip
While the generic and positive models are given by smooth maps and have codimension two discriminant loci, our model for the negative fibration is piecewise smooth and has mixed codimension one and two discriminant: it is an ``amoeba'' whose three legs are pinched down to codimension 2 (cf. Figure \ref{interpol}). In fact it can be described as a perturbation of Gross' negative fibration, localized in a small neighborhood of the `figure eight' (i.e. the singular locus of the negative fibre), which forces the singularities of the fibres to become isolated points and the discriminant locus to jump
to codimension one near the vertex. The topology of the total space is unchanged by this perturbation.
Joyce \cite{Joyce-SYZ}  had already conjectured that special Lagrangian fibrations should be in general
piecewise smooth and should have codimension 1 discriminant locus. 
Over the codimension 1 part of the discriminant locus, our model has exactly the topology which Joyce 
proposed as the special Lagrangian version of Gross' negative fibre.

\medskip

Our first attempt to construct a model of a Lagrangian negative fibration produces a fibration which fails to be smooth along a large codimension one subset, a whole plane containing the discriminant locus (cf. Example \ref{thin leg}). This model is not suitable for the symplectic compactification. This is essentially due to the fact that piecewise smooth fibrations in general do not induce integral affine structures on the base. The affine structure induced by fibrations of this sort consists of two pieces separated by the codimension one wall. Piecewise smooth fibrations of this type are called \textit{stitched} and have been studied in great detail by the authors \cite{CB-M-stitched, CB-M-torino}. It turns out that the information on the lack of regularity of these fibrations can be encoded into certain invariants. This allows us to have good control on the regularity of stitched fibrations. In particular, we are able to modify Example \ref{thin leg} to a Lagrangian fibration which induces an integral affine structure on the complement of a closed 2-disc containing the codimension one component of the discriminant. Moreover, away from this `bad disc', where the fibration fails to be smooth, the induced integral affine structure is simple.

\medskip
\todo{i removed some info, which was not essential}
Given a simple integral affine $3$-manifold with singularities $(B, \Delta, \mathscr A)$
a \textit{localized thickening} of $\Delta$ is given by the data $(\Delta^{\blacklozenge}, \{ D_{p^-} \}_{p^- \in \mathcal N})$ where: 
\begin{itemize}
\item[(i)] $\Delta^{\blacklozenge}$ is the closed subset obtained from $\Delta$ after replacing a 
           neighborhood of each negative vertex with a shape of the type depicted in 
           Figure~\ref{local:thicken} (an ``amoeba'' with thin legs).
\item[(ii)] $\mathcal N$ is the set of negative vertices and for each $p^-\in \mathcal N$, $D_{p^-}$ is a  disk containing the codimension $1$ component of $\Delta^{\blacklozenge}$ around  $p^-$ (depicted as the gray area in Figure~\ref{local:thicken}).
\end{itemize}
Given a localized thickening define
\[
B_{\blacklozenge} = B - \left( \Delta \cup \bigcup_{p^- \in \mathcal N} D_{p^-} \right).
\]
and denote by $\mathscr A_{\blacklozenge}$ the restriction of the affine structure on $B_{\blacklozenge}$

\medskip
The main result of this paper is the following (cf. Theorem \ref{the_symplectic_comp}):

\medskip
\noindent\textbf{Theorem.}
Given a compact simple integral affine $3$-manifold with singularities $(B, \Delta, \mathscr A)$,
all of whose negative vertices are straight. There is a localized thickening $(\Delta_{\blacklozenge}, \{ D_{p^-} \}_{p^- \in \mathcal N})$  and a smooth, 
compact symplectic $6$-manifold $(X, \omega)$ together with a piecewise smooth Lagrangian fibration 
$f: X \rightarrow B$ such that 
\begin{itemize}
\item[(i)] $f$ is smooth except along $\bigcup_{p^- \in \mathcal N} \,  f^{-1}(D_{p^{-}})$;
\item[(ii)] the discriminant locus of $f$ is $\Delta_{\blacklozenge}$;
\item[(iii)] there is a commuting diagram
\begin{equation*} 
\begin{CD}
X(B_{\blacklozenge}, \mathscr A_{\blacklozenge})  @>\Psi>> X\\
@Vf_0VV  @VVfV\\
B_{\blacklozenge} @>\iota>> B
\end{CD}
\end{equation*}
where $\psi$ is a symplectomorphism and $\iota$ the inclusion;
\item[(iv)] over a neighborhood of a positive vertex of $\Delta_{\blacklozenge}$ the fibration is 
             positive, over a neighborhood of a point on an edge the fibration is generic-singular, over
             a neighborhood of $D_{p^-}$ the fibration is Lagrangian negative. 
\end{itemize}

\medskip
As a corollary of Theorem \ref{the_symplectic_comp} and Gross' topological compactification \cite{TMS}, when $(B, \Delta, \mathscr A)$ is as in Example \ref{quint:aff}, the symplectic manifold obtained is homeomorphic to the quintic Calabi-Yau 3-fold. Applying the Legendre transform to Example \ref{quint:aff} produces a compact simple integral affine manifold with singularities $(\check B,\check{\Delta},\check{\mathscr A})$ \cite{G-Siebert}. The latter induces a bundle $X(\check B_0)$, dual to $X(B_0)$. By applying the Theorem we  obtain a compact symplectic manifold $(\check X, \check\omega)$ homeomorphic to Gross' topological compactification $X(\check B_0)$, therefore homeomorphic to a mirror of the quintic.

\medskip
The affine structures we consider here satisfy a property called \emph{simplicity}. Essentially, our notion of simplicity coincides with Gross and Siebert's simplicity in dimensions $n=2$ and $3$. Theorem \ref{the_symplectic_comp} should produce pairs of compact symplectic manifolds fibering over Gross and Siebert's integral affine manifolds, therefore producing a vast number of examples of dual Lagrangian $T^3$ fibrations.  For example, in \cite{Gross_Batirev}, Gross shows that 
to the pairs of Calabi-Yau's constructed with the method of Batyrev and Borisov as complete 
intersections in dual Fano toric varieties, one can associate a pair of simple affine manifolds with 
singularities which, when compactified, give back a pair of manifolds homeomorphic to the two Calabi-Yau's. The latter statement is the content of 
\cite{Gross_Batirev}{Theorem 0.1}, which is proved in \cite{GroSie_Tor} by Gross and Siebert. 
Combining this with our result, we obtain
a construction of symplectic manifolds fibred by Lagrangian tori, which are homeomorphic to the Batyrev 
and Borisov mirror pairs of 
Calabi-Yau manifolds.  Also, another source of examples may come from the structures constructed 
in \cite{Haase-Zharkov, Haase-ZharkovII, Haase-ZharkovIII}, provided they are simple.

\medskip
We should mention at this point that Lagrangian $T^3$ fibrations of Calabi-Yau manifolds have been constructed before by Ruan \cite{Ruan2, Ruan3, Ruan-III}. Ruan's construction does not use integral affine geometry, rather, it depends on a gradient flow argument. In particular Ruan's construction depends 
on the embedding inside an ambient manifold.  We suspect that Ruan's fibrations share many similarities with our symplectic compactifications but we haven't been able to verify this. It is not clear what kind of regularity Ruan's fibrations have, therefore whether they induce integral affine structures on the base.  One interesting aspect of our method is that it makes explicit connection with the formulation of Mirror Symmetry in \cite{Kontsevich-Soibelman} and \cite{G-Siebert2003}, where affine geometry is essential.

\medskip
The main motivation of this paper is Mirror Symmetry but we expect interesting applications in symplectic topology to emerge from the results we present here. Our construction of Lagrangian fibrations has a flavor similar to the work on almost toric symplectic geometry of Leung and Symington \cite{LeungSym}. A theory on almost toric 6-folds could emerge from the methods applied in this article. On the other hand, being our construction so explicitly connected to affine geometry, it is possible that the construction in 
Theorem \ref{the_symplectic_comp} will eventually shed light onto the new methods in symplectic 
enumerative problems arising from tropical geometry. \todo{check if you still like this sentence.}

\medskip
The material of this paper is organized as follows. We start giving in \S \ref{sect. topology} the description of Gross' compactification of topological $T^n$ bundles with semi-stable monodromy. Here we explain how to modify Gross' negative fibration to a fibration with a localized thickening near the negative vertex. In \S \ref{section: aff mfld & lag fib} we introduce the integral affine manifolds we use in the rest of the paper. We formalize our notion of simplicity by means of standard models of affine manifolds with singularities with prescribed holonomy. Our notion of simplicity coincides with the one in \cite{G-Siebert} in dimension $n=2$ and $3$. Simplicity is, essentially, 
a condition which guarantees that the induced Lagrangian $T^n$ bundles have semi-stable monodromy that can be compactified. We describe some examples of non-compact and compact simple integral affine manifolds with singularities. As an illustration of some of the methods we use, we show in Theorem \ref{thm. the K3} how, in dimension $n=2$, one can produce symplectic manifolds diffeomorphic to K3 surfaces. In \S \ref{section: aff mfld & lag fib} we describe Lagrangian models of positive and generic fibrations and prove that they induce integral affine structures which are simple. These models can be used to produce semi-stable symplectic compactifications over simple affine manifolds \emph{without} negative vertices (cf. Theorem \ref{thm:big_positive}). This is not enough, in general, to construct symplectic manifolds homeomorphic to Calabi-Yaus --such as a quintic and its mirror-- as one should normally include negative vertices. In any event, given the existence of simple affine bases with positive vertices only (or without any vertices at all) Theorem \ref{thm:big_positive} tells us how to construct a symplectic manifold together with a Lagrangian fibration over it. In this case, the Lagrangian fibrations obtained are everywhere smooth and the thickening of the discriminant is not necessary. There are explicit Examples of integral affine manifolds structures with no vertices \cite{Gross_Batirev} and Theorem \ref{thm:big_positive} can be used to produce symplectic compactifications. In \S \ref{sec:pwfibr} we move on to piecewise smooth fibrations. We give concrete examples of piecewise smooth Lagrangian $T^3$ fibrations. In particular, in Example \ref{thin leg} we explicitly construct a Lagrangian version of the topological negative fibration with fat discriminant given in \S \ref{sect. topology}. This model is piecewise smooth over a large region. In \S \ref{stitched fibr} we review some of the techniques we developed in \cite{CB-M-stitched}, which allow us to make certain non-smooth Lagrangian fibrations into smoother ones, such as the one in Example \ref{thin leg}.  The material of this section is rather technical and the reader may skip it in a first reading. In \S \ref{negative} we construct Lagrangian fibrations of negative type. These are local models whose discriminant is a localized thickening of a 3-valent negative vertex $p-$.  The fibration is smooth away from a 2-disc $D_{p-}$ containing the codimension 1 component of the discriminant. Away from $D_{p-}$, the affine structure is integral and simple. Finally, in \S \ref{section:compact} we prove Theorem \ref{the_symplectic_comp}.

\medskip
\textit{Aknowledgments:} The second author was partially funded by an EPSRC 
Research Grant GR/R44041/01, UK. Both authors would like to thank Mark Gross and Richard Thomas for 
useful discussions. Moreover they thank the following institutions for hosting them while working on 
this project: the Abdus Salam 
International Centre for Theoretical Physics in Trieste, the Department of Mathematics at Imperial College
in London, the Max Plank Institute in Leipzig and in Bonn, the Dipartimento di Scienze e Tecnologie Avanzate 
of the University of Piemonte Orientale in Alessandria (Italy), the Department of Mathematics of the 
University of Pavia (Italy), the IHES in Paris. 

\section{The topology.}\label{sect. topology}
In this section we review Mark Gross' Topological Mirror Symmetry \cite{TMS}, which is the starting point for the results of this paper. Gross developed a method to compactify certain $T^n$ bundles over $n$-dimensional manifolds to obtain topological models of Calabi-Yau manifolds. We now outline how this method works. Along the way, we discuss how Gross' method can be modified to produce topological fibrations with mixed codimension one and two discriminant locus. We focus in dimension $n=2$ and $3$.

\medskip
A \textit{topological $T^n$ fibration} $f:X\rightarrow B$ is 
a continuous, proper, surjective map between smooth manifolds, $\dim X=2n$, 
$\dim B=n$, such that for a dense open set $B_0\subseteq B$ and for all $b\in B_0$ the fibre $X_b=f^{-1}(b)$ is homeomorphic to an $n$-torus. We call the set $\Delta:=B-B_0$ the \textit{discriminant locus of $f$}. 
Sometimes we will denote a topological fibration by a triple $\mathcal F = (X, f, B)$. Notice that this notion of fibration differs from the usual differential geometric one in the sense that here $\mathcal F$ is allowed to have singular fibres over points in $\Delta$. Allowing singular fibres is necessary if we aim at obtaining total spaces with interesting topology, such as Calabi-Yau manifolds other than complex tori.
When $X$ is a symplectic manifold, with symplectic form $\omega$, 
a topological $T^n$-fibration is said to be Lagrangian if $\omega$ restricted to the smooth part 
of every fibre vanishes.

\begin{defi}\label{def. symp eq}
Let $\mathcal F=(X, f,B)$ and $\mathcal F'=(X',f',B')$ be a pair of 
topological fibrations with discriminant loci $\Delta$ and $\Delta'$ 
respectively. We define the following notions of \textit{conjugacy} between $\mathcal F$ and $\mathcal F'$:
\begin{itemize}
\item[(i)] We say that $\mathcal F$ is \textit{conjugate} to $\mathcal F'$ 
if there exist a homeomorphism $\psi :X \rightarrow X'$ and a homeomorphism 
$\phi :B\rightarrow B'$ sending $\Delta$ to $\Delta'$ homeomorphically, such 
that $f'\circ\psi=\phi\circ f$. We shall say that $\mathcal F$ is
$(\psi ,\phi)$\textit{-conjugate} to $\mathcal F'$ whenever the
specification is required. 
\item[(ii)] If in addition $X$ and $X^{\prime}$ are symplectic manifolds and 
the fibrations are Lagrangian, we will say that 
$\mathcal F$ is \textit{symplectically conjugate} to $\mathcal F'$
if $\psi$ is a $C^\infty$ symplectomorphism and $\phi$ is a $C^\infty$ diffeomorphism. 
\item[(iii)] Given points $b \in \Delta$ and 
$b^{\prime} \in \Delta^{\prime}$, we shall say that $\mathcal F$ is 
\textit{(symplectically) conjugate to} $\mathcal F^{\prime}$ \textit{over} $\Delta$ 
(or over $b$ and $b^{\prime}$) if there are neighborhoods $U$ and $U^{\prime}$
of $\Delta$ and $\Delta^{\prime}$ (or of $b$ and $b^{\prime}$) 
respectively, such that $(f^{-1}(U), f,U)$ is (symplectically) conjugate to
$((f^{\prime})^{-1}( U^{\prime}), f^{\prime}, U^{\prime})$.
\end{itemize}
\end{defi}
Part (iii) can  also be found in the literature as \textit{semi-global (symplectic) equivalence} as it involves a fibred neighborhood of a fibre but not the total space. When $\mathcal F$ carries additional specified data, --e.g. a (Lagrangian) section  or a choice of basis of $H_1(X,\Z)$-- one may also consider a slightly stronger version of (i)-(iii) which requires that the specified data is preserved, e.g. that $\phi$ sends the section of $f$ to the section of $f'$ and a basis of $H_1(X,\Z)$ to a basis of $H_1(X',\Z)$ . Clearly all three notions define equivalence relations. The corresponding equivalence classes will be called \textit{germs} of fibrations. Throughout this article 
we will often use conjugation to topologically or symplectically glue together fibred sets in order to obtain larger fibred sets and eventually produce compact (symplectic) manifolds. 

Given a topological (or Lagrangian) fibration $\mathcal F = (X, f, B)$ and a subset $U \subset B$, we will often use the notation $\mathcal F|_{U}$ to denote the fibration $(f^{-1}(U), f, U)$ and 
we will refer to it as the restriction of $\mathcal F$ to $U$.

\medskip
The topological fibrations considered by Gross have everywhere codimension two discriminant. For $n=2$, $\Delta$ is a finite collection of points and the singular fibres are \textit{nodal}. For $n=3$, $\Delta$ is a connected trivalent graph with vertices labeled `positive' or `negative'. There are three types of singular fibres in this case: \textit{generic-singular fibres}, i.e. the product of a nodal fibre with $S^1$; \textit{positive fibres}, i.e. a 3-torus with a 2-cycle collapsed to a point; and \textit{negative fibres}, singular along a `figure eight'. For a more detailed description of these singular fibres we refer the reader to Examples \ref{ex top ff}, \ref{ex. (2,2)}, \ref{ex. (2,1)} and \ref{ex. (1,2)} below or to \cite{TMS} for further details.

\medskip
In this article, we will allow $\Delta$ to jump dimension, i.e. $\Delta$ will include the region $\Delta_a\subseteq\Delta$, which may be regarded as a ``fattening" of a graph near negative vertices. We also propose a new model with discriminant locus of type $\Delta_a$ (cf. Example \ref{ex. alt (2,1)}) which is an alternative to Gross' negative fibration and, in some sense, it is a more generic version of it. The idea of using models with codimension one discriminant was first suggested by Joyce \cite{Joyce-SYZ}{\S 8}, based on his knowledge of special Lagrangian singularities. Ruan's Lagrangian fibrations \cite{Ruan, Ruan2, Ruan3, Ruan-III} also have codimension one discriminant loci.

\medskip
Consider the following three closed subsets of $\R^3$: 
\[ \mathscr C_{e} = \{ x_1 = x_2 = 0 \}, \]
\[ \mathscr C_{d} = \{ x_1 = x_2 = 0, \, x_3 \leq 0 \} \cup \{ x_1 = x_3 = 0, \, x_2 \leq 0 \} 
                                       \cup \{ x_1=0, \, x_2 = x_3 \geq 0 \}, \]
\[ \mathscr C_{a} = \mathscr C_{d} \cup \left\{ x_1 = 0, \, x_2^2 + x_3^2 \leq \half \right\}. \]
Clearly $\mathscr C_{d}$ is a model of a neighborhood of a vertex in a three valent graph and 
$\mathscr C_{a}$ can be regarded as a fattening of $\mathscr C_{d}$ around the vertex. We also 
denote by $D^3$ the open unit ball in $\R^3$.

\begin{figure}[!ht]
\begin{center}
\input{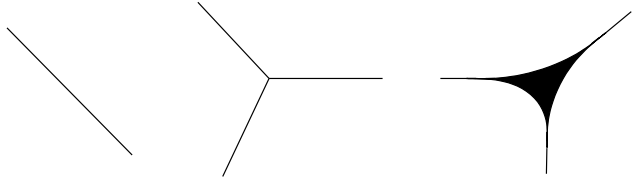}
\caption{The three possibilities for $U_i \cap \Delta$, $n=3$.}\label{deltas}
\end{center}
\end{figure}

In this paper, we consider fibrations satisfying the following topological properties:

\begin{ass} \label{top ass} \todo{I think now this is better}
Let $\mathcal F = (X, f, B)$ be a topological $T^n$ fibration with discriminant locus $\Delta\subseteq B$ and fibre $X_b$ over $b\in B$. We assume that $\mathcal F$ satisfies the following conditions:

\begin{enumerate}
\item for $n=2$, $\Delta$ is a finite union of points and given a 
small neighborhood $U$ of a point in $\Delta$, the fibration $\mathcal F|_U$ is topologically 
conjugate to a nodal fibration (see Example~\ref{ex top ff});

\item for $n=3$, there is a finite covering $\{ U_i \}$ of $\Delta$ with open subsets of 
      $B$ such that one of the following three possibilities occur (see also Figure~\ref{deltas}): 

\begin{enumerate}
\item the pair $(U_i, U_i \cap \Delta)$ is homeomorphic to $(D^3, D^3 \cap \mathscr C_{d})$ and 
      $\mathcal F|_{U_i}$ is topologically conjugate to either a positive or a negative fibration 
      (see Examples~\ref{ex. (1,2)} and \ref{ex. (2,1)});
\item the pair $(U_i, U_i \cap \Delta)$ is homeomorphic to $(D^3, D^3 \cap \mathscr C_{a})$ 
      and $\mathcal F|_{U_i}$ is topologically conjugate to an alternative negative fibration 
      (see Example~\ref{ex. alt (2,1)});
\item the pair $(U_i, U_i \cap \Delta)$ is homeomorphic to $(D^3, D^3 \cap \mathscr C_{e})$ 
      and $\mathcal F|_{U_i}$ is topologically conjugate to a generic-singular fibration 
      (see Example~\ref{ex. (2,2)});
\end{enumerate} 
We denote by $\Delta_d$ the set of points in $\Delta$ belonging to a $U_i$ satisfying $(a)$, 
which are the vertices of $U_i \cap \Delta$. We call these points vertices of $\Delta$. 
We denote by $\Delta_a$ the union of the sets $U_i \cap \Delta$, where $U_i$ satisfies $(b)$; we 
can assume these sets to be pairwise disjoint.  A point in $\Delta$ admitting open neighborhood $U$ of 
$B$ such that $(U, U \cap \Delta)$ is homeomorphic to $(D^3, D^3 \cap \mathscr C_{e})$ is called an 
\textit{edge point}. We denote by $\Delta_g$ the set of edge points.
\end{enumerate}

\end{ass}

We denote by $\Sigma$ the locus formed by the singularities of all the fibres, therefore sometimes $\Sigma$ will also be denoted by $\Crit(f)$; when $f$ is smooth,  $\Crit (f)$ will indeed coincide with the set of critical points of $f$. We insist, however, that $f$ is not a priori required to be a smooth map. In fact,  we will see that, near $\Delta_a$, our fibrations are not smooth. Inspired by tropical geometry, we refer to a connected component of $\Delta_a$ as a \textit{3-legged amoeba} (with thin ends). As we will see later when we will 
introduce affine structures, an important property of $\Delta_a$ is that it is locally planar, i.e. each 
connected component of $\Delta_a$ is contained, in some sense, in a 2-plane.

\begin{defi} \todo{also this has been fixed}
Let $f:X\rightarrow B$ be a topological $T^n$ fibration and let $U\subset B$ be an open contractible 
neighborhood of $b\in\Delta$ such that $U\cap\Delta = \{ b \}$, when $n=2$; or else, when $n=3$, such that $U$ satisfies $(a)$, $(b)$ or $(c)$ in point $2$ of Assumption~\ref{top ass}. 
Let $X_{b_0}$ be a fibre over $b_0\in U-\Delta$. Consider the monodromy representation
\[
\mathcal{M}_b:\pi_1(U-\Delta, b_0)\rightarrow SL(H_1(X_{b_0},\Z)).
\] 
The image of $\mathcal{M}_b$ is called the \textit{local monodromy group about $X_b$} (also denoted by $\mathcal M_b$). 
\end{defi}

\medskip
Now we review the local models of these fibrations. For the details we refer the reader to \cite{TMS}{\S 2}. 
The construction of the local models relies on the following:

\begin{prop}\label{prop. S1-bundle}
Let $Y$ be a manifold of dimension $2n-1$. Let $\Sigma\subseteq Y$ be an oriented submanifold of codimension three and let $Y'=Y-\Sigma$. Let $\pi':X'\rightarrow Y'$ be a principal $S^1$-bundle over $Y'$ with Chern class $c_1=\pm 1$. For each triple $(Y,\Sigma,\pi')$ there is a unique compactification $X=X'\cup \Sigma$ extending the topology of $X'$, making $X$ into a  manifold and such that
\[\begin{array}{ccc}
X' & \hookrightarrow & X \\ 
\downarrow & \  & \downarrow \\ 
Y' & \hookrightarrow & Y
\end{array}\] commutes, with $\pi :X\rightarrow Y$ proper and $\pi|_\Sigma :\Sigma\rightarrow\Sigma$ the identity. 
\end{prop}

\begin{rem}\label{rem. pi local mod} One can explicitly describe the above compactification as follows. For any point $p\in\Sigma$ there is a neighborhood $U\subset Y$ of $p$ such that $U\cong\R^3\times\C^{n-2}$ and $U\cap\Sigma$ can be identified with $\{0\}\times\C^{n-2}$. By unicity of $\pi$, there is a commutative diagram
\begin{equation}
\begin{CD}
\pi^{-1}(U) @>\cong >> \C^2\times\C^{n-2} \\ 
@V\pi VV    @V\bar\pi VV \\ 
U @>\cong >> \R^3\times\C^{n-2}
\end{CD}
\end{equation}
where $\bar\pi (z_1,z_2,\zeta )=(|z_1|^2-|z_2|^2,z_1z_2,\zeta)$, $\zeta\in\C^{n-2}$. 
\end{rem}

\begin{figure}[!ht]
\begin{center}
\input{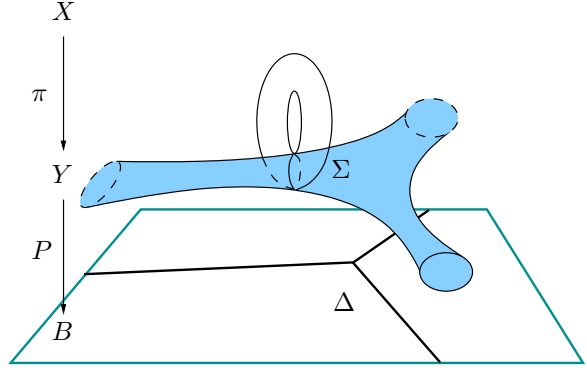}
\caption{Negative fibration.}\label{pants_thin}
\end{center}
\end{figure}

The constructions of topological $T^n$ fibrations in this section are based on the following 
basic principle. One starts with a manifold $Y=B\times T^{n-1}$ with $\dim B =n$, a submanifold $\Sigma\subset Y$ and a map $\pi:X\rightarrow Y$ as in Proposition \ref{prop. S1-bundle}. The trivial $T^{n-1}$ fibration $P:Y\rightarrow B$ can be lifted to a $T^n$ fibration $f:=P\circ\pi :X\rightarrow B$ with discriminant locus $\Delta:=P(\Sigma )$. One can readily see that for $b\in\Delta$, the singularities of the fibre $X_b$ occur along $\Sigma\cap P^{-1}(b)$. The set $\Sigma$ --which is the locus of singular fibres of $\pi$-- can be regarded as the locus where the vanishing cycles of the fibres of $f$ collapse (cf. Figure \ref{pants_thin}).

\begin{ex}[Nodal fibration] \label{ex top ff} 
This example is the topological model for the fibration over a point of
$\Delta$ in the case $n=2$. Let $D$ be the unit disc in $\C$ and  $D^\ast=D-\{0\}$. Let 
$f_0:X_0\rightarrow D^\ast$ be a $T^2$-bundle 
with monodromy generated by $\left(\begin{array}{cc}
1 & 0 \\ 
1 & 1
\end{array}\right)$. We can use Proposition \ref{prop. S1-bundle} to compactify $X_0$ as follows. The monodromy invariant cycle, $L\in H_1(f_0^{-1}(b),\Z )$, induces a fibre preserving $T(L)$ action, with $T(L)=L\otimes\R\slash L$. The quotient modulo this action yields an $S^1$-bundle $\pi_0:X_0\rightarrow Y_0=D^\ast\times S^1$. One can verify that $\pi_0$ extends to an $S^1$-bundle $\pi':X'\rightarrow Y'=D\times S^1-\{(0,p)\}$, where $p\in S^1$. Furthermore $c_1(\pi')=\pm 1$. Then Proposition \ref{prop. S1-bundle} ensures that $X'$ compactifies to a manifold $X=X'\cup\{pt\}$ and that there is a proper map $\pi :X\rightarrow Y=D\times S^1$ extending $\pi'$. Defining $P:Y\rightarrow D$ as the projection map, we obtain a fibration $f=P\circ\pi :X\rightarrow D$ extending $f_0$. The only singular fibre, $f^{-1}(0)$, is homeomorphic to $T^2=S^1\times S^1$ after $S^1\times\{x\}\subset T^2$ is collapsed to $x$. We denote this fibre by $I_1$, following Kodaira's notation for singular fibres of elliptic fibrations. In Hamiltonian mechanics, a Lagrangian fibration with this topology is known as a \textit{focus-focus fibration.}
\end{ex}

\begin{ex}[Generic singular fibration]\label{ex. (2,2)}
This example is the model for the fibration over a neighborhood of an edge point of $\Delta$ --in \cite{TMS} this is called $(2,2)$ fibration. Let $B=D\times (0,1)$, where $D \subset \C$ is the unit disc, and let $Y=T^2\times B$. Define $\Sigma\subset Y$ to be the cylinder sitting above $\{0\}\times (0,1)\subset B$ defined as follows. Let $e_1,e_3$ be a basis of $H_1(T^2,\Z)$.
Let $S^1\subset T^2$ be a circle representing the homology class 
$e_3$. Define $\Sigma=S^1\times\{ 0\}\times (0,1)$. Now let $\pi':X'\rightarrow Y':=Y-\Sigma$ be an $S^1$-bundle with Chern class $c_1=1$. Then $X'$ compactifies to a manifold $X=X'\cup\Sigma$ and there is a proper map $\pi :X\rightarrow Y$ extending $\pi'$. We can now define $f=P\circ\pi :X\rightarrow B$ where $P:Y\rightarrow B$ is the projection. Then it is clear $f$ is a $T^3$ fibration with singular fibres homeomorphic to $I_1\times S^1$ lying over $\Delta:=\{ 0\}\times (0,1)$.
If $e_2$ is an orbit of $\pi$, one can take $e_1,e_2, e_3$ as a basis of $H_1(X_b, \Z)$, where $X_b$ is a regular fibre. In this basis, $e_2$ and $e_3$ are monodromy invariant and a generator of the monodromy group 
of $f$ about $\Delta$ is represented in this basis by
\begin{equation}\label{eq matrix g}
T=\left( \begin{array}{ccc}
                 1 & 0 & 0 \\
                 1 & 1  & 0 \\
                 0 & 0  & 1 \end{array} \right).
\end{equation}
\end{ex}

\begin{figure}[ht]
\psfrag{g1}{$g_{1}$}
\psfrag{g2}{$g_2$}
\psfrag{g3}{$g_3$}
\psfrag{b}{$b$}
\begin{center}
\epsfig{file=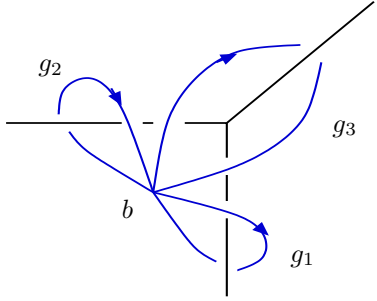, width=5cm,height=4cm}
\caption{Loops $g_1$, $g_2$ and $g_3$, such that $g_1 g_2 g_3 = 1$.}\label{loops}
\end{center}
\end{figure}

\begin{ex}[Negative fibration]\label{ex. (2,1)}
This example is one of the two models over a neighborhood of
a point in $\Delta_d$ --in \cite{TMS} this is called $(2,1)$ fibration. Let $Y=T^2\times B$ with $B$ homeomorphic to a 3-ball. 
Let $\Delta\subset B$ be a cone over three distinct, non-collinear points. We write $\Delta=\{b_0\}\cup\Delta_1\cup\Delta_2\cup\Delta_3$ where $b_0$ is the vertex of $\Delta$ and the $\Delta_i$ are the legs of $\Delta$. Fix a basis $e_2$, $e_3$ for $H_1(T^2,\Z)$. Define $\Sigma \subset T^2\times B$ to be a pair of pants lying over $\Delta$ such that for $i=1,2,3$, $\Sigma \cap (T^2\times\Delta_i)$  is a leg of $\Sigma $ which is the cylinder generated by \todo{to get monodromy convention uniform for us i changed 
the sign of $e_2$} $-e_3$, $-e_2$ and $e_2+e_3$ respectively. These legs are glued together along a nodal curve or `figure eight' lying over $b_0$. Now consider an $S^1$-bundle $\pi':X'\rightarrow Y'=Y-\Sigma $ with Chern class $c_1=1$. This bundle compactifies to $\pi:X\rightarrow Y$. Now consider the projection map $P:Y\rightarrow B$. The composition $f=P\circ\pi$ is a proper map. The generic fibre of $f$ is a 3-torus. For $b\in\Delta$ the fibre $f^{-1}(b)$ is singular along $P^{-1}(b)\cap \Sigma $, which is a circle when $b\in\Delta_i$, or the aforementioned figure eight when $b=b_0$. Thus the fibres over $\Delta_i$ are homeomorphic to $I_1\times S^1$, whereas the central fibre, $X_{b_0}$, is singular along a nodal curve.
A regular fibre can be regarded as the total space of an $S^1$-bundle over $P^{-1}(b)$. We can take as 
a basis of $H_1(X_b,\Z )$, $e_1(b), e_2(b),e_3(b)$, where $e_2$ and $e_3$ are the 1-cycles in $P^{-1}(b)=T^2$ as before and $e_1$ is a fibre of the $S^1$-bundle. The cycle $e_1(b)$ vanishes as $b\rightarrow\Delta$. In this basis, the matrices 
generating the monodromy group corresponding to loops $g_i$ about $\Delta_i$ with $g_1g_2g_3=1$, 
(cf. Figure \ref{loops}) are \todo{the matrices changed correspondingly}
\begin{equation}\label{eq matrix neg}
T_1=\left( \begin{array}{ccc}
                 1 & 1 & 0 \\
                 0 & 1  & 0 \\
                 0 & 0  & 1 
              \end{array} \right),\quad 
T_2=\left( \begin{array}{ccc}
                 1 & 0 & -1 \\
                 0 & 1  & 0 \\
                 0 & 0  & 1 
              \end{array} \right),\quad 
T_3=\left( \begin{array}{ccc}
                 1 & -1 & 1 \\
                 0 & 1  & 0 \\
                 0 & 0  & 1 
\end{array}\right). 
\end{equation}
\end{ex}

\begin{figure}[!ht] 
\begin{center}
\input{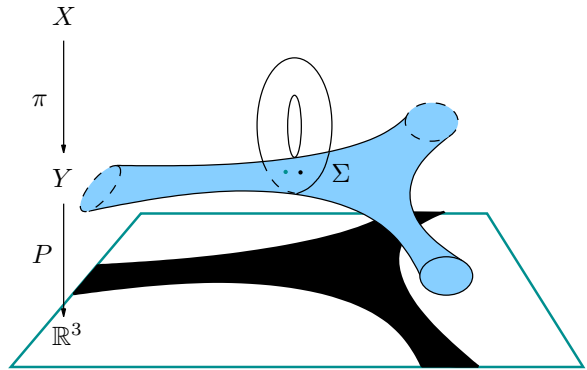}
\caption{Negative fibration with amoeba-like discriminant.}\label{pants}
\end{center}
\end{figure}

\begin{ex}[Alternative negative fibration]\label{ex. alt (2,1)}
This is the local model for a fibration over a neighborhood of a component of $\Delta_a$. Consider $Y$ and $\Sigma$ as in Example~\ref{ex. (2,1)}.
Now think of making a small perturbation of $\Sigma$ just in a neighborhood of the ``figure eight'' --i.e. where the three cylinders forming $\Sigma$ are joined together-- and leaving the rest unchanged. A generic perturbation will be such that, near the fibre over $b_0$, $\Sigma$ will intersect the fibres of $P: Y \rightarrow B$ in isolated points. Then $P(\Sigma)$ will have the shape of a $3$-legged amoeba. One then 
constructs the bundle $\pi':X'\rightarrow Y'=Y-\Sigma $ with Chern class $c_1=1$ and compactifies it to $\pi:X\rightarrow Y$. The total fibration is $f= P \circ \pi$. 

\medskip
We can give an explicit construction of a fibration of this type, following ideas in \cite{Gross_spLagEx}{\S 4}. 
Consider $(\C^{\ast})^{2}$ with the $T^2$ fibration 
$\Log: (v_1, v_2) \mapsto (\log |v_1|, \log |v_2| )$.
Let $Y = \R \times (\C^{\ast})^{2} $ and $P$ be the fibration
\[ P: (t, v) \rightarrow (t, \Log v ), \]
where $t \in \R$ and $v = (v_1, v_2) \in (\C^{\ast})^{2}$. 
Define a surface $\Sigma^{\prime}$ in $(\C^{\ast})^2$ to be
\[  \Sigma^{\prime} = \{ v_1 + v_2 + 1 = 0 \}, \]
and view it as a surface in $\{ 0 \} \times (\C^{\ast})^2 \subset Y$.
Clearly $P(\Sigma^{\prime})$ is $\{ 0 \} \times \Log( \Sigma^{\prime})$
and one can compute that it has the shape depicted in Figure~\ref{fig: amoeba}. 
Images by $\Log$ of algebraic curves in $(\C^{\ast})^{2}$ are known in the 
literature as amoebas, and this explains the name we gave to the components of 
$\Delta_a$.

\begin{figure}[!ht]
	\centering
	\includegraphics[width=4cm,height=4cm,bb=0 0 258 258]{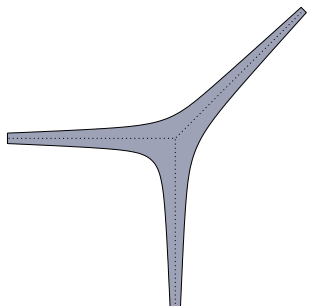}
	\caption{Amoeba of $v_1+v_2+1=0$}
	\label{fig: amoeba}
\end{figure}

As a surface in $\C^2$,  $\Sigma^{\prime}$ intersects 
$\{ v_1 = 0 \}$ in 
$q_1= (0, 0, -1)$ and $\{ v_2 = 0 \}$ in $q_2 = (0, -1,0)$. One can 
see that in a small neighborhood of $q_1$ one can twist $\Sigma^{\prime}$ 
slightly, so to make it coincide, in a smaller neighborhood, with 
$\{ v_2 = -1 \}$. 
Similarly one can twist $\Sigma^{\prime}$ near $q_2$, so to make it coincide with 
$\{ v_1 = -1 \}$. Finally, when $|v_1|$ and $|v_2|$ are both big, we can
twist $\Sigma^{\prime}$ so to coincide with $\{ v_1 + v_2 = 0 \}$. Let
$\Sigma$ be this new twisted version of $\Sigma^{\prime}$. A schematic description of these twistings is described in Figure \ref{interpol}, where $\Sigma'$ is the light-colored diagonal line and $\Sigma$ is the over-imposed twisted dark line. It is clear that $P(\Sigma) = \{ 0 \} \times \Log(\Sigma)$ will have the shape 
of a $3$-legged amoeba whose legs have been pinched to $1$-dimensional 
segments toward the ends, as depicted in the right-hand side of Figure \ref{interpol}  (Mikhalkin \cite{mikh_pants} also defines a similar construction and calls this shape a localized amoeba). The bundle $\pi':X'\rightarrow Y'=Y-\Sigma $ 
with Chern class $c_1=1$ and its compactification $\pi:X\rightarrow Y$ can
again be constructed. The fibration is $f = P \circ \pi$ and 
$\Delta =  P(\Sigma)$.

\begin{figure}[!ht]
\begin{center}
\input{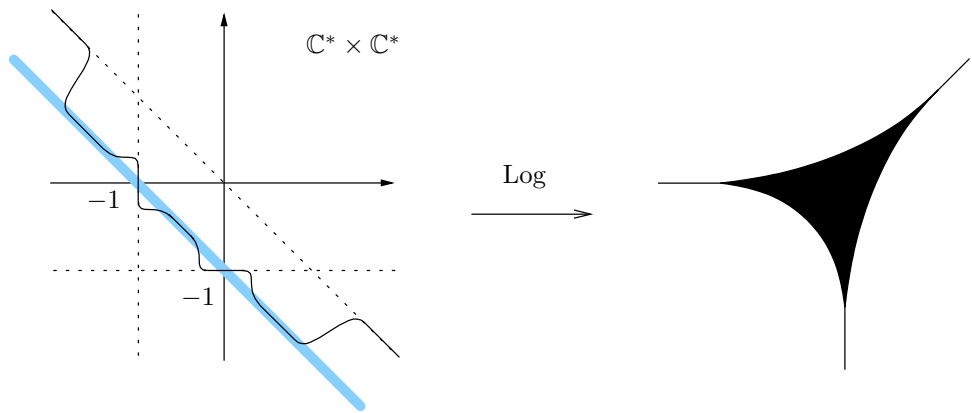}
\caption{The twisted $\Sigma$ gives and amoeba with thin legs.}\label{interpol}
\end{center}
\end{figure}

We give a description of the fibration over the codimension $1$ part of 
$\Delta$. One can see that the fibres of $\Log$ over a point 
in the interior of the amoeba intersect $\Sigma$ in two distinct points. 
These two points come together to a double point as the base point approaches
the boundary of the amoeba. 
If $p_1$ and $p_2$ are two points on $T^2$ --which may coincide-- then the 
singular fibres of $f$ look like $S^1 \times T^2$ after $S^1 \times \{ p_j \}$
is collapsed to a point. This behavior is topologically the same as the one conjectured by 
Joyce \cite{Joyce-SYZ} for special Lagrangian $T^3$ fibrations.
Moreover, the singularities of the fibres are modeled on those of an explicit example of 
a special Lagrangian fibration with non-compact fibres (cf. Joyce \cite{Joyce-SYZ}{\S 5}). 

\medskip
In view of Proposition \ref{prop. S1-bundle} and Remark \ref{rem. pi local mod}, 
the total space $X$ in this example is diffeomorphic to the one in 
Example \ref{ex. (2,1)}, although the fibrations differ. 
In both cases the singularities of the fibres occur along the 
intersection of the critical surface $\Sigma$ with the fibres of $P$. 
But the intersections happen in a different way. 
In Example \ref{ex. (2,1)} they occur either along circles, or along a figure 
eight. Here they occur along circles when the fibre is over a point in 
the codimension $2$ part of $\Delta$ and as isolated points when the fibre is 
over a point in the codimension $1$ part. 
As argued by Joyce, the isolated singularities are more generic in certain 
sense (cf. \cite{Joyce-SYZ}{\S 3}). A schematic description of the fibration over the codimension $1$ part 
of $\Delta$ is depicted in Figure \ref{pants}. It can be compared with 
Figure~\ref{pants_thin}. We remark that over the codimension $2$ part of 
$\Delta$, the fibration has the same topology of the generic 
singular fibration of Example~\ref{ex. (2,2)}. It follows that the monodromy 
around the legs is same as the monodromy of Example \ref{ex. (2,1)}, 
i.e. it is represented by the matrices (\ref{eq matrix neg}).
\end{ex}

Most of the effort in this paper is devoted to the construction of a fibration 
as in the previous example which is also \emph{Lagrangian} with respect to a symplectic
form on $X$. In the process we will also make more explicit the twistings
which allow us to deform $\Sigma^{\prime}$ into $\Sigma$. 

\medskip
Observe that in the above examples there is a fibre-preserving $S^1$-action, 
induced by the $S^1$ bundle $\pi'$. One can use the same principle to construct 
$T^2$-invariant fibrations starting from suitable compactifications of 
$T^2$-bundles:

\begin{ex}[Positive fibration]\label{ex. (1,2)}
This model is the other possible fibration over a neighborhood of a point in $\Delta_d$ --in \cite{TMS} this is called $(1,2)$ fibration.
Let $Y=S^1\times B$ with $B$ and $\Delta\subset B$ as in Example \ref{ex. (2,1)}. Let  $Y'=Y\setminus (\{ p\}\times\Delta )$, where $p\in S^1$. Let $L\cong\mathbb{Z}^2$ and define $T(L)=L\otimes_{\mathbb{Z}}\R\slash L$. Now consider a principal $T(L)$-bundle $\pi' :X'\rightarrow Y'$. Under some mild assumptions on $\pi'$ (cf. \cite{TMS} Prop. 2.9), there is a unique 
manifold $X$ with $X'\subset X$ extending the topology of $X'$ and a proper extension 
$\pi :X\rightarrow Y$ of $\pi'$. 
The composition of $\pi$ with the projection $Y\rightarrow B$ defines a topological $T^3$-fibration, $f:X\rightarrow B$. The fibre of $f$ over $b\in B\setminus\Delta$ is $T^3$. The fibre over $b\in\Delta_i$ is homeomorphic to $S^1\times I_1$, whereas the fibre over the vertex $b_0\in\Delta$ is homeomorphic to $S^1\times T^2\slash (\{ point\}\times T^2)$. It is proved in
\cite{TMS} that the monodromy group of this model is generated, in some basis, by the inverse transpose of the matrices (\ref{eq matrix neg}). The reader should not worry, at this point, for the lack of details in this 
description as we will give explicit Lagrangian models for this example later on. 
\end{ex}

Notice that the monodromy representation of the above models is \emph{semi-stable}, in other words the monodromy matrices of $\mathcal M_b$ are unipotent. This terminology is imported from the classical theory of elliptic fibrations.  The topological models described above may be regarded as 3-dimensional topological analogues to semi-stable singular elliptic fibres. We are now ready to state Gross' result. We refer the reader to \cite{TMS}{\S 2} for the details:

\begin{thm}[Gross]\label{thm. semi-stable}
Let $B$ be a 3-manifold and let $B_0\subseteq B$ be a dense open set such that $\Delta:=B-B_0$ is a trivalent graph, i.e. such that $\Delta = \Delta_d\cup\Delta_g$. Assume that the vertices of $\Delta$ are labeled, i.e. $\Delta_d$ decomposes as a union $\Delta_+\cup\Delta_-$ of positive and negative vertices. Suppose there is a $T^3$ bundle $f_0:X (B_0)\rightarrow B_0$ such that its local monodromy $\mathcal M_b$ is generated by 
\begin{enumerate}
\item $T$ as in (\ref{eq matrix g}), when $b\in\Delta_g$;
\item $T_1,T_2,T_3$ as in (\ref{eq matrix neg}), when $b\in\Delta_-$;
\item $(T_1^t)^{-1},(T_2^t)^{-1},(T_3^t)^{-1}$, when $b\in\Delta_+$.
\end{enumerate}
Then there is a $T^3$ fibration $f:X\rightarrow B$ and a commutative diagram
\[\begin{array}{ccc}
X(B_0) & \hookrightarrow & X \\ 
\downarrow & \  & \downarrow \\ 
B_0 & \hookrightarrow & B.
\end{array}\]
Over connected components of $\Delta_g$, $(X, f, B)$ 
is conjugate to the generic singular fibration, over points of
$\Delta_+$ it is conjugate to the positive fibration
and over points of $\Delta_-$ to the negative fibration.   
\end{thm}

A topological manifold $X$ obtained from $X(B_0)$ as in Theorem \ref{thm. semi-stable} is called a \textit{topological semi-stable compactification}. Fibrations arising from semi-stable compactifications satisfy the so-called \emph{topological simplicity} property (cf. \cite{TMS}\S 2). This is intimately related to the \emph{affine simplicity} of the subsequent sections. It is due to simplicity that Theorem \ref{thm. semi-stable} may be used to produce dual $T^n$ fibrations of manifolds homeomorphic to mirror pairs of Calabi-Yau manifolds. In \S \ref{section: aff mfld & lag fib} we shall review Gross' construction of a $T^3$ bundle $X(B_0)$ which compactifies to a smooth manifold $X$ 
homeomorphic to the quintic hypersurface in $\mathbb{P}^4$. The compactification 
of the dual bundle produces a manifold, $\check X$, which is homeomorphic to the 
mirror quintic. This construction gives evidence that the 
SYZ duality should indeed explain Mirror Symmetry.

\medskip
Theorem \ref{thm. semi-stable} could be stated and proved, with little effort, replacing $\Delta_-$ with $\Delta_a$, i.e. 
replacing a neighborhood of each negative vertex, with a $3$-legged amoeba. Over connected components of $\Delta_a$, the resulting fibration would then be conjugate to the alternative negative fibration of Example~\ref{ex. alt (2,1)} but the topology of the total space remains the same.
In fact we can do more: the main result of this paper is the proof that there exist \emph{symplectic} semi-stable compactifications with respect to which the fibres are 
Lagrangian. The starting point for this compactifications will be the Lagrangian $T^3$ bundles obtained
from affine $3$-dimensional manifolds.

\section{Affine manifolds and Lagrangian fibrations} \label{section: aff mfld & lag fib}
Let us denote by $\aff (\R)=\R^{n} \rtimes \Gl(n,\R)$ the group of \textit{affine linear} transformations, i.e. elements in $\aff (\R)$ are maps $A:\R^n\rightarrow\R^n$, $A(x)=L(x)+v$, where $L\in\Gl(n,\R)$ and $v\in\R^n$.
The subgroup of $\aff (\R)$ consisting of affine linear transformations with integral linear part will be denoted by: 
\[
\aff_\R(\Z)=\R^{n} \rtimes \Gl(n,\Z).
\] 
Let us denote by 
\[
\aff (\R^n,\R^{n'})=\R^{n^{\prime}} \times \Hom ( \R^{n},\R^{n^{\prime}})
\] 
and by 
\[
\aff_\R (\Z^n,\Z^{n'})=\R^{n^{\prime}} \times \Hom ( \Z^{n},\Z^{n^{\prime}}).
\]
\begin{defi}  Let $B$ be a topological $n$-dimensional manifold.
\begin{itemize}
\item[(i)] An \textit{affine manifold} is a pair $(B,\mathscr{A})$ where $B$ is an $n$-dimensional manifold
and $\mathscr{A}$ is a maximal atlas on $B$ whose transition maps are $\aff(\R)$ transformations. We call $\mathscr A$ an \textit{affine structure} on $B$. 

\item[(ii)] An affine manifold $(B, \mathscr{A})$ is \textit{integral} if the transition maps of the affine structure $\mathscr A$ are $\aff_\R(\Z)$ transformations. We call $\mathscr A$ an \textit{integral affine structure} on $B$.
 
\item[(iii)] A continuous map $\alpha: B \rightarrow B^{\prime}$ is \textit{(integral) affine} if on each local coordinate chart, $\alpha$ is an element of ($\aff_\R(\Z^n,\Z^{n'})$) $\aff(\R^n,\R^{n'})$. Two (integral) affine manifolds $B$ and $B^{\prime}$ are said to be \textit{(integral) affine isomorphic} if there is an (integral) affine homeomorphism between them.
\end{itemize}
\end{defi}
It is becoming standard to call an affine manifold as in (ii) \textit{tropical manifold} \cite{Gross_SYZ_rev}. Though convenient for various good reasons, this is not a well established terminology at the time this paper is being written, so we prefer to stick to definition (ii) instead. Our convention coincides with that in  \cite{K-S-archim} and \cite{Haase-Zharkov} and differs from \cite{G-Siebert2003}. Affine manifolds whose structure group is $\aff (\Z)=\Z^{n} \rtimes \Gl(n,\Z)$ will be denoted $\aff(\Z)$\textit{-manifolds} (these are called integral affine in \cite{G-Siebert2003}).

\medskip
Given an affine manifold $(B,\mathscr A)$, consider a chart $(U,\phi)\in\mathscr A$ with affine coordinates $\phi=(u_1, \ldots, u_n)$. The cotangent bundle $T^{\ast}_B$ admits a flat connection $\nabla$ defined by
\[ \nabla du_{j} = 0, \]
for all $j=1, \ldots, n$ and all charts $(U,\phi)\in\mathscr A$. When $(B,\mathscr A)$ is integral affine we
can also define a maximal integral lattice $\Lambda \subset T^{\ast}_B$ by
\[ \Lambda|_{U} = \spn_{\Z} \inner{ du_1, \ldots}{du_n} \]
for all $(U,\phi)\in\mathscr A$. Therefore to every integral affine manifold $(B,\mathscr A)$ we can 
associate the $2n$-dimensional manifold
\[ X(B,\mathscr A) = T^{\ast}_B / \Lambda, \]
which together with the projection $f: X(B,\mathscr A) \rightarrow B$ forms
a $T^n$ fibre bundle. Also notice that the standard symplectic form
$\omega$ on $T^{\ast}_B$ descends to $X(B,\mathscr A)$ and the fibres of $f$
are Lagrangian. 

\medskip
The flat connection $\nabla$ on $T^{\ast}_B$ of an integral affine manifold $(B,\mathscr A)$ has a holonomy representation $\rho^\ast: \pi_1(B, b) \rightarrow\Gl(n, \Z)$ obtained by parallel transport along closed paths. A choice of basis of $\Lambda_b$ is identified naturally with a choice of basis of $H_1( f^{-1}(b), \Z)$. Under this identification, the holonomy representation $\rho^\ast$ coincides with the monodromy representation of the bundle $X(B,\mathscr A)\rightarrow B$. More precisely, if $g\in\pi_1(B,b)$ and $\mathcal M_b(g)$ is the corresponding monodromy matrix, then $\mathcal M_b(g)=\rho^\ast(g)$. The integral affine manifold $(B,\mathscr A)$ also induces a flat connection on $T_B$ whose holonomy representation, $\rho$, is dual to $\rho^\ast$, i.e. the matrix $\rho (g)$ is the inverse transpose of $\rho^\ast (g)$. In what follows, unless otherwise stated, ``holonomy representation" should be understood as the holonomy representation of the aforementioned flat connection on the \emph{cotangent} bundle $T^\ast_B$.

\medskip
It is well known that affine manifolds arise naturally from Lagrangian fibrations. This is the classical theory of action-angle coordinates in Hamiltonian mechanics.

\subsection{Action-angle coordinates.}
We review here some standard facts about Lagrangian fibrations
which we will use in the next Sections. For details we refer to Duistermaat \cite{Dui}.
Assume we are given a $2n$-dimensional symplectic manifold $X$
with symplectic form $\omega$, a smooth $n$-dimensional manifold $B$ 
and a proper smooth submersion $f: X \rightarrow B$ whose fibres are 
connected Lagrangian submanifolds.
For every $b \in B$, denote by $F_b$ the fibre of $f$ at $b$.
\begin{prop}[Arnold-Liouville] \label{lattice}
In the above situation, for every $b \in B$, $T^{\ast}_{b}B$ acts
transitively on $F_b$. In particular
there exists a maximal sub-lattice $\Lambda_b$ of $T^{\ast}_{b}B$ 
such that $F_b$ is naturally diffeomorphic to $T^{\ast}_{b}B / \Lambda_b$,
therefore $F_b$ is an $n$-torus.
\end{prop}
\begin{proof}
To every $\alpha \in T^{\ast}_{b}B$ we can associate a vector field
$v_{\alpha}$ on $F_b$ by
\[ \iota_{v_{\alpha}} \omega = f^{\ast} \alpha. \]
Let $\phi_{\alpha}^{t}$ be the flow of $v_{\alpha}$ with time
$t \in \R$. Then we define the action $\theta_{\alpha}$ of $\alpha$ 
on $F_b$ by
\[ \theta_{\alpha}(p) = \phi_{\alpha}^{1}(p), \]
where $p \in F_b$. One can check that such an action is well defined
and transitive. Then, $\Lambda_b$ defined as
\[ \Lambda_{b} = \{ \lambda \in T^{\ast}_{b}B \ | \ \theta_{\lambda}(p) = p,
\ \text{for all} \ p \in F_b \} \]
is a closed discrete subgroup of $T^\ast_bB$, i.e. a lattice.
From the properness of $F_b$ it follows that $\Lambda_b$ is maximal
(in particular homomorphic to $\Z^{n}$) and that $F_b$ is diffeomorphic
to $T^{\ast}_bB/ \Lambda_b$.
\end{proof}
We denote $\Lambda = \cup_{b \in B} \Lambda_b$. Given the presheaf on $B$ defined by 
$U \mapsto  H_1(f^{-1}(U), \Z)$, the associated sheaf is a locally constant sheaf. We 
can identify it with $\Lambda$ as follows. Let $U\subseteq B$ be 
a contractible open set. For every $b \in U$, $H_1(F_b, \Z)$ can be naturally 
identified with $H_1(f^{-1}(U), \Z)$. 
To every $\gamma \in H_1(f^{-1}(U), \Z)$, we can associate a $1$-form $\lambda$ on $U$ as 
follows. For every vector field $v$ on $U$, if we denote by $\tilde{v}$ a lift, define
\begin{equation} \label{per:def}
 \lambda(v) = - \int_{\gamma} \iota_{\tilde{v}} \omega. 
\end{equation}
It turns out that this identifies the above sheaf with $\Lambda\subset T^\ast_B$. 
If $\gamma_1, \ldots, \gamma_n$ are a basis for $H_1(F_b, \Z)$, then (\ref{per:def}) 
gives us a $\Z$-basis $\lambda_1, \ldots, \lambda_n$ of $\Lambda$ over a contractible open 
neighborhood $U$ of $b$. 

In particular,  one can read the monodromy of $f:X\rightarrow B$ from the monodromy of 
$\Lambda$. We state now the fundamental theorem of smooth proper Lagrangian submersions: 

\begin{thm}[Duistermaat] \label{lag:fund}
Given a basis $\{ \gamma_1,\ldots ,\gamma_n \}$ of $H_1(F_b, \Z)$, then the 
corresponding 1-forms $\lambda_1, \ldots, \lambda_n$ defined on a contractible 
open neighborhood $U$ of $b$ are closed and locally generate $\Lambda$. In particular, $\Lambda$ is Lagrangian with respect to the standard symplectic structure in $T^\ast_B$. A choice of functions $a_j$ such that $\lambda_j = da_j$ defines coordinates $a = (a_1, \ldots, a_n)$ called \textit{action coordinates}. A covering $\{ U_i \}$ of $B$ by contractible open sets and a choice of action coordinates on each $U_i$ defines an integral affine structure $\mathscr A$ on $B$. Moreover, if $f$ has a Lagrangian section $\sigma: U \rightarrow X$ over an open set $U\subseteq B$, then there is a natural symplectomorphism
\begin{equation}
\Theta:T^{\ast}_U / \Lambda\rightarrow f^{-1}(U). 
\end{equation}
If $\sigma$ is a global section then $X(B,\mathscr A)$ is symplectically conjugate to $X$. If in addition the monodromy of $\Lambda$ is trivial $X$ is symplectically conjugate to $B\times T^n$. The map $\Theta$ is called the \textit{period map} or \textit{action-angle} coordinates map.
\end{thm}
\begin{proof} We just give a sketch of the proof. Using the 
Weinstein neighborhood theorem one can show that in a 
sufficiently small tubular neighborhood of a fibre $F_b$,
the symplectic form is exact, i.e $\omega = - d\eta$
for some 1-form $\eta$. Notice that $\eta|_{F_b}$ is
a closed 1-form. Define functions $a_j$ on $U$ by
\[ a_j = \int_{\gamma_j} \eta. \]
One can show that
\[ \lambda_j = da_j \]
and therefore $\lambda_j$ is closed. It is clear that the coordinates $a = (a_1, ..., a_n)$
are well defined up to an integral affine transformation
and therefore they define an integral affine structure on $B$
inducing the lattice $\Lambda$ in $T^{\ast}_B$. 
Finally, notice that given a section $\sigma: U \rightarrow X$
we have a covering map
\[ \begin{array}{rcl}
        T^{\ast}_U & \rightarrow & f^{-1}(U) \\
          \alpha & \mapsto & \theta_{\alpha}(\sigma(\pi(\alpha)))
    \end{array}
                  \]
This map induces a diffeomorphism between $T^{\ast}_U / \Lambda$ and
$f^{-1}(U)$. One can check that in the case $\sigma$ is Lagrangian
this map is a symplectomorphism. For the proof of the last statement we refer the reader to \cite{Dui}. 
\end{proof}

\begin{cor} \label{aff:symp}
Let $\mathcal F=(X,f,B)$ and $\mathcal F'=(X',f',B')$ be smooth proper Lagrangian fibrations inducing integral affine 
structures $\mathscr{A}$ and $\mathscr{A}^{\prime}$ on $B$ and $B^{\prime}$ respectively. Assume there exist Lagrangian
sections $\sigma$ and $\sigma^{\prime}$ of $f$ and $f^{\prime}$ respectively. Then an integral affine isomorphism $\phi$ between $B$ and 
$B^{\prime}$ induces a symplectic $(\psi,\phi)$-conjugation between $\mathcal F$ and $\mathcal F'$ such that $\psi \circ \sigma = \sigma' \circ \phi$. 
\end{cor} 
\begin{proof} Let $\Lambda \subset T^{\ast}_B$ and 
$\Lambda^{\prime} \subset T^{\ast}_{B^{\prime}}$ be the
lattices induced from the integral affine structures on 
$B$ and $B^{\prime}$, respectively. From Theorem \ref{lag:fund} it follows that
$X$ and $X^{\prime}$ are symplectomorphic to 
$T^{\ast}_B / \Lambda$ and $T^{\ast}_{B^{\prime}} / \Lambda^{\prime}$, respectively.
Given an integral affine isomorphism $\phi$ between $B$ and $B^{\prime}$,
clearly $\phi^{\ast}$ is a symplectomorphism between $T^{\ast}_{B'}$ 
and $T^{\ast}_{B}$ inducing an isomorphism between $\Lambda^{\prime}$
and $\Lambda$. Therefore $\phi^{\ast}$ descends to a symplectomorphism
$\tilde\psi$ between $T^{\ast}_{B^{\prime}} / \Lambda^{\prime}$ and 
$T^{\ast}_B / \Lambda$. Defining $\psi=\Theta'\circ(\tilde\psi)^{-1}\circ\Theta^{-1}$ the claim follows.
\end{proof}
The following is an easy but important consequence of Arnold-Liouville-Duistermaat theorem:
\begin{cor}\label{smth:sg:equiv}
Proper Lagrangian submersions do not have semi-global symplectic invariants. In other words, all such fibrations are symplectically conjugate to $U\times T^n$ when restricted to a small enough neighborhood $U$ of a base point.  
\end{cor}
It is clear that \emph{smoothness} of the fibration map plays a crucial role in the above result. Semi-global invariants do arise for certain piecewise $C^\infty$ Lagrangian fibrations \cite{CB-M-stitched}. We say more about this in \S \ref{stitched fibr}.

\subsection{Affine manifolds with singularities.}
When a Lagrangian fibration has singular fibres, its base is no longer an affine manifold but an affine manifold \emph{with singularities.} These singularities can be a priori rather complicated. The topological properties described in \S \ref{sect. topology} motivate the following: 

\begin{defi}
An (integral) affine manifold with singularities is a
triple $(B, \Delta , \mathscr{A})$, where $B$ is a topological $n$-dimensional
manifold, $\Delta  \subset B$ a set which is locally a finite union of locally
closed submanifolds of codimension at least $2$ and $\mathscr{A}$ is an (integral) 
affine structure on $B_0 = B - \Delta $. A continuous map between (integral) affine manifolds with singularities
\[
\alpha: B\rightarrow B'
\] 
is (integral) affine if $\alpha^{-1}(B_0')\cap B_0$ is dense in $B$ and the restriction $\alpha_0=\alpha|_{\alpha^{-1}(B_0')\cap B_0}$:
\[
\alpha_0:\alpha^{-1}(B_0')\cap B_0\rightarrow B_0'
\]
is an (integral) affine map. We say that $\alpha$ is an (integral) affine isomorphism if $\alpha$ is an homeomorphism and $\alpha_0$ is an (integral) affine isomorphism of (integral) affine manifolds.
\end{defi}

From now on we restrict to  dimension $n=2$ or $3$. Let $(B, \Delta, \mathscr A)$ be an affine manifold with singularities and let $(B_0,\mathscr A)$ be the corresponding affine manifold. Let $X(B_0,\mathscr A)$ be the Lagrangian $T^n$ bundle over $B_0$ as introduced at the beginning of this section. We shall start imposing conditions on the singularities of the affine structure which, in particular, will imply that $X(B_0,\mathscr A)$ is of the topological type described in \S \ref{sect. topology}, e.g. such that $X(B_0,\mathscr A)$ will have semi-stable monodromy as in Theorem \ref{thm. semi-stable}. 

\medskip
We start defining local models of integral affine manifolds with singularities. In dimension 2, the allowed behavior is described in the following:

\begin{ex}[The node] \label{ff:aff} We define an affine structure with singularities on 
$B = \R^2$. Let $\Delta  = \{ 0 \}$ and let $(x_1,x_2)$ be the standard
coordinates on $B$. As the covering $\{U_i \}$ of $B_0 =\R^2 - \Delta$
we take the following two sets
\[ U_1 = \R^2 - \{ x_2 = 0 \ \text{and} \ x_1 \geq 0 \}, \]
\[ U_2 = \R^2 - \{ x_2 = 0 \ \text{and} \ x_1 \leq 0 \}. \]
Denote by $H^+$ the set $\{ x_2 > 0 \}$ and by $H^-$ the set $\{ x_2 < 0 \}$. 
Let $T$ be the matrix
\begin{equation} \label{ff:mon}
  T = \left( \begin{array}{cc}
                 1 & 0 \\
                 1 & 1
              \end{array} \right).
\end{equation}
The coordinate maps $\phi_1$ and $\phi_2$ on $U_1$ and $U_2$ are
defined as follows
\begin{eqnarray*}
 \phi_1 & = & \mathrm{Id} \\
 \ & \ & \\
 \phi_2 & = & \left \{ \begin{array}{ll}
                     \mathrm{Id} & \text{on} \ \bar H^+ \cap U_2, \\
                     (T^{-1})^t & \text{on} \ H^- 
                   \end{array} \right. 
\end{eqnarray*}
The atlas $\mathscr{A} = \{ U_i , \phi_i \}_{i=1,2}$
is clearly an affine structure on $B_0$. It is easy to check
that given a point $b \in B_0$, we can chose a basis of 
$T^{\ast}_{b} B_0$ with respect to which the holonomy 
representation $\rho^\ast$ sends the anti-clockwise oriented generator of 
$\pi_1 (B_0)$ to the matrix $T$. 
\end{ex}

In dimension $3$ we have the following models.

\begin{ex}[The edge] \label{edge:aff} Let $I \subseteq \R$ be an open interval. 
Consider $B = \R^2 \times I$ and $\Delta  = \{ 0 \} \times I$.
On $B_0 = (\R^2 - 0) \times I$ we take the product affine structure
between the affine structure on $\R^2 - 0$ described in the previous 
example and the standard affine structure on $I$. 
\end{ex}

\begin{ex}[A variation] \label{edge:affvar}
In the previous example the discriminant locus $\Delta$ was a straight line. 
We can slightly perturb $\Delta$ so that it becomes a smooth 
curve. More precisely, let $B = \R^2 \times I$ as before and  consider a
smooth function $\tau: I \rightarrow \R$.  Let 
\[ \Delta_\tau = \{ (\tau(s),0, s), s \in I \} \subset B \]
and define a covering $\{ U_i \}$ of $B_0 = B - \Delta_\tau$ to be
\[ U_1 = (\R^2 \times I) - \{ (x_1, 0 , s) \ | \ x_1 \geq \tau(s) \},\]
\[ U_2 = (\R^2 \times I) - \{ (x_1, 0, s) \ | \ x_1 \leq \tau(s) \}.\]
Now let $H^+ = \{ x_2 > 0 \}$ and $H^- = \{ x_2 < 0 \}$. Take the following 
matrix
\begin{equation*}
  T = \left( \begin{array}{ccc}
                 1 & 0 & 0 \\
                 1 &  1 & 0 \\
                 0 &  0 & 1
              \end{array} \right)
\end{equation*}
and define maps $\phi_j$ on $U_j$ to be
\begin{eqnarray*}
 \phi_1 & = & \mathrm{Id} \\
 \ & \ & \\
 \phi_2 & = & \left \{ \begin{array}{ll}
                     \mathrm{Id} & \text{on} \ \bar H^+ \cap U_2, \\
                     (T^{-1})^t & \text{on} \ H^-. 
                   \end{array} \right. 
\end{eqnarray*}
Clearly $\mathscr A = \{ U_i, \phi_i \}_{i =1,2}$ defines an affine structure 
on $B_0 = B - \Delta_\tau$. When $\tau = 0$, this example coincides with 
the previous one. Notice that the curve $(\tau(s), 0 , s)$ is contained 
inside the $2$-plane $\{x_2=0\}$, which can be viewed as an integral surface of the distribution
spanned by the vectors in $TB_0$ which are invariant with respect to the 
holonomy representation $\rho$ on $TB_0$. Two different curves give non-isomorphic 
singular affine structures, unless the curves can be taken one into the other
via an integral affine transformation. 
\end{ex}

\begin{ex}[Positive vertex] \label{pos:aff} Take $B = \R \times \R^2$, with coordinates 
$(x_1, x_2, x_3)$ and identify $\R^2$ with $\{ 0 \} \times \R^2$. 
Inside $\R^2$ consider the cone over three points: 
 \[  \Delta  =  \{ x_2 = 0, \, x_3 \leq 0 \} \cup 
              \{ x_3=0, \, x_2 \leq 0 \} 
               \cup \{x_2 = x_3, \, x_3 \geq 0 \}. \] 
Now define closed sets in $B$
\begin{eqnarray*}
R & = & \R \times \Delta, \\
R^+ & = & \R_{\geq 0} \times \Delta, \\
R^- & = & \R_{\leq 0} \times \Delta, 
\end{eqnarray*}
and consider the following cover $\{ U_i \}$ of $\R^3 - \Delta$:
\begin{eqnarray*}
  U_1 & = & \R^3 - R^+,  \\ 
  U_2 & = & \R^3 - R^-.
\end{eqnarray*}
It is clear that $U_1 \cap U_2$ has the following three connected components
\begin{eqnarray*}
  V_1 & = & \{ x_2 < 0, \ x_3 < 0 \},  \\ 
  V_2 & = & \{ x_2 > 0, \ x_2 > x_3 \},  \\
  V_3 & = & \{ x_3 > 0, \ x_3 > x_2 \}. 
\end{eqnarray*}
Take two matrices
\begin{equation}
 T_1  =  \left( \begin{array}{ccc}
                 1 & 1 & 0 \\
                 0 & 1  & 0 \\
                 0 & 0  & 1 
              \end{array} \right), \ \ \ 
 T_2   = \left( \begin{array}{ccc}
                 1 & 0 &  -1 \\
                 0 & 1  & 0 \\
                 0 & 0  & 1 
              \end{array} \right). \label{t12}
\end{equation}
Now on $U_1, U_2$ we define coordinate maps $\phi_1$, $\phi_2$ as 
follows
\begin{eqnarray*}
\phi_1 & = & \I, \\
\phi_2 & = & \left \{ \begin{array}{ll}
                     \I & \text{on} \  \bar V_1 \cap U_2, \\
                     T_1^{-1}  & \text{on} \  \bar V_2 \cap U_2 \\
                     T_2 & \text{on} \  \bar V_3 \cap U_2 
                   \end{array} \right. 
\end{eqnarray*}
Again we see that $\mathscr{A} = \{ U_i, \phi_i \}_{i = 1,2}$ 
gives an affine structure on $B_0 = \R^3 - \Delta $. One can compute
that given a point $b \in B_0$ and closed paths $g_1$, $g_2$ and
$g_3$ as in Figure~\ref{loops}, we can choose a basis of 
$T^{\ast}_{b} B_0$ with respect to which the holonomy matrices satisfy $\rho^\ast(g_j)=(T_j^{-1})^t$ for $j=1,2,3$. 
\end{ex}

\begin{ex}[A variation] \label{pos:affvar}
In the previous example, $\Delta$ was a graph with three edges meeting in 
one vertex. All three edges were straight lines. In the spirit of Example~\ref{edge:affvar}
we can perturb each edge of $\Delta$ to a smooth curve starting at the vertex.
Each straight edge of the previous example is contained in a $2$-plane which is an 
integral plane of the distribution spanned by the vectors which are invariant with respect to 
the holonomy around that edge. For example, consider the edge $E_1 = \{x_1=x_2 = 0, x_3 \leq 0 \}$ 
of $\Delta$. Then $E_1$ is contained inside the half plane, $P_1=\{x_2 = 0, x_3 \leq 0 \}$, 
whose tangent vectors are $T_1$ invariant, where $T_1=\rho (g_1)$ is the holonomy of $T_{B_0}$ with respect to $E_1$. An analogous thing happens with the other two edges. The union of all three half planes
gives $R$. The new perturbed edges, $E_j'$, must be curves inside the half planes $P_j$.
More precisely, let  $\tau$ be a function on $\Delta$ which 
is the restriction of a smooth function defined on an open neighborhood of $\Delta$, such that $\tau(0) =0$. If we let $R$ be as in the previous example, define
\begin{eqnarray*}
 \Delta_{\tau} & = & \{ (\tau(q), q) \in \R \times \Delta \} \\
  R^+ & = & \{ ( x_1 , q) \in \R \times \Delta \ | x_1 \geq \tau(q) \} \\
  R^- & = & \{ ( x_1 , q) \in \R \times \Delta \ | x_1 \leq \tau(q) \}
\end{eqnarray*}
Now charts $\mathscr A = \{ U_i, \phi_i \}_{i=1,2}$ on $B - \Delta_{\tau}$ can be defined like in 
the previous example, but with these new definitions of $R^+$ and $R^-$. It is clear that
$( B, \Delta_{\tau}, \mathscr A)$ defines an affine manifold with singularities. Two different 
choices of functions $\tau$ define non-isomorphic integral affine manifolds with singularities,
unless their graphs inside $R$ can be mapped one to the other via an integral affine map. 
\end{ex}

\begin{ex}[Negative vertex] \label{neg:aff} Let $B$ and $\Delta $ be as in Example~\ref{pos:aff}. 
Clearly, $\R^2 - \Delta$ has three connected components, which we denote $C_1, C_2$ and 
$C_3$. Let $\bar C_j=C_j\cup\partial C_j$. Viewing $\R^2$ embedded in $B$ as $\{0 \} \times \R^2$, consider the following three open 
subsets of $B_0$: 
\begin{eqnarray*}
  U_1 & = & \R^{3} - ( \bar C_2 \cup \bar C_3),\\ 
  U_2 & = & \R^{3} - 
               ( \bar C_1 \cup \bar C_3),  \\
  U_3 & = & \R^{3} - 
               ( \bar C_1 \cup \bar C_2 ). 
\end{eqnarray*}
Let 
\begin{eqnarray*}
  V^+ & = & \{ x_1 > 0 \},  \\ 
  V^- & = & \{ x_1 < 0 \}. 
\end{eqnarray*}
Clearly $U_i \cap U_j = V^+ \cup V^-$ when $i \neq j$. 
If $T_1$ and $T_2$ are as in (\ref{t12}),
define the following coordinate charts on $U_1$, $U_2$, $U_3$ respectively:
\begin{eqnarray*}
\phi_1 & = & \I, \\
\phi_2 & = & \left \{ \begin{array}{ll}
                     (T_1^{-1})^t & \text{on} \  \bar V^+ \cap U_2 \\
                     \I & \text{on} \  \bar V^- \cap U_2
                   \end{array} \right.  \\
\phi_3 & = & \left \{ \begin{array}{ll}
                     \I & \text{on} \  \bar V^+ \cap U_3 \\
                     (T_2^{-1})^t & \text{on} \  \bar V^- \cap U_3
                   \end{array} \right. 
\end{eqnarray*}
We can check that the affine structure defined by these charts is such that,
for fixed $b \in B_0$, there exists a basis of $T^{\ast}_{b} B_0$ with respect 
to which the holonomy representation is such that $\rho^\ast(g_j)=T_j$, where $g_j$ are as
in Figure \ref{loops}. In particular, the holonomy is given by the inverse transpose matrices of the
holonomy in the previous example. 
\end{ex}

\begin{ex}[A variation] \label{neg:affvar}
Again, we can perturb the above example by replacing the straight edges of 
$\Delta$ with smooth curves starting at the origin. This time these curves
have to be contained inside $\{ 0 \} \times \R^2$, which is the integral 
surface (containing $\Delta$) of the distribution spanned by the $\rho$-holonomy invariant 
vectors in $TB_0$. The perturbed $\Delta$, which we could denote $\Delta_{\tau}$,
still separates $\R^2$ in three connected
components $C_1, C_2$ and $C_3$. Then the definition of the affine structure 
carries through just like in the previous example and we denote it by $\mathscr A_{\tau}$. 
\end{ex}

We are now ready to give a definition of the specific affine structures 
with singularities which we will consider.

\begin{defi}\label{def:simple} A $2$-dimensional affine manifold with singularities
$(B,\Delta , \mathscr{A})$ is said to be \textit{simple} if $\Delta$ consists
of a finite union of isolated points and a neighborhood 
of each $p\in\Delta$ is affine isomorphic to a neighborhood of $0\in\R^2$ as in Example \ref{ff:aff}. We call $p\in\Delta$ a \textit{node}.
A $3$-dimensional affine manifold with singularities $(B,\Delta , \mathscr{A})$
is \textit{simple} if it satisfies:
\begin{itemize}
\item[(i)] $\Delta $ is a trivalent graph; 
\item[(ii)] a neighborhood of each vertex of $\Delta$ is affine isomorphic to a neighborhood of $0\in\R^3$ in either Examples \ref{pos:aff} or \ref{pos:affvar}, in which case we call it a \textit{positive vertex}; or to a neighborhood of $0\in\R^3$ in either Examples \ref{neg:aff} or \ref{neg:affvar}, in which case we call it a \textit{negative vertex};
\item[(iii)] a neighborhood of each edge of the graph is affine isomorphic to a neighborhood of $\Delta$ in Example \ref{edge:aff}; or a neighborhood of $\Delta_\tau$ in Example \ref{edge:affvar} for a suitable $\tau$. 
\end{itemize}
\end{defi}

The following is direct consequence of the above definition and Theorem \ref{thm. semi-stable}:

\begin{cor} Let $(B,\Delta , \mathscr{A})$ be a simple affine manifold with singularities and let $(B_0,\mathscr A)$ be the underlying integral affine manifold. Then
\[
f_0:X(B_0,\mathscr A)\rightarrow B_0
\]
is a $T^n$ bundle with semi-stable monodromy as in Theorem \ref{thm. semi-stable}. In particular, there is an $2n$-manifold $X$ and a topological semi-stable compactification $X(B_0,\mathscr A)\hookrightarrow X$. Furthermore, the topological fibration $f:X\rightarrow B$ obtained is topologically simple.
\end{cor}

\subsection{Examples}
Here we give some examples of affine manifolds with singularities
and then we prove the $2$-dimensional version of the main theorem
of this article.

\begin{figure}[!ht] 
\begin{center}
\input{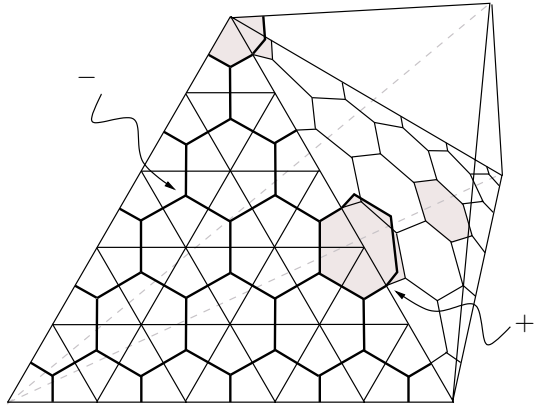}
\caption{Affine $S^3$ with singularities.}\label{fig:affine_s3}
\end{center}
\end{figure}

\begin{ex} \label{k3:aff} 
In $\R^3$ consider the $3$-dimensional simplex $\Xi$
spanned by the points
\[ P_0=(-1,-1,-1), \ \ P_1 = (3, -1, -1), \ \ P_2 = (-1, 3, -1),
      \ \ P_3 = (-1, -1, 3). \]
Let $B = \partial \Xi$. We explain how to construct a simple affine
structure with singularities on $B$. 
Each edge $\ell_j$ of $\Xi$ has $5$ integral points (i.e. belonging
to $\Z^n$), which divide $\ell_j$ into $4$ segments. 
For each $j =1, \ldots, 6$ denote by $\Delta ^{j}_{k}$, $k=1, \ldots, 4$ 
the four barycenters of these four segments. We let
\[ \Delta  = \{ \Delta ^{j}_{k}; j= 1 \ldots 6 \ \text{and} \ k=1, \ldots, 4 \}. \]
A covering of $B_0=B-\Delta$ can be defined as follows. The first four
open sets consist of the four open faces $\Sigma_i$, $i=1 \ldots, 4$
with the affine coordinate maps $\phi_i$ induced by their affine embeddings 
in $\R^3$. Denote
by $I$ the set of integral points of $B$ which lie on an edge.
For every $Q \in I$ we can choose a small open set $U_Q$ in $B_0$
such that $\{\Sigma_i \}_{i=1, \ldots, 4} \cup \{ U_{Q} \}_{Q \in I}$
is a covering of $B_0$. Let $R_Q$ denote the $1$-dimensional subspace
of $\R^3$ generated by $Q \in I$. 
One can verify that if $U_Q$ is small enough, the projection 
$\phi_Q: U_Q \rightarrow \R^{3} / R_{Q}$ is an homeomorphism.
A computation shows that the atlas $\mathscr{A} = 
\{\Sigma_i, \phi_i \}_{i=1, \ldots, 4} \cup \{ U_{Q}, \phi_Q \}_{Q \in I}$
defines an affine structure on $B_0$ making $(B,\Delta, \mathscr A)$ simple. 
\end{ex}
\begin{ex} \label{quint:aff} 
This three dimensional example is taken from \cite{GHJ} \S 19.3. Let $\Xi$ be the $4$-simplex in $\R^3$ spanned by
 \begin{eqnarray*}
      & P_0=(-1,-1,-1, -1), \ P_1 = (4,-1,-1, -1), \ P_2 = (-1, 4,-1,-1), 
                                              \\
& P_3 = (-1, -1, 4, -1), \ \
                      P_4 = (-1, -1, -1, 4).
\end{eqnarray*}
Let $B = \partial \Xi$. Denote by $\Sigma_j$ the open $3$-face of $B$
opposite to the point $P_j$ and by $F_{ij}$ the closed $2$-face separating 
$\Sigma_i$ and $\Sigma_j$. Each $F_{ij}$ contains $21$ integral points
(including those on its boundary). These form the vertices of a
triangulation of $F_{ij}$ as in Figure \ref{fig:affine_s3}. By joining the barycenter
of each triangle with the barycenters of its sides we form
a trivalent graph as in Figure \ref{fig:affine_s3}. Define the set $\Delta $ to be the union
of all such graphs in each $2$-face. Denote by $I$ the set of
integral points of $B$. Just as in the previous example,
we can form a covering of $B_0 = B - \Delta $ by taking the open $3$-faces
$\Sigma_j$ and small open neighborhoods $U_Q$ inside $B_0$ of $Q \in I$.
A coordinate chart $\phi_i$ on $\Sigma_i$ can be obtained from its 
affine embedding in $\R^4$. If we denote again by $R_Q$ the linear
space spanned by $Q \in I$, as a chart on $U_Q$ we take the projection
$\phi_Q: U_Q \rightarrow \R^4 / R_Q$. A computation shows that this affine structure is simple. In fact the vertices of $\Delta $ which are contained in the interior of each $2$-face are of negative type and those which are contained in the $1$-faces are of positive type. 
\end{ex}

\begin{ex}[A variation]\label{quint:affvar} 
In the previous example, all edges of $\Delta$ were straight lines, but one can perturb them 
in the sense of Examples~\ref{edge:affvar}, \ref{pos:affvar} and \ref{neg:affvar}. In fact we can 
form a new $\Delta$ by keeping the vertices fixed and connecting them through smooth 
curves, which are small perturbations of the straight edges of the previous example. 
If these curves stay inside the $2$-faces of $B$, then the affine structure on $B - \Delta$ 
can be defined just like above. 
\end{ex}

In some cases, such as in Examples~\ref{k3:aff} and \ref{quint:aff} given an affine manifold with singularities, one can define a second affine structure $\check{\mathscr A}$ on $B$, via a discrete Legendre transform of $\mathscr A$ (cf. Gross and Siebert \cite{G-Siebert, G-Siebert2003}). Here we shall not give details about how this process works. Though it is important to mention that this method produces a second integral affine manifold with singularities, $(B,\check\Delta ,\check{\mathscr A} )$ which coincides topologically with $(B,\Delta ,\mathscr A )$ but with holonomy representation $\check\rho$ dual to $\rho$. In dimension 3 this means, in particular, that the positive vertices of $\Delta$ become negative vertices of $\check \Delta$ and vice-versa. 

\medskip
These examples of singular affine manifolds are very important. The bundles associated to them satisfy the hypothesis of Theorem \ref{thm. semi-stable} so they can be used to produce topological semi-stable compactifications which are homeomorphic to well known examples of Calabi-Yau manifolds:

\begin{thm}[Gross \cite{TMS}] \label{quint:mirr}
Let $(B, \Delta , \mathscr{A})$ be the integral affine manifold with singularities
described in Example \ref{quint:aff} and let
\[
(B, \Delta , \mathscr{A})\rightarrow (B,\check\Delta ,\check{\mathscr A} )
\] 
be its Legendre transform. Let $X(B_0,\mathscr A)\hookrightarrow X$ and $X(B_0,\check{\mathscr A})\hookrightarrow\check X$ be the corresponding topological semi-stable compactifications. Then $X$ is homeomorphic to the quintic hypersurface and $\check X$ is homeomorphic to its mirror. 
\end{thm}

Later in this article we show that there are \emph{symplectic} semi-stable compactifications recovering the quintic and its mirror. These compactifications rely deeply on the existence of suitable local models of Lagrangian fibrations with singular fibres. The construction of such models is a highly delicate issue. 

\subsection{The focus-focus fibration}
In dimension 2 it is much easier to produce symplectic semi-stable compactifications. Now we will show how Example \ref{k3:aff} gives rise to a symplectic semi-stable compactification diffeomorphic to a K3 surface. This will require a local model of Lagrangian $T^2$ fibration with a semi-stable singular fibre, such as the one in the following:

\begin{ex} \label{smooth:ff}
Let $X = \C^{2} - \{ z_1 z_2 + 1 = 0 \}$ and let $\omega$ be the
restriction to $X$ of the standard symplectic form on $\C^2$. One
can easily check that the following map $f:X\rightarrow\R^2$ is a Lagrangian fibration: 
\begin{equation} \label{sm:ff}
f(z_1,z_2)=\left( \frac{|z_1|^2-|z_2|^2}{2}, \, \log |z_1z_2 +1| \right).
\end{equation}
The only singular fibre is $f^{-1}(0)$, which has the topology of a $I_1$ fibre.  It follows that this fibration is conjugate to the topological fibration in  Example \ref{ex top ff}. 
\end{ex}

Lagrangian fibrations with semi-stable singular fibres, e.g., conjugate to the fibration in Example \ref{ex top ff}, are called \textit{focus-focus fibrations}. They have been studied extensively in Hamiltonian Mechanics \cite{Dui}, \cite{Zung1} --where they got their name-- and more recently in symplectic topology \cite{LeungSym}, \cite{San} and Mirror Symmetry \cite{G-Wilson2}.

\medskip
Let $\arg: \C^{\ast} \rightarrow \R$ be the multi-valued function 
$\rho e^{i \theta} \mapsto \theta$. Denote by $D\subseteq\C$ the unit open disk and let $D^{\ast} = D -\{ 0 \}$. Let $\mathcal F=(X,\omega, f, D)$ be a focus-focus fibration.
It has been shown \cite{San} that there are coordinates 
$b=(b_1, b_2)$ on $\R^2$, 
with values in $D$, a smooth function 
$q:  D \rightarrow \R$ such that $q(0)=0$ and a choice of generators of
$H_1(f^{-1}(b), \Z)$ with respect to which the periods $\lambda_1$
and $\lambda_2$ of $\mathcal F$ can be written as 
\begin{eqnarray*}
\lambda_1 & = & - \log |b| \ db_1 + \arg b \ db_2 + dq \\
\lambda_2 & = & 2 \pi  \, db_2.
\end{eqnarray*}
Clearly $\lambda_1$ is multi-valued and blows up as $b\rightarrow 0$. The lattice
\[ \Lambda = \spn_{\Z} \inner{\lambda_1}{\lambda_2} \]
has monodromy given by  $T$ as in (\ref{ff:mon}).
We now describe the affine structure induced on $D^{\ast}$.
Consider the two open subsets
 \begin{eqnarray*}
 U_1 & = & D - \{ \im b = 0 \ \text{and} \re b \geq 0 \}, \\
 U_2 & = & D - \{ \im b = 0 \ \text{and} \re b \leq 0 \}.
 \end{eqnarray*}
On $U_1$ we chose the branch of $\arg$ with values in 
$(0, 2\pi)$ and we denote it by $\arg_1$. On $U_2$
we chose the branch with values in $(- \pi, \pi)$ which
we denote by $\arg_2$. 
Clearly on $U_1 \cap U_2$ we have $\arg_1 = \arg_2 + 2\pi$. 
A computation shows
that the maps $\psi_j: U_j \rightarrow \R^2$ given
by
\[ \psi_j(b) = (-b_1 \log |b| + b_1 + q(b) + 
                      b_2 \arg_j b, \, 2 \pi b_2 ), \]
with $q(0) = 0$, are a choice of affine coordinates
associated to $\lambda_1$ and $\lambda_2$. 

\medskip
It is easy to check that the map $\psi_1$ (or $\psi_2$)
extends continuously to $D$. 
Call $\alpha: D \rightarrow \R^2$ the extended map. 
On a sufficiently small neighborhood $V\subseteq D$ of $0$, the map $\alpha$ 
is a homeomorphism of $V$ onto $\alpha(V)$. The reader may verify that 
$0\in V$ is a node with respect to the affine structure given by 
$\{ U_j, \psi_j \}$. In other words, the map $\alpha$ restricted to 
$V^{\ast} = V - \{ 0 \}$ is an affine isomorphism between $V^\ast$ and the affine manifold $\alpha(V^{\ast})$ whose affine structure is the restriction of the one in Example~\ref{ff:aff}.
The affine structure with singularities on $D$ induced by a focus-focus 
fibration is therefore simple. In particular, the affine structure induced by 
Example \ref{smooth:ff} is simple.

\medskip
\begin{rem}\label{rem:san}
Germs of focus-focus fibrations --with respect to symplectic conjugation-- are classified by formal power series in two variables $\R[\![x,y]\!]$ with vanishing constant term \cite{San}. Such series correspond to the Taylor coefficients of functions $q\in C^\infty (D)$ as above evaluated at $0\in\R^2$. This means that there is an infinite number of \emph{different} germs of focus-focus fibrations, all inducing simple affine manifolds with singularities, i.e. inducing \emph{the same} singular affine structure on the base.  In \S \ref{pos:gener} we will see that a similar phenomenon happens in higher dimensions.
\end{rem}

\subsection{The K3 surface.}

\begin{thm}\label{thm. the K3} Let $(B, \Delta , \mathscr{A})$ be the affine manifold
with singularities in Example~\ref{k3:aff} and let $X(B_0,\mathscr{A})$ be the associated $T^2$ bundle with symplectic structure $\omega_0$ and projection $f_0$ induced by the standard ones in 
$T_{B_0}^\ast$.  There exists a compact symplectic manifold $(X, \omega)$, a Lagrangian
fibration $f: X \rightarrow B$ and an embedding $\iota: X(B_0,\mathscr A) \rightarrow X$ such that 
$\iota^{\ast} \omega = \omega_0$ and $f \circ \iota = f_0$.
Moreover $X$ is diffeomorphic to a smooth $K3$ surface. 

\end{thm}
\begin{proof} Let $f_V: X_V \rightarrow V$ be a focus-focus fibration over a small open neighborhood $V$ of its node $0\in V$. Let $V^\ast=V-\{0\}$ and denote by $(V^{\ast},\mathscr A_V)$ the integral affine manifold induced by $f_V$. Let $X(V^\ast,\mathscr A_V)$ be the associated Lagrangian $T^2$ bundle over $V^\ast$. It can be shown that $f_V$ has a Lagrangian section $s: V \rightarrow X_V$ such that $s(V)\cap\Crit (f_V)=\varnothing$. Then from Theorem~\ref{lag:fund} it follows that $f^{-1}_V(V^\ast)\subset X_V$ is symplectically conjugate to $X(V^\ast, \mathscr A_V)$.

\medskip
Now let $P \in \Delta $ and let $U\subset B$ be a small neighborhood of $P$. Denote by $U^{\ast} = U - P$ and by $X(U^\ast,\mathscr A)$ the Lagrangian $T^2$ bundle over $U^\ast$ given by the restriction of $X(B_0,\mathscr A)$ to $U^\ast$.  Recall that both $U$ and $V$ are simple affine manifold with singularities. Then, after taking $U$ and $V$ small enough, there exists an integral affine isomorphism $V^{\ast}\cong U^{\ast}$.  From Corollary \ref{aff:symp}, the latter isomorphism induces is a symplectic conjugation,
\[
f^{-1}_V(V^\ast)\cong X(V^\ast, \mathscr A_V)\cong X(U^\ast,\mathscr A),
\] 
which can be used to symplectically glue $X_V$ to $X(B_0)$. Define $(X,\omega )$ to be the symplectic manifold obtained after applying this gluing over all points $P \in \Delta $ and $f:X\rightarrow B$ the resulting fibration. It is clear that $(X,\omega)$ is a semi-stable compactification of $(X(B_0,\mathscr A),\omega_0)$ such that $\iota^\ast\omega =\omega_0$. It is easy to check that $(X, \omega, f, B)$ is topologically conjugate to a simply connected elliptic fibration with 24 singular fibres of type $I_1$. It follows that $X$ is diffeomorphic to a K3 surface.
\end{proof} 

\begin{cor}
In view of Remark \ref{rem:san}, given $(B, \Delta , \mathscr{A})$ as in Example~\ref{k3:aff}, a compactification $X(B_0,\mathscr A)\hookrightarrow (X,\omega)$ as above is uniquely determined up to symplectic conjugation by a choice of 24 formal power series in two variables:
\[
\mathfrak q_1, \ldots ,\mathfrak q_{24}\in\R[\![x,y]\!]
\]
corresponding to germs of focus-focus fibrations $\mathcal F_1,\ldots \mathcal F_{24}$. In particular, there are infinitely many Lagrangian fibrations of a symplectic K3 surface, fibering over $(B, \Delta , \mathscr{A})$, which are all topologically conjugate but not \emph{symplectically} conjugate. 
\end{cor}

The space $\R[[x,y]]$ being contractible, implies that every two focus-focus fibrations can be connected with a path in $\R[[x,y]]$. The standard Moser's argument implies that the corresponding total spaces are symplectomorphic. Similarly, any two symplectic structures obtained using Theorem \ref{thm. the K3} can be connected with a path in $\R[[x,y]]^{24}$. Moser's argument implies that all such manifolds are symplectomorphic. \todo{can you find the reference?}

\medskip
Following an alternative approach, Zung obtained a Lagrangian fibration of a symplectic 4-manifold  which is also diffeomorphic to a K3 surface (cf. \cite{Zung1-II}{Example 4.19}). Leung and Symington \cite{LeungSym} use affine geometry as starting point to construct and classify --up to diffeomorphism-- the so-called \textit{almost toric symplectic 4-manifolds.} The fibration we obtained in Theorem \ref{thm. the K3} coincides with one of the list in \cite{LeungSym}. 

\medskip
 Other ways of constructing affine manifolds with singularities
have been proposed by Gross and Siebert \cite{G-Siebert, G-Siebert2003}, Hasse and Zharkov \cite{Haase-Zharkov, Haase-ZharkovII, Haase-ZharkovIII}. 
In \cite{Gross_Batirev}, Gross finds a combinatorial method to obtain simple affine manifolds 
with singularities out of the geometry of the polytopes which Batyrev and Borisov use to 
construct 
pairs of Calabi-Yau varieties as complete intersections inside Fano toric varieties. 
From Theorem 0.1 of \cite{Gross_Batirev} (proved by Gross and Siebert in \cite{GroSie_Tor}) 
it follows that these affine manifolds give rise to topological semi-stable compactifications 
homeomorphic to the two Batyrev-Borisov's Calabi-Yau varieties. We shall see in this paper that similar compactifications can be carried out in the symplectic category.

\section{Positive and generic-singular fibrations.}
\label{pos:gener}

We describe some of the local models needed to produce symplectic compactifications. These models may be regarded as 3-dimensional analogues to focus-focus fibrations. The arguments given here can be generalized to dimension $n>3$.   All fibrations in this Section are given by smooth maps.

\begin{defi} Let $\mathcal F=(X,\omega,f,B)$ be a Lagrangian fibration.
\begin{itemize}

\item[(i)] A \textit{Lagrangian generic-singular fibration} is a smooth Lagrangian fibration $\mathcal F$, 
with non-degenerate singularities (in the sense of \cite{Zung6}) which is conjugate to a topological $T^3$ 
fibration of generic type (cf. Example \ref{ex. (2,2)}). 

\item[(ii)] A \textit{Lagrangian positive fibration} is a Lagrangian fibration $\mathcal F$ which is conjugate to a topological $T^3$ fibration of positive type (cf. Example \ref{ex. (1,2)}).
\end{itemize}
\end{defi}

The non-degeneracy condition implies that the singularity is of rank-1 focus-focus type, such singularities are normalized \cite{Zung6}.

\subsection{Examples}

We start giving examples of non-proper Lagrangian fibrations describing the singular behavior of (i) and (ii) near $\Crit (f)$. Let $D^k\subseteq \R^k$ be the standard open ball.

\begin{ex}\label{ex. focus-focus}
Consider $\R^4$ with standard coordinates $(x_1, x_2, y_1, y_2)$ and let $D^4\subseteq\R^4$. Let $D^{1}\times S^1$ have coordinates $(r,\theta )$. Define $V=D^4\times D^{1}\times S^1$ with the standard symplectic structure and $F(x_i,y_i,r,\theta)=(b_1,b_2,b_3)$ where
\begin{equation}\label{coord}
\begin{array}{lll}
b_1=x_1y_1+x_2y_2, &b_2=x_1y_2-x_2y_1, & b_3=r_3.
\end{array}
\end{equation}
The reader may verify that $\mu=(b_2,b_3)$ is the moment map of a Hamiltonian action of $T^{2}$ and that $F$ is a $T^{2}$ invariant Lagrangian fibration of $V$ over $D^2\times D^{1}$. The singular fibres are homeomorphic to $\R\times S^1\times S^1$ after $\{ p\}\times S^1\times S^1$ is collapsed to $\{p\}\times S^1$. 
\end{ex}

\begin{ex}\label{ex. HL}
Consider $\C^3$ with canonical coordinates $z_1,z_2,z_3$. 
Define $F(z)=(b_1,b_2,b_3)$, where
\begin{equation}
\begin{array}{lll}
b_1=\im z_1z_2z_3,& b_2=|z_1|^2-|z_2|^2,& b_3=|z_1|^2-|z_3|^2.
\end{array}
\end{equation}
Here $\mu (z_1,z_2,z_3)=(b_2,b_3)$ is the moment map of a $T^2$-action,
furthermore the above functions Poisson commute, so the fibres of
$F$ are Lagrangian.  The critical locus of $F$
is $\Crit  (F)=\bigcup_{ij}\{ z_i=z_j=0\}$ and its discriminant
locus is $\Delta=\{ b_1=0, b_2=b_3\geq 0\}\cup\{ b_1=b_2=0,
b_3\leq 0\}\cup\{ b_1=b_3=0, b_2\leq 0\}$, i.e. a cone over three
points with vertex at $0\in\R^3$. The regular fibres are
homeomorphic to $\R\times T^2$. The singular fibre over
$0\in\Delta$ is homeomorphic to $\R\times T^2$ after $\{ p\}\times
T^2$ is collapsed to $p\in\R$. All the other singular fibres are
homeomorphic to $\R\times T^2$ after a two cycle $\{ p \} \times T^2 \subset \R \times T^2$ is collapsed to a circle. This is one of the
examples of special Lagrangian fibrations by Harvey and
Lawson \cite{HL}. 
\end{ex}

Now we give explicit examples of Lagrangian positive and generic-singular fibrations.
\begin{ex} \label{proper:pos}
Let $X=\C^3-\{ 1+ z_1z_2z_3=0\}$ with canonical coordinates $z_1, z_2, z_3$ and the standard symplectic structure. Consider the $T^2$-action on $X$ given by
$(z_1,z_2,z_3)\mapsto (e^{i\theta_1}z_1,e^{i\theta_2}z_2,e^{-i(\theta_1+\theta_2)}z_3)$. We obtain $f:X\rightarrow\R^3$ given by $f=(f_1,f_2,f_3)$ where
\begin{center}
\begin{tabular}{lll}
$f_1=\log |1+z_1z_2z_3|$, & $f_2=|z_1|^2-|z_2|^2$, & $f_3=|z_1|^2-|z_3|^2$.
\end{tabular}
\end{center}
It is straightforward to check that the above functions Poisson commute, hence the fibres of $f$ are Lagrangian. It follows that $f$ is modeled on Example \ref{ex. HL} near $\Crit (f)$. In particular, the discriminant locus is a cone over three points which coincides with the one in Example \ref{ex. HL}. 
This example has the topology of a positive fibration.
\end{ex}

\begin{ex} \label{smooth:generic} 
Let $X'=\C^2-\{ z_1z_2 - 1=0\}$ and let $X=X'\times\C^\ast$ with the standard symplectic structure. Define $f:X\rightarrow\R^3$ by $f=(f_1,f_2,f_3)$ where 
\begin{center}
\begin{tabular}{lll}
$f_1=\frac{|z_1|^2-|z_2|^2}{2}$, & $f_2=\log |z_3|$, & $f_3=\log |z_1z_2-1|$.
\end{tabular}
\end{center}
Again, these functions Poisson commute, hence $f$ is Lagrangian. The singular fibres of $f$ are lying over  $\Delta =\{ (0, r, 0) \mid r\in\R\}$. The reader may verify that the above gives a generic-singular fibration. 
\end{ex}

The reader should be aware that the above are just examples of Lagrangian positive and generic-singular fibrations. In fact, there are infinitely many germs of such fibrations \cite{RCB1}.

\subsection{The affine structures.}

Now we describe the integral affine structures induced by the above models by giving their period lattices explicitly. For the details we refer the reader to \cite{RCB1}.  Fibrations with generic-singular fibres can be normalized near $\Crit (f)$ according to the following: 

\begin{thm}\label{thm. normal form}  Let $\mathcal F= (X,\omega, f,B)$ be a generic-singular fibration. Assume that $\Sigma= \Crit (f)$ is non-degenerate. Then there is a $T^{2}$ invariant neighborhood $U\subseteq X$ of $\Sigma$ and a commutative diagram
\begin{equation}\label{diag:norm} 
\begin{CD}
U @>\Psi>> D^4\times D^{1}\times S^{1}\\
@Vf|_{U}VV  @VVFV\\
B @>\psi>> D^2\times D^{1}
\end{CD}
\end{equation}
where coordinates $(x,y)$ on $D^4$ and $(r,\theta)$ on $D^{1}\times S^{1}$ define standard symplectic coordinates, the map $\Psi$ is a symplectomorphism, $\psi$ is a diffeomorphism sending $\Delta$ to $\{0\}\times D^{1}$ and $F$ is given by (\ref{coord}). Furthermore $\Psi$ can be taken to be $T^{n-1}$ equivariant.
\end{thm}
The above is a corollary of a result due to Miranda and Zung \cite{Zung6}; we refer the reader to \cite{RCB1}\S 3 for the details. 
\begin{rem}
For convenience we shall assume that $B=f(U)$ where $U$ is as in Theorem \ref{thm. normal form}. We can think of the above normalization as providing $U$ with canonical coordinates and $B\cong D^2\times D^1$ with coordinates $b_1, b_2, b_3$ such that the Hamiltonian vector fields of $b_i\circ f|_U$ are linear. This linearization will be used to compute the action coordinates explicitly. This is crucial to understand the singularities of the affine structure in the base.
\end{rem}

\begin{prop}\label{prop. generic lattice}
Let $\mathcal F=(X,\omega ,f,B)$ be any generic-singular
fibration and $F_{\bar b}= f^{-1}(\bar b)$ a smooth fibre. 
There is a basis of $H_1(F_{\bar b},\mathbb{Z})$ whose corresponding basis
$\lambda_1, \lambda_2, \lambda_3$ of 
the period lattice $\Lambda$ of $\mathcal F$, in the coordinates
$b=(b_1,b_2,b_3)$ on $B\cong D^2\times D^1$ given by
Theorem \ref{thm. normal form}, can be written as
\begin{equation}\label{eq. generic periods}
\lambda_1=\lambda_0+dH,\qquad \lambda_2=2\pi db_2,\qquad \lambda_3=db_3,
\end{equation}
where $H\in C^\infty (B)$ is such that $H(0)=0$ and
$\lambda_0=-\log|b_1+ib_2| db_1+\Arg (b_1+ib_2)db_2$. The monodromy of $\Lambda$ is given by
\begin{equation}
\left( \begin{array}{ccc}
                 1 & 0 & 0 \\
                 1 & 1  & 0 \\
                 0 & 0  & 1 
              \end{array} \right) .
\end{equation}
\end{prop}
\begin{proof} 
The proof is the same as in \cite{RCB1} Proposition 3.10. Let $s=b_1+\sqrt{-1}b_2$ and $r_3=b_3$. Roughly speaking, one considers the maps given by $\sigma_1(s,r)=\left(\bar s\slash\epsilon, r,\theta_0\right)$ and $\sigma_2(s,r)=(\epsilon,s\slash\epsilon,r,\theta_0)$ for small $\epsilon >0$ and $\theta_0\in S^{1}$ fixed; these define sections of $f|_U=F$ disjoint from $\Crit(F)$, where $F$ is as in (\ref{coord}). The Hamiltonian vector fields $\eta_i$ of $F_i$ extend to $X\setminus U$. One can define a basis $\gamma$ of  $H_1(F_{\bar{b}},\mathbb{Z})$ in terms of suitable composition of the integral curves of $\eta_i$. The period $\lambda_1$ is obtained by integrating along the path $\gamma_1$ starting at $\sigma_1(s,r)$, passing through $\sigma_2(s,r)$ and going back to $\sigma_1(s,r)$. The contribution of $\gamma_1\cap U$ to the period $\lambda_1$ is $\lambda_0$, whereas the contribution of $\gamma_1\cap X\setminus U$ is $dH$. The remaining periods can be computed integrating along classes in $H_1(F_{\bar{b}},\Z)$ represented by integral curves of $\eta_2$ and $\eta_3$, respectively.
\end{proof}

As in the 2-dimensional focus-focus fibration, one can choose suitable branches of $\lambda_0$ and define action coordinates on these branches. One can easily verify that this defines a simple singular affine structure on $B$. We have:

\begin{cor}\label{cor:gen_simple}
A generic-singular fibration $\mathcal F=(X,\omega,f,B)$ induces a simple affine structure with singularities
on $B$. 
\end{cor}
\begin{proof} 
Consider the coordinates $(b_1, b_2, b_3)$ on $B=D^2\times D^1$ and the period lattice as in Proposition~\ref{prop. generic lattice}. With respect to these
coordinates $\Delta = \{ b_1 = b_2 = 0 \}$. Define open subsets of 
$B_0 = B - \Delta$:
\begin{eqnarray*}
       V_1 & = & B - \{ (b_1, 0 , b_3) \ | \ b_1 > 0 \}, \\
       V_2 & = & B - \{ (b_1, 0, b_3) \ | \ b_1 < 0 \}.
\end{eqnarray*}
On $V_j$ the action coordinates have the form
\[ A_j(b_1, b_2, b_3) = ( \psi_j(b_1, b_2) + H(b_1, b_2, b_3), 2 \pi b_2, b_3), \]
where $\psi_j$ is a choice of primitive of $\lambda_0$. Then $\mathscr A =\{U_j,A_j\}$ gives the integral affine structure on $B_0$.
As in the focus-focus case, for either $j=1,2$, the map $A_j$ extends to a homeomorphism, $A:B\rightarrow A(B)\subseteq\R^2\times\R$  such that $A(0)=0$. It is easy to show that, if $\tau(t) = H(0,0,t)$, then $A$ is an isomorphism between 
$(B,\Delta,\mathscr A)$ and a neighborhood of $\Delta_{\tau}$ in 
the affine manifold with singularities of Example~\ref{edge:affvar}.
\end{proof}

The case of Lagrangian fibrations of positive type is analogous.  Positive fibrations are locally modeled on the fibration in Example \ref{ex. HL} in a neighborhood of its critical locus. One can use this local description to compute the periods. We have (cf. \cite{RCB1}{Theorem 4.19}):

\begin{prop} \label{per:pos}
Let $\mathcal F=(X,\omega,f, B)$ be a Lagrangian fibration of positive type and $F_{\bar b}= f^{-1}(\bar b)$ a smooth fibre. 
Then there is a basis of $H_1(F_{\bar b},\mathbb{Z})$ and 
local coordinates $(b_1,b_2,b_3)$ on $B$ around $\bar b$, such that the corresponding period 1-forms are:
\begin{align} \label{pos:per}
\lambda_1 =\lambda_0 +dH ,\quad \lambda_2 = 2\pi db_2,\quad  \lambda_3 =
2\pi db_3 
\end{align}
where $H$ is a smooth function on $B$ such that $H(0)=0$ and $\lambda_0$ is multi-valued 1-form blowing up at $\Delta\subset B$, where 
\[
\Delta=\{ b_1=0, b_2=b_3\geq 0\}\cup\{ b_1=b_2=0,b_3\leq 0\}\cup\{ b_1=b_3=0, b_2\leq 0\}. 
\]
In the basis $\lambda_1,\lambda_2,\lambda_3$ of $\Lambda$ and for suitable generators of $\pi_1 (B-\Delta )$ satisfying $g_1g_2g_3=I$ (cf. Figure~\ref{loops}), the monodromy representation of $\mathcal F$ is
generated by the matrices: 
\begin{center}
$T_1 = \begin{pmatrix}
  1 & 0 & 0\\
  -1 & 1 & 0\\
  0 & 0 & 1
\end{pmatrix}$,\quad 
$T_2 = \begin{pmatrix}
  1 & 0 & 0\\
  0 & 1 & 0\\
  1 & 0 & 1
\end{pmatrix}$,\quad 
$T_3 = \begin{pmatrix}
  1 & 0 & 0\\
  1 & 1 & 0\\
  -1 & 0 & 1
\end{pmatrix}$.
\end{center}
\end{prop}

\medskip
We now prove that the affine structure on the base of a positive fibration is simple.

\begin{prop}\label{prop:pos_simple}
A Lagrangian fibration $\mathcal F=(X,\omega,f, B)$ of positive type induces on $B$ the structure of a simple affine manifold with singularities with positive vertex.
\end{prop}

\begin{proof} 
Let $(b_1, b_2, b_3)$ be the coordinates on $B$ and $\Delta\subseteq B$ as in Proposition~\ref{per:pos}. To avoid cumbersome notation let us assume $B = \R \times \R^2$. We may identify $\R^2$ with $\{ 0 \} \times \R^2$. Then $\Delta \subset \R^2$. Let $\lambda_1,\lambda_2,\lambda_3$ be the periods of $\mathcal F$ as in (\ref{pos:per}). 
We want to show that the affine structure on $B- \Delta$ induced by $\mathcal F$ is isomorphic to the one 
given in Examples~\ref{pos:aff} or \ref{pos:affvar}. To do this we will consider the locally defined map
$A = (A_1, A_2, A_3)$, where each $A_j$ is a suitable branch of a primitive of $\lambda_j$ such that $A_j(0)=0$. First we will show that --perhaps after replacing $B$ by a smaller neighborhood of $0$-- the map $A$ extends to a homeomorphism $A:B\rightarrow A(B)\subseteq\R^3$. Let
\begin{eqnarray*}
  R & = & \R \times \Delta           \\ 
  R^+ & = & \R_{\geq 0 } \times \Delta  \\ 
  R^- & = & \R_{\leq 0} \times \Delta
\end{eqnarray*}
and take the open cover $\{U_1, U_2 \}$ of $B-\Delta$ where
\begin{eqnarray}
 U_1 & = & B - R^{+},  \nonumber                  \\
U_2 & = & B - R^{-}. \label{cov:pos}  
\end{eqnarray}
On $U_1$ we can choose an affine coordinates map given by
\[ A(b_1, b_2, b_3) = ( \psi_1(b_1, b_2, b_3), 2\pi b_2, 2\pi b_3), \]
where $\psi_1$ is a primitive of $\lambda_1$. Clearly $A(R^-) \subset R$. 
We now show that $A$ extends continuously to $B$. The key observation 
is that the symplectic form $\omega$ is exact in a neighborhood of the 
singular fibre over the vertex of $\Delta$. This is straightforward in the case of 
Example~\ref{proper:pos}, where $\omega$ is the standard symplectic form on 
$\C^3$ but it is also true in general.
So assume $\omega = d\eta$ for some 1-form $\eta$. Now let us fix a basis $e = (e_1, e_2, e_3)$ of
$H^1(f^{-1}(U_1), \Z)$, corresponding to the periods $\lambda_1, \lambda_2$ and 
$\lambda_3$ respectively. Recall that action coordinates can be computed by
\[ A(b) = \left( - \int_{e_1(b)} \eta, \  - \int_{e_2(b)} \eta, \  
                                       - \int_{e_3(b)} \eta \right),\]
where $e_j(b)$ is a $1$-cycle, contained in $f^{-1}(b)$, representing $e_j$. 
We prove first that $A$, as a map, extends continuously to $B - \Delta$.
Notice that $e_2$ and $e_3$ are monodromy invariant, so we may assume that 
$e_2(b)$ and $e_3(b)$ are well defined for all $b \in B - \Delta$ and that
\begin{equation} \label{easy:coo}
  - \int_{e_j(b)} \eta = 2 \pi b_j, 
\end{equation}
for $j = 2, 3$. In particular, $A_2$ and $A_3$ are defined on $B$. Let us study
\[ \psi_1(b) = - \int_{e_1(b)} \eta. \]
Suppose that $\psi_1(\bar b) = 0$ for a fixed point $\bar b \in U_1$. Given another
point $b \in U_1$ let $\Gamma: [0,1] \rightarrow U_1$ be a path such that $\Gamma(0) = \bar b$ and $\Gamma(1) = b$. 
Consider the cylinder $S$ inside $f^{-1}(U_1)$ spanned by the cycles $e_1(\Gamma(t))$. Then 
one can see that 
\begin{equation} \label{psi:one}
 \psi_1(b) = \int_{S} \omega. 
\end{equation}
We may use (\ref{psi:one}) to define $\psi_1(b)$ for $b \in R^+ - \Delta$.  Since $B - \Delta$ is not simply connected, this expression of $\psi_1$ is well defined provided that it is independent of the chosen path $\Gamma$. Suppose that $\Gamma_1$ and $\Gamma_2$ are two different paths from $\bar b$ to $b$ such that $\Gamma_1 - \Gamma_2$ 
is not homotopically trivial in $B - \Delta$, then we have to show that if $S_1$ and $S_2$ are the 
corresponding cylinders, then
\[  \int_{S_1 - S_2} \omega = 0. \]
Denote by $e_1^+(b)$ and $e_1^-(b)$ those boundary components of $S_1$ and $S_2$ respectively,
which lie on top of $b$ (the endpoint of both $\Gamma_1$ and $\Gamma_2$).
Then
\[ \partial (S_1 - S_2) = e_1^+(b) - e_1^-(b), \]
and
\[ \int_{S_1 - S_2} \omega = \int_{e_1^+(b) - e_1^-(b)} \eta. \]
Because of monodromy, $e_1^+(b)$ and $e_1^-(b)$ may not coincide and it is not obvious 
that the above integral vanishes. Nevertheless, we know that $b \in R^+$ and there are three
cases: if $b = (b_1, b_2, b_3)$ then either $b_2 = 0$, $b_3 =0$ or $b_2 = b_3$. Let 
us look at that the latter case.
With respect to the basis $e=(e_1,e_2, e_3)$ as above, the monodromy matrices $T_1$, $T_2$ and $T_3$ corresponding respectively to generators $g_1$, $g_2$ and $g_3$ of $\pi_1(B - \Delta)$ as depicted in Figure~\ref{loops} are those given in Proposition~\ref{per:pos}.
\begin{figure}[ht]
\psfrag{g1}{$\Gamma_{1}$}
\psfrag{g2}{$\Gamma_{2}$}
\psfrag{b}{$\bar b$}
\psfrag{b1}{}
\psfrag{b2}{}
\psfrag{b3}{$b$}
\begin{center}
  \epsfig{file=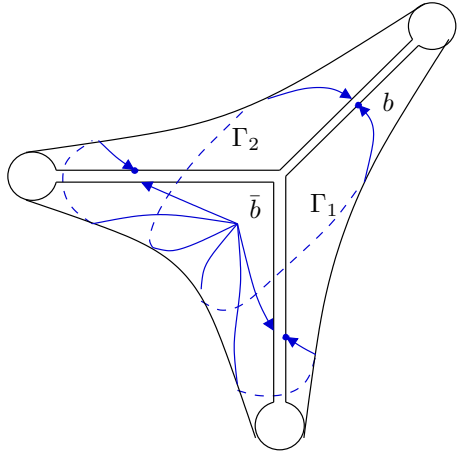, width=6cm,height=6cm}
\end{center}
\caption{The cut pair of pants are wrapping around $\Delta$ and give a schematic picture for 
         $U_1 = B - R^+$, the cut represents $R^+$.
         Here $b \in R^+$ and $\Gamma_1$ and $\Gamma_2$ are two possible paths from $\bar b$ to 
         $b$.}\label{pos_paths}
\end{figure}

Let $\bar b$, $b$, $\Gamma_1$ and $\Gamma_2$ 
be given as in Figure \ref{pos_paths}, then one can see that 
\todo{this should now be correct} $\Gamma_1 - \Gamma_2 = g_1^{-1} g_2^{-1}$. This implies that
\[ e_1^+(b) = e_1^{-}(b) - e_2(b) +  e_3(b) \]
and therefore that 
\[ \int_{e_1^+(b) - e_1^-(b)} \eta = \int_{- e_2(b) + e_3 (b)} \eta = 2 \pi (b_2 - b_3) = 0, \]
where in the second equality we have used (\ref{easy:coo}). Similarly one treats the cases 
$b_2 = 0$ or $b_3=0$ using monodromy matrices $T_1$ and $T_2$ respectively.
This shows that $\psi_1$ extends continuously to $B - \Delta$. It can be easily seen that it 
also extends continuously to points in $\Delta$. In fact one can use (\ref{psi:one}) as a definition
of $\psi_1(b)$ when $b \in \Delta$. This makes sense since the cycles $e_1(b)$ spanning $S$ can be extended as cycles on 
singular fibres when $b \in \Delta$, e.g. when $b = 0$,  $e_1(0)$ is a homologically non trivial closed 
curve passing through the singularity of $f^{-1}(0)$, in particular $e_1(0)$ is the generator of 
$H_1(f^{-1}(0), \Z) = \Z$. 

\medskip
We argue that $A$ is injective onto its image, at least when restricted to a smaller neighborhood of $b=0$. 
This would imply that $A$ is a homeomorphism. Clearly, $A$ is injective if and only if for 
fixed values of $b_2$ and $b_3$, the function $\psi_1(\, \cdot\, , b_2, b_3)$ is injective in a neighborhood
of $b= 0$. Since $d \psi_1 = \lambda_1$, this holds if the coefficient of $db_1$ in $\lambda_1$ is 
never zero in a neighborhood of $b = 0$. In fact, it was shown in \S 4 of \cite{RCB1} that this coefficient
blows up to infinity as $b\rightarrow 0$, in particular it never vanishes. 

\medskip
One can easily check that $A$ defines an isomorphism between the affine structure with 
singularities induced on $B$ by the fibration $\mathcal F$ and the one described in 
Example~\ref{pos:affvar}, where $\tau: \Delta \rightarrow \R$ is given by $\tau = \psi_1|_{\Delta}$.
We only need to verify that $\tau$ is smooth. In fact, it turns out that $\tau = H|_{\Delta}$ where $H$ is the smooth function in (\ref{pos:per}); this follows from the computation of $\lambda_0$ given in \cite{RCB1}\S4. Consider the fibration $F: \C^3 \rightarrow \R^3$ of Example \ref{ex. HL}. This is the local model for the singularity 
of a positive fibration. Consider two sections $\sigma_-$ and $\sigma_+$ of $F$, 
disjoint from $\Crit(F)$ and such that for every $b \in \Delta$, $\sigma_-(b)$ and $\sigma_+(b)$ lie on distinct 
connected components of the smooth part of the fibre over $b$. For every $b \in \R^3$ consider a curve $\gamma(b)$
contained $F^{-1}(b)$ joining $\sigma_-(b)$ to $\sigma_+(b)$ and define the function
\[ a_0(b) = - \int_{\gamma(b)} \eta. \]
Then $\lambda_0 = d a_0$. Clearly $a_0$ can be continuously defined on $\R^3$. Using the fact that $F$ satisfies 
$F(-z_1, z_2, z_3) = (-b_1, b_2, b_3)$, where $F(z_1, z_2, z_3) = (b_1, b_2, b_3)$, one can show that 
$a_0$ satisfies $a_0(-b_1, b_2, b_3) = - a_0(b_1, b_2, b_3)$ and therefore that $a_0|_{\Delta} = 0$. This 
proves that $\tau = H|_{\Delta}$.
\end{proof}

\subsection{Gluing over the discriminant locus}

Given a simple affine manifold with singularities, we show how to symplectically glue singular fibres of positive or generic type to the associated $T^3$ bundle. This gives us a (partial) symplectic compactification over positive and generic points of the singular locus. 

\medskip
Consider a cylinder $D^2\times I$ inside $\R^2\times\R$, where 
$I$ is an open interval, and let $\Delta =\{0\}\times I$. Let $H$ be a smooth 
real-valued function on $D^2\times I$. The \textit{germ} of $H$ along $\Delta$, denoted $H_\Delta$, is the Taylor expansion series of $H$ along $\Delta$. This is a formal power series in two variables whose coefficients are smooth functions on $I$. 

\begin{rem}\label{rem:power} For any given formal power series in two variables $h = \sum h_{ij} x_1^i x_2^j$ whose coefficients are smooth functions $h_{ij}=h_{ij}(r)$ on $I$, there is a function $H$ on $D^2\times I$ whose germ along $\Delta$ is $h$. An analogous statement in the case of a 
formal power series in one variable with real coefficients is standard (cf. \cite{Rudin} Exercise 
13, page 384). It is an exercise to check that it is also true in two variables with coefficients 
depending on a parameter.
\end{rem}

\medskip
Recall that the generators of the period lattice of a generic-singular fibration may be written as $\lambda_1=\lambda_0+dH$, $\lambda_2=2\pi db_2$ and $\lambda_3=db_3$, where $(b_1,b_2,b_3)$ are coordinates in in $D^2\times I$, $\lambda_0$ as (\ref{eq. generic periods}) and $H$ a smooth function. One can prove the following (cf. \cite{RCB1}):

\begin{thm}\label{thm. main RCB} For any smooth function $H$ over 
$B = D^2 \times I$, there is a generic-singular fibration $\mathcal F_{H} = (X, \omega, f, B)$ whose period lattice is generated by 1-forms as in (\ref{eq. generic periods}). Furthermore, two generic-singular fibrations $\mathcal F_{H}$ and $\mathcal F_{H'}$ are symplectically conjugate in a neighborhood of $\Delta$ if and only if $H_\Delta =H'_\Delta$.
\end{thm}

We call $H_{\Delta}$ the \textit{invariant} of the fibration $\mathcal F_H$. 
We proved in Corollary \ref{cor:gen_simple} that the affine base of a generic-singular fibration is always simple, isomorphic to Example \ref{edge:affvar}. Furthermore, the shape of its discriminant locus (in affine coordinates), as well as the isomorphism class of its singular affine base is determined by the function $\tau (r)=H(0,0,r)$  which is the restriction of $H$ to $\Delta$. In other words, by the zero order term of the germ $H_\Delta$. In the special case when the zero order term of $H_\Delta$ vanishes, the base is affine isomorphic to the  product of an affine disc with a node times the standard affine interval, in this case we call the associated fibration $\mathcal F_{H}$ \textit{straight}, in all other cases we call it \textit{twisted}.

\begin{lem}\label{lem. deform} 
Given any function $\tau\in C^\infty(\Delta)$ on an edge $\Delta \subset D^2\times I$ 
with $\tau (0)=0$, there is a generic-singular fibration whose base is locally affine isomorphic to the affine manifold with singularities $(\R^2\times I,\Delta_{\tau}, \mathscr A)$ of Example \ref{edge:affvar}.
\end{lem}
\begin{proof}
In view of Remark \ref{rem:power}, we can certainly find a smooth function $H$ on $D^2\times I$ such that $H|_\Delta=\tau$. We can then form $\mathcal F_H$ using Theorem~\ref{thm. main RCB}.
\end{proof}

\medskip
Analogously, positive fibrations are also classified by germs $H_\Delta$, where in this case $\Delta\subset D^3$ is a trivalent vertex and $H$ a smooth function on $D^3$ as in Proposition \ref{per:pos}; for the details we refer to \cite{RCB1}. Given a positive fibration, Proposition \ref{prop:pos_simple} tells us that its base is locally isomorphic to $(\R^3,\Delta_\tau,\mathscr A)$ as in Example \ref{pos:affvar}. A particular case is when $\tau=0$ which gives a straight vertex. More generally we showed (cf. proof of Proposition \ref{prop:pos_simple}) that $\tau = H|_{\Delta}$. In particular, we have:

\begin{lem}\label{lem. deform3} Given any function $\tau\in C^\infty(\Delta)$ on a trivalent vertex $\Delta\subset D^3\subset\R^3$ with $\tau (0)=0$, there is a positive fibration whose base is locally affine isomorphic to the affine manifold with singularities $(\R^3, \Delta_{\tau}, \mathscr A)$ 
of Example  \ref{pos:affvar}.
\end{lem}

We stress that the constructions described in Lemmas \ref{lem. deform} and \ref{lem. deform3} only 
involve the zero order term of $H_\Delta$, which is enough for determining the affine structure. From \cite{RCB1} it follows that we have many possible choices of $H_\Delta$ giving the same affine structure:

\begin{cor}  Given a prescribed affine manifold with singularities $(B, \Delta_{\tau}, \mathscr A)$ either as in Example \ref{edge:affvar} in the generic case or as in Example  \ref{pos:affvar} in the positive case, there are infinitely many non symplectically conjugate germs of Lagrangian fibrations whose 
bases are locally affine isomorphic to $(B, \Delta_\tau,\mathscr A)$.
\end{cor}
Observe that the above result holds also in the case when $\tau\equiv0$, i.e. when the discriminant is completely straight. Exploiting the flexibility given by Lemmas~\ref{lem. deform} and \ref{lem. deform3}, we can show that we can always
locally compactify a torus bundle given by simple affine manifolds with 
singularities near a positive or generic point of the discriminant locus: 
\begin{prop}\label{prop:aff_glue} Let $(B,\Delta, \mathscr A)$ be a given 
simple affine 3-manifold with singularities. 
Then we have the following
\begin{itemize}
\item[(i)] if $J \subseteq \Delta_g$ is an edge of $\Delta$, then 
there is a generic-singular fibration $\mathcal F$, with affine base 
$(B', \Delta', \mathscr A')$ and neighborhood $U \subseteq B$ of $J$
 such that there exists an integral affine isomorphism 
$(B',\Delta', \mathscr A') \cong (U, J, \mathscr A)$ inducing a symplectic conjugation $X(B'_0,\mathscr A') \cong X(U-J, \mathscr A)$;
\item[(ii)] if $p \in \Delta_d$ is a positive vertex of $\Delta$, then there 
is positive fibration $\mathcal F$ with base $(B', \Delta', \mathscr A')$ and a neighborhood $U\subseteq B$ of $p$ such that there exists an integral affine isomorphism 
$(B', \Delta', \mathscr A')\cong(U, U \cap \Delta,\mathscr A)$ 
inducing a symplectic conjugation $X(B'_0,\mathscr A') 
\cong X(U-(U \cap \Delta), \mathscr A)$.
\end{itemize}
Moreover, using the symplectic conjugations in (i) and (ii), we can 
symplectically glue the germ of $\mathcal F$ into $X(B_0, \mathscr A)$.
\end{prop}
\begin{proof} It is just a matter of applying Lemmas \ref{lem. deform} and \ref{lem. deform3} to find suitable $\mathcal F$. Since both positive and generic singular fibrations have a Lagrangian section, the result follows from Corollary 
\ref{aff:symp}.
\end{proof}

\subsection{Gluing legs}
While for the gluing in Proposition \ref{prop:aff_glue} it is sufficient to consider the zero order term of $H_{\Delta}$, to glue two singular Lagrangian 
fibrations $\mathcal F$ and $\mathcal F'$ along their legs one should take into account \emph{all} terms. This is essentially due to the fact that, gluing legs also involves gluing them along their singular fibres. We will see that 
Theorem \ref{thm. main RCB} also takes care of this.

\medskip
Suppose we are given a simple affine $3$-manifold with singularities
 $(B, \Delta, \mathscr A)$ and two points $p$ and $p'$ of $\Delta$ connected 
by an edge $J$ ($p$ and $p'$ may be generic, positive or 
negative points). Let us assume that we have glued to $X(B_0, \mathscr A)$ the germs of singular Lagrangian 
fibrations $\mathcal F$ and $\mathcal F'$ fibering over disjoint neighborhoods $V$ and $V'$ of $p$ and $p'$ 
respectively (e.g. using Proposition~\ref{prop:aff_glue}, if $p$ and $p'$ are positive 
or generic). We do not consider only the case when $p$ and $p'$ are either positive of generic, since we want 
the arguments here to hold also for negative points onto which we can glue  fibrations like the ones in \S\ref{negative}.  We only assume here that $\mathcal F$ and $\mathcal F'$ have legs with generic-singular fibres on their ends and these ends are connected by $J$. 
We now explain how to glue to $X(B_0, \mathscr A)$ a generic singular 
fibration along $J$ in such a way that this gluing is made compatible with the gluing of $\mathcal F$ and $\mathcal F'$.

\begin{figure}[!ht] 
\begin{center}
\input{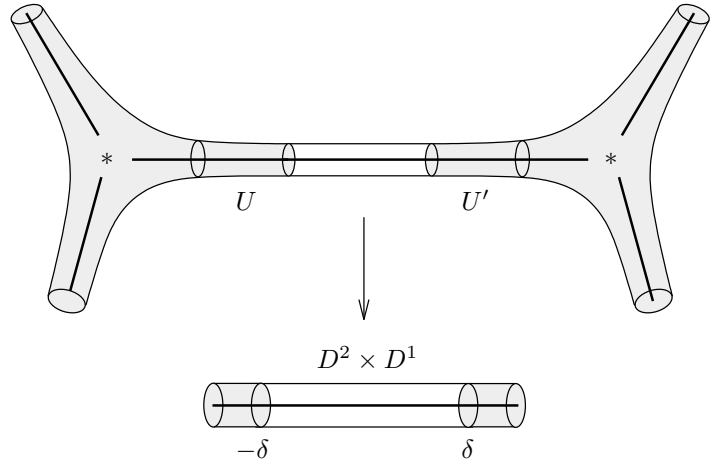}
\caption{The gluing of two legs along their ends. The asterisk represents components of the discriminant of $\mathcal F$ and $\mathcal F'$, which can be of either positive or negative type (or void).}\label{glue_leg}
\end{center}
\end{figure}

\medskip
We can assume that there are disjoint neighborhoods $U$ and $U'$ of the ends of 
$J$, as in Figure \ref{glue_leg}, and generic-singular fibrations $\mathcal L =\mathcal F|_U$  and $\mathcal L' = \mathcal F'|_{U'}$ over $U$ and $U'$. Let $H_{\Delta}$ and $H'_{\Delta}$ be, respectively, the invariants of $\mathcal L$ and 
$\mathcal L'$ as in Theorem \ref{thm. main RCB}.

\medskip
Since $J$ is an edge of $\Delta$, there is a neighborhood $W$ of $J$, with $W \cap \Delta = J$, such that 
$(W, J)$ is (locally) affine isomorphic to $(D^2 \times I, \Delta_{\tau})$ as in Example~\ref{edge:affvar}. Without loss of generality, we can assume $I = (-1, 1)$ and that 
there exists $\delta \in (0,1)$ such that $U\cong D^2 \times (-1, -\delta)$ and $U'\cong D^2 \times (\delta, 1)$.
Denote $I_{-\delta} = (-1, -\delta)$ and $I_{\delta} = (\delta, 1)$. 
Clearly, we can interpret $H_{\Delta}$ and $H'_{\Delta}$ as formal power series along $I_{-\delta}$ and 
$I_{\delta}$ respectively. By the arguments of the previous section, we must have that the zero order terms of 
$H_{\Delta}$ and $H'_{\Delta}$ coincide with $\tau|_{I_{-\delta}}$ and $\tau|_{I_{\delta}}$ respectively.

\medskip
It is now clear that we can choose a formal power series $\tilde{H}_{\Delta}$ along $I$ 
such that
\begin{itemize} 
 \item[(a)] the zero order term of $\tilde{H}_{\Delta}$ is $\tau$;
 \item[(b)]$\tilde{H}_{\Delta}$ coincides with $H_{\Delta}$ and $H'_{\Delta}$ along 
           $I_{-\delta}$ and $I_{\delta}$ respectively.
\end{itemize}
This can be done using cut-off functions. For this purpose it may be necessary to shrink $I_{-\delta}$
and $I_{\delta}$ by taking a slightly bigger $\delta$.

We can now apply Remark~\ref{rem:power} and the first part of Theorem~\ref{thm. main RCB} to find 
the germ of a generic-singular Lagrangian fibration $\tilde{\mathcal L}$ fibering over $W$ whose invariant 
is $\tilde{H}_{\Delta}$. The second part of Theorem~\ref{thm. main RCB} and condition $(b)$ above
imply that $\tilde{\mathcal L}|_{U} \cong \mathcal L$ and $\tilde{\mathcal L}|_{U'} \cong \mathcal L$, 
moreover condition $(a)$ implies that $\tilde{\mathcal L}$ can be glued to $X(B_0, \mathscr A)$ along $J$. 
It is clear that the symplectic conjugations $\tilde{\mathcal L}|_{U} \cong \mathcal L$ and 
$\tilde{\mathcal L}|_{U'} \cong \mathcal L'$ coincide with the map gluing $\tilde{\mathcal L}$
to $X(B_0, \mathscr A)$. 

\medskip
We have proved:
\begin{prop} \label{prop:def:leg} Let $(B, \Delta, \mathscr A)$ be a simple 
affine 3-manifold with singularities and let $p, p' \in \Delta$ be points 
connected by an edge $J$. Suppose there are disjoint neighborhoods $V$ and 
$V'$ of $p$ and $p'$ respectively and a neighborhood $W$ of $J$, with $W \cap \Delta = J$, 
such that the following conditions hold
\begin{itemize}
\item[(i)] if $\tilde B = B_0 \cup (V \cup V')$, there exists a Lagrangian fibration 
$\mathcal F = (X, \omega, f, \tilde B)$ and a commuting diagram
\begin{equation*} 
\begin{CD}
X(B_0, \mathscr A)  @>\Psi>> X\\
@Vf_0VV  @VVfV\\
B_0 @>\iota>> \tilde B
\end{CD}
\end{equation*}
where $\Psi$ is a symplectomorphism and $\iota$ the inclusion. 
\item[(ii)] $\mathcal F|_{W \cap V}$ and $\mathcal F|_{W \cap V'}$ are generic-singular 
            fibrations.
\end{itemize} 
Then, if we let $\tilde B' = \tilde B \cup W$, there exists a Lagrangian fibration 
$\mathcal F' = (X', \omega', f', \tilde B')$ and a commuting diagram
\begin{equation*} 
\begin{CD}
X(B_0, \mathscr A)  @>\Psi'>> X'\\
@Vf_0VV  @VVf'V\\
 B_0 @>\iota >> \tilde B'
\end{CD}
\end{equation*}
where $\Psi'$ is also a symplectomorphism.
\end{prop}

\medskip
The upshot of the results of this Section is that: 1) we can construct local models of generic and positive singular fibres; 2) we know how to glue them onto any given simple affine manifold with generic and positive singularities; 3) these gluings can be made compatible over common intersections. In fact, we can show:

\begin{thm}\label{thm:big_positive}
Let $(B,\Delta, \mathscr A)$ be a compact simple integral affine 3-manifold with singularities 
without negative vertices.
Then there is a compact smooth symplectic 6-manifold $(X,\omega)$ and a $C^\infty$ Lagrangian fibration 
$f:X\rightarrow B$ with discriminant locus $\Delta$, which is a semi-stable compactification of the $T^3$ bundle $X(B_0,\mathscr A)\rightarrow B_0$.
\end{thm}

The proof is an application of the above preparation results. Using Proposition~\ref{prop:aff_glue} we 
can first glue in the positive vertices, then using Proposition~\ref{prop:def:leg} we glue in the generic-singular 
fibres over the edges. Theorem \ref{thm:big_positive} is a particular case of our more general result we shall prove 
in \S\ref{section:compact}, where we also include negative fibrations. 
We emphasize that the fibration obtained in Theorem \ref{thm:big_positive} is \emph{smooth}. This will not happen if $\Delta$ includes negative vertices. In that case, the resulting fibration will be \emph{piecewise} smooth only.

\medskip
As a further remark we point out that Theorem \ref{thm:big_positive} can be generalized to dimension 
$n\geq 3$, since there are natural generalizations of generic and positive singularities and the analysis of their affine structures carries through as in the $n=3$ case. Our notion of simplicity can also be generalized to higher dimensions, though for $n>3$ it may no longer coincide with the notion of simplicity in the sense of Gross and Siebert \cite{G-Siebert2003}.

\section{Piecewise smooth fibrations}\label{sec:pwfibr}
It is now commonly accepted that to produce Lagrangian fibrations of the type described in \S\ref{sect. topology} one should also allow \emph{piecewise smooth} fibrations  (cf. \cite{Gross_spLagEx}, \cite{Joyce-SYZ}, \cite{Ruan3}).  Here we present a simple way to produce local models of piecewise smooth Lagrangian fibrations. We suspect that models of the sort presented here are also implicit in Ruan's 
fibrations  but we have been unable to verify this. Our method is inspired by ideas of Gross \cite{Gross_spLagEx}, Goldstein 
\cite{Goldstein1} and Joyce \cite{Joyce-SYZ}.
\subsection{Fibrations with torus symmetry.}
Let $(X,\omega )$ be a symplectic $2n$-manifold and let  
$\mu :(X,\omega)\rightarrow\mathfrak t^\ast$ be the moment map of a 
Hamiltonian $T^k$-action. Let $t\in\mu (X)$ and let 
$\pi_t:\mu^{-1}(t)\rightarrow X_t$ be the projection modulo the $T^k$ 
action. When $t$ is a regular value of $\mu$,  $X_t$ is a smooth
manifold and the symplectic form $\omega$ descends to a symplectic form 
$\omega_t$ on $X_t$. When $t$ is a critical value of $\mu$, $X_t$ may be a singular 
space and $\omega_t$ will be only defined on the smooth part of $X_t$. 
The  space $(X_t,\omega_t)$ is the Marsden-Weinstein reduced 
space at $t$. 

\begin{rem} We shall denote by 
\[ \omega_{\C^m}=\frac{i}{2}\sum_{k} dz_k \wedge d\zbar_k \]
the standard symplectic structure on $\C^m$ and $\omega_0$ will denote the reduced symplectic form of the reduced space $(X_t,\omega_t)$ at time $t=0$.
\end{rem}

\medskip
Goldstein \cite{Goldstein1} and Gross \cite{Gross_spLagEx} used 
reduced spaces to construct  $T^k$-invariant (special) Lagrangian fibrations. 
The following is a particular case of \cite{Gross_spLagEx}{Thm. 1.2}:
\begin{prop} \label{prop. lifting}
Let $T^k$ act effectively on $X$, $k\leq n-1$. Suppose that there is a 
continuous map $G:X\rightarrow M$ to an $(n-k)$-dimensional manifold $M$ such 
that $G(T\cdot x)=G(x)$ for all $T\in T^k$. Suppose that for $t$ in a dense subset of $\mu (X)$ the induced maps 
$G_t:X_t\rightarrow M$ have fibres that are Lagrangian with respect to 
$\omega_t$. Then $f:X\rightarrow \mu (X)\times M$ given by:
\begin{equation} \label{eq. lifted-fibr}
f=(\mu , G)
\end{equation}
defines a $T^k$-invariant Lagrangian  fibration. 
\end{prop}
When the $T^k$-action has fixed points, the construction of Proposition 
\ref{prop. lifting} will produce
fibrations with interesting singular fibres. We will 
give some explicit examples shortly.

\begin{rem} \label{onedim:red}
In the extremal case when $k=n-1$, constructing Lagrangian 
fibrations using Proposition~\ref{prop. lifting} is very easy. In 
this situation, the reduced spaces $X_t$ are two dimensional and every 
map $G_t: X_t \rightarrow \R$ with $1$-dimensional level sets defines a 
Lagrangian fibration on $X_t$. In particular, any $T^{n-1}$-invariant
continuous map $G :X\rightarrow \R$ which, on each $X_t$, descends to a map $G_t$  
with $1$-dimensional level sets can be used to construct Lagrangian fibrations. 
We will make much use of this fact later on.
\end{rem}

\subsection{The reduced geometry.}
Consider the following $S^{1}$ action on $\C^n$, with $n \geq 2$:
\begin{equation} \label{action}
e^{i\theta}(z_1,z_2, z_3, \ldots, z_n)
          =(e^{i\theta}z_1,e^{-i\theta}z_2,z_3, \ldots, z_n).
\end{equation}
This action is Hamiltonian with respect to $\omega_{\C^n}$. Clearly it
is singular along the $2(n-2)$ dimensional symplectic submanifold 
$\Crit (\mu)= \{ z_1=z_2=0 \}$. 
The moment map is:
\begin{equation}\label{eq mu}
\mu (z_1, \ldots, z_n)=\frac{|z_1|^2-|z_2|^2}{2}.
\end{equation}
The only critical value of $\mu$ is $t=0$ and 
$\Crit (\mu) \subset \mu^{-1}(0)$. 

\medskip
Now consider the map $\bar\pi$ as in Remark \ref{rem. pi local mod}.  Recall that $\bar\pi$ is given by
\begin{equation}\label{barpi}
\begin{array}{rll} 
 \bar\pi:\C^n &\rightarrow &\R \times \C^{n-1} \\
 (z_1, \ldots, z_n) & \mapsto & (\mu, z_1 z_2, z_3, \ldots, z_n).
\end{array}
\end{equation}

When restricted to $\C^n - \Crit (\mu)$, the above is an $S^1$-bundle onto 
$(\R \times \C^{n-1})- \bar\pi(\Crit (\mu))$ with Chern class $c_1 = 1$.
Let $\pi_t$ be the restriction to $\mu^{-1}(t)$ of the map
\begin{equation}\label{pi_t}
(z_1, \ldots, z_n) \mapsto (z_1z_2, z_3, \ldots, z_n).
\end{equation}
Then $\pi_t$ can be used to identify the reduced space $\mu^{-1}(t) / S^1$ with 
$\C^{n-1}$. Under this identification, i.e. letting the coordinates
$u_1 = z_1 z_2$ and $u_j = z_{j+1}$ when $2 \leq j \leq n-1$, the reduced symplectic form
$\omega_t$ can be written as:
\begin{equation} \label{om:red}
  \omega_t = \frac{i}{2} \left( \frac{1}{2 \sqrt{ t^{2} + |u_{1}|^{2}}} 
           \,  du_{1} \wedge d\overline{u}_{1} +  \sum_{j=2}^{n-1} \, du_{j} 
                  \wedge d\overline{u}_{j} \right).
\end{equation}  
Clearly, away from $t=0$, the reduced spaces are smooth manifolds. 

\medskip
On the other hand, at $t=0$ the reduced form $\omega_0$ blows up
along the hyperplane 
\[
\Sigma:=\pi_0(\Crit (\mu))=\{u_1=0\},
\] 
so the reduced space $(\C^{n-1},\omega_0)$
is singular. However, it was observed by Guillemin and Sternberg 
in \cite{Gui-Stern-bi}, that it can be smoothed out, i.e. it can be identified
with $(\C^{n-1},\omega_{\C^{n-1}})$. 
Indeed, the identification is given by the following 
\begin{equation}\label{eq Gamma_0}
\Gamma_0:(u_1,u_2, \ldots, u_{n-1}) \mapsto 
          \left ( \frac{u_1}{\sqrt{|u_1|}},u_2, \ldots, u_{n-1} \right ).
\end{equation}
The map $\Gamma_0$ is continuous, smooth away from $u_1=0$ and 
such that $\Gamma_0^\ast\omega_{\C^{n-1}}=\omega_0$. 
One can do more: one can identify all the reduced spaces with 
$(\C^{n-1}, \omega_{\C^{n-1}})$ at once. Consider the map
\begin{equation}\label{eq. Gamma_t}
\Gamma_t: (u_1, u_2, \ldots, u_{n-1}) \mapsto 
          \left( \frac{u_1}{\sqrt{|t|+\sqrt{t^2+|u_1|^2}}}, u_2, \ldots, u_{n-1} \right).
\end{equation}
One can verify that $\Gamma_t$ is a symplectomorphism between 
$(\C^{n-1}, \omega_t)$ and the standard symplectic space $\C^{n-1}$. However, this identification has the problem that, although continuous and smooth for fixed $t\in\R$, it is not smooth in $t$ when $t=0$. In fact
one can show that it cannot be otherwise.

\subsection{A construction}
We now illustrate a general method to construct piecewise smooth Lagrangian
fibrations using Proposition \ref{prop. lifting} and the observations about
the reduced geometry with respect to the $S^1$ action as in (\ref{action}).

\medskip 
Let $\Log :(\C^\ast )^{n-1}\rightarrow\R^{n-1}$ be the map defined by 
\begin{equation}\label{eq log fibr.}
\Log (v_1,\ldots, v_{n-1})=(\log |v_1|, \ldots, \log |v_{n-1}|).
\end{equation}
Clearly, the above map is a Lagrangian fibration with respect to the 
restriction of $\omega_{\C^{n-1}}$ to $(\C^\ast )^{n-1}$. Moreover, it defines a trivial 
$T^{n-1}$-bundle over $\R^{n-1}$. Let the map
\[
\Phi: \C^{n-1} \rightarrow \C^{n-1}
\]
be a smooth symplectomorphism of the standard $\C^{n-1}$. Let $X_t$ be the open and dense subsets 
of $(\C^{n-1}, \omega_t)$ defined by
\[ X_t = \Gamma_t^{-1} \circ \Phi^{-1} \left( (\C^{\ast})^{n-1} \right). \] 
Denote, with slight abuse of notation, 
\[
\Sigma :=\{u_1=0\}\cap X_0.
\]
Then examples of maps $G_t : X_t \rightarrow \R^{n-1}$ as in Proposition 
\ref{prop. lifting} can be defined by
\[ G_t = \Log \circ \Phi \circ \Gamma_t. \]
This clearly makes sense also when $t=0$. It is also clear that, for
all fixed $t \in \R$, $G_t$ is a Lagrangian fibration with respect to the
reduced symplectic form (\ref{om:red}).  We summarize this in the following:

\begin{prop}\label{prop. piecewise smoothness}
Let $\Phi$, $X_t$ and $G_t$ be as defined above. Let $Q$ be the map given by
\begin{equation} \label{qiu}
 Q(t,u_1, \ldots, u_{n-1}) = (t, G_t(u_1, \ldots, u_{n-1})). 
\end{equation}
Then $Q$ is defined on the dense open subset $Y\subseteq\R \times \C^{n-1}$
defined by
\[ Y = \{ (t,u_1, \ldots, u_{n-1}) \in \R \times \C^{n-1} \ | \ 
                                   (u_1, \ldots, u_{n-1}) \in X_t \}. \]
Letting $\bar \pi$ be as in (\ref{barpi}) and  
\[ X = (\bar \pi)^{-1}(Y) \]
with the standard symplectic form induced from $\C^n$, the map
$f: X \rightarrow \R^n$ given by
\[ f =  Q \circ \bar\pi  \]
is a piecewise smooth Lagrangian fibration of $X$ which fails to
be smooth on the $(2n-1)$-dimensional subspace $\mu^{-1}(0) \cap X$.
\end{prop}
It is clear that all the singular fibres of $f$ must lie in $\mu^{-1}(0) \cap X$. 
In fact, the singular fibres are all the lifts of fibres of $G_0$ in $X_0$ 
which intersect $\Sigma$. The topology of the
singularity depends on the topology of this intersection. The discriminant
locus of the fibration is therefore the set $\Delta \subset \R^{n}$ given by
\[ \Delta = \{ 0 \} \times \left(\Log \circ \Phi \circ \Gamma_0 (\Sigma)\right).\]
Given a point $b = (0, b_1, \ldots, b_{n-1}) \in \Delta$, 
the fibre $f^{-1}(b)$ looks like $S^1 \times G_{0}^{-1}(b_1, \ldots, b_{n-1})$ 
after the circles over all points in 
$G_0^{-1}(b_1, \ldots, b_{n-1}) \cap \Sigma$ have been collapsed to points 
(cf. Figure \ref{pants}).

\subsection{Examples}
In the following examples we use the above construction with $n=2$ or $3$.
Define the piecewise smooth map $\gamma: \C^2 \rightarrow \C$ by
\begin{equation}\label{eq. g}
\gamma (z_1,z_2)= \begin{cases}
\frac{z_1z_2}{|z_1|},\quad\text{when}\ \mu(z_1, z_2)\geq 0\\
\\
\frac{z_1z_2}{|z_2|},\quad\text{when}\ \mu(z_1, z_2) <0.
\end{cases}
\end{equation}
If $\pi_t$ is the restriction of the map (\ref{pi_t}) to $\mu^{-1}(t)$, then one can easily see that for all $(z_1, z_2, z_3) \in \mu^{-1}(t)$,
the map $\Gamma_t \circ \pi_t$ is given by
\[ \Gamma_t \circ \pi_t: (z_1, z_2, z_3) \mapsto (\gamma(z_1z_2), z_3). \]

From  Proposition \ref{prop. piecewise smoothness}, we see that $\Gamma_t \circ \pi_t$ can be twisted by a symplectomorphism $\Phi$. The topology of the resulting fibration depends on how we choose $\Phi$.

\begin{ex}[The amoeba] \label{ex amoebous fibr}
Take as a symplectomorphism $\Phi$ the linear map
\begin{equation}\label{eq. linear sympl}
\Phi (u_1,u_2)=\frac{1}{\sqrt{2}}\left ( u_1-u_2,u_1+u_2-\sqrt{2}\right ).
\end{equation}
Then the fibration resulting from Proposition 
\ref{prop. piecewise smoothness} can be written explicitly in the 
coordinates of the total space. We obtain:
\begin{equation} \label{eq. the fibration}
f(z_1, z_2, z_3) = \left(\frac{1}{2}\left(|z_1|^2-|z_2|^2\right), \log \frac{1}{\sqrt{2}}\left| \gamma - z_3\right|, 
                  \log \frac{1}{\sqrt{2}} \left|\gamma + z_3 - \sqrt{2}\right|\right),
\end{equation}
where $\gamma$ is as in (\ref{eq. g}). It is not difficult to see that $\Phi \circ \Gamma_0$ sends $\Sigma$ to the surface in 
$(\C^\ast)^2$ given by 
\[ \Sigma^{\prime} = \{ v_1 + v_2 + 1 = 0 \}, \]
which is, topologically, a pair of pants. Then the discriminant locus is 
\[ \Delta = \{ 0 \} \times \Log ( \Sigma^{\prime}), \]
which has the shape in Figure \ref{fig: amoeba}. This example
is topologically conjugate to the one in Example~\ref{ex. alt (2,1)}, before
the surface $\Sigma^{\prime}$ has been twisted. For the discussion of
the topology of the fibres in this example we refer to 
Example~\ref{ex. alt (2,1)}.
\end{ex}

In dimension $n=2$ we have the following:
\begin{ex} [Stitched focus-focus] \label{broken focus focus}
Using Proposition~\ref{prop. piecewise smoothness} we can obtain the following 
piecewise smooth fibration:
\begin{equation} \label{fib ff}
f(z_1,z_2)=\left( \frac{|z_1|^2-|z_2|^2}{2}, \, \log |\gamma (z_1,z_2)+1| \right).
\end{equation}
It is clearly well defined on 
$X = \{ (z_1, z_2) \in \C^2 \ | \ \gamma(z_1, z_2) + 1 \neq 0 \}.$
Observe that $f$ is topologically conjugate to a focus-focus fibration,
 hence to Example~\ref{ex top ff}. 
The only singular fibre is $f^{-1}(0)$ and it is a pinched torus. The fibration 
fails to be smooth on $\mu^{-1}(0)$. This example consists of the union of two smooth Lagrangian fibrations meeting along the ``stitch", $\mu^{-1}(0)$. We study this kind of piecewise smoothness in detail in \S \ref{stitched fibr}.
\end{ex}

Notice that in this example we are in the extremal case of Proposition~\ref{prop. lifting}, i.e. the reduced spaces are 2-dimensional and Remark~\ref{onedim:red} applies. 
In particular, the second component of 
$f$ in (\ref{fib ff}) could be replaced by any $T^2$ invariant function $G$, i.e. depending on $t =  \frac{1}{2}\left(|z_1|^2-|z_2|^2\right)$ and $u_1 = z_1 z_2$, subject to the condition that all maps $G_t$ have
$1$-dimensional level sets. 
Using this idea it is easy to construct everywhere smooth fibrations, such as the one in Example~\ref{smooth:ff} where $G(t, u_1) = \log |u_1 + 1|$.
Of course the topology of the resulting fibration depends on the topology of 
the maps $G_t$.

\medskip
We have an analogous model in dimension $n=3$:

\begin{ex}[The leg] \label{leg}
Consider the following affine symplectomorphism of $(\C^2, \omega_{\C^2})$
\begin{equation} \label{sympl leg}
 \Phi: (u_1 , u_2) \mapsto ( -u_2, u_1 -1). 
\end{equation}
The surface $\Sigma$ is sent by $\Phi \circ \Gamma_0$ to 
                $\Sigma^{\prime} = \{ v_2 + 1 = 0\}$.
The amoeba of $\Sigma^{\prime}$ is just a straight line. 
The resulting fibration $f$ is
\begin{equation}\label{eq leg fibr} 
f(z_1, z_2, z_3) = \left( \frac{|z_1|^2-|z_2|^2}{2}, \, \log |z_3|, \, 
             \log |\gamma (z_1,z_2)-1| \right).
\end{equation}
The discriminant locus is $\{ 0 \} \times \R \times \{0 \} \subset \R^3$, a 
horizontal line in the plane $\{ 0 \} \times \R^2$. The fibration is a 
piecewise smooth version of the generic singular fibration in 
Example~\ref{smooth:generic}.
Notice that this fibration is invariant under the Hamiltonian $T^2$-action
\begin{equation} \label{t2:action}
(e^{i\theta_1}, e^{i \theta_2}) \cdot (z_1,z_2, z_3) = 
     (e^{i\theta_1}z_1,\, e^{-i\theta_1}z_2,\, e^{2i \theta_2} z_3),
\end{equation}
whose moment map is
\[ (z_1, z_2, z_3) \mapsto \left( \frac{|z_1|^2-|z_2|^2}{2}, \, |z_3|^2 \right). \]

There are other choices of symplectomorphisms
$\Phi$ giving piecewise smooth generic fibrations. 
Although not very different from the  previous one, we will write other two
for convenience, since we will need them in the next example. 
The first one is
\begin{equation} \label{sympl leg 2}
 \Phi: (u_1 , u_2) \mapsto ( u_1 - 1, u_2 - \sqrt{2}). 
\end{equation}
It gives the fibration
\begin{equation}\label{eq leg fibr 2} 
f(z_1, z_2, z_3) = \left( \frac{|z_1|^2-|z_2|^2}{2}, 
             \, \log |\gamma (z_1,z_2)-1|, \, \log \left|z_3 - \sqrt{2}\right|
              \right),
\end{equation}
whose discriminant locus is the vertical line $\{0 \} \times \{ 0 \} \times \R \subset \R^3$. Also in this case it is clearly invariant under a $T^2$ action. The last choice of $\Phi$ is
\begin{equation} \label{sympl leg 3}
 \Phi: (u_1, u_2) \mapsto 
                  \frac{1}{\sqrt{2}} ( u_1-u_2, \, u_1+u_2 ), 
\end{equation}
giving
\begin{equation}\label{eq leg fibr 3} 
f(z_1, z_2, z_3) = \left( \frac{|z_1|^2-|z_2|^2}{2}, 
             \, \log |\gamma (z_1,z_2)- z_3|, \, \log |\gamma (z_1,z_2)+ z_3|
              \right),
\end{equation}
whose discriminant is the slope +1 diagonal through zero in $\{ 0 \} \times \R^2$. The $T^2$ action
in this case is given by
\begin{equation} \label{t2:diag}
(e^{i\theta_1}, e^{i \theta_2}) \cdot (z_1,z_2, z_3) = 
     (e^{i(\theta_2 + \theta_1)}z_1,\, e^{i(\theta_2 - \theta_1)}z_2,
                                                \, e^{2i \theta_2} z_3),
\end{equation} 
whose moment map is
\[ (z_1, z_2, z_3) \mapsto \left( \frac{|z_1|^2-|z_2|^2}{2}, \,
 \frac{|z_1|^2+|z_2|^2}{2}+|z_3|^2 \right). \]
In the above examples, the reduced spaces are all 2-dimensional. Using 
Remark~\ref{onedim:red} we can construct variations of (\ref{eq leg fibr}) by replacing the last component of 
(\ref{eq leg fibr}) with any function depending on  $t= \frac{|z_1|^2-|z_2|^2}{2}$, $s = |z_3|^2$ and
$u_1 = z_1z_2$,  subject to the condition that all the maps $G_t$ 
have $1$-dimensional level sets. A choice providing an example of a smooth fibration 
is given by $G = \log|u_1 - 1|$, which gives us 
Example~\ref{smooth:generic}. One can do more. In fact, one can take a function 
$G$ which gives an interpolation between the piecewise smooth fibration in (\ref{eq leg fibr}) and the smooth one in Example~\ref{smooth:generic}. This can be done by taking $G$ depending also on $s$, such that $G$ 
is equal to $\log |\gamma (z_1,z_2)-1|$ when $s$ is big and equal to 
$\log|u_1 -  1|$ when $s$ is small. We will say more about this later 
on, as this idea is useful in an important step of the main construction 
of the paper.
\end{ex}

\begin{ex}[The amoeba with thin legs] \label{thin leg}
We now construct an example which interpolates Example \ref{ex amoebous fibr}
and \ref{leg}. Consider the smooth function:
\[ H_0 = \frac{\pi}{4} \im (u_1 \overline{u}_2) \]
and let $\eta_{H_0}$ be  the Hamiltonian vector field associated to $H_0$.
If $\Phi_s$ is the flow generated by $\eta_{H_0}$, then the Hamiltonian
symplectomorphism associated to $H_0$ is defined to be $\Phi_{H_0} = \Phi_{1}$.
One computes that in our case
\[ \Phi_{H_0}: (u_1, u_2) \mapsto \frac{1}{\sqrt{2}} (u_1 - u_2, u_1 + u_2). \]
It maps $ \{ u_1 = 0 \}$ to $\{ v_1 + v_2 = 0 \}$.
We now want a symplectomorphism which acts like $\Phi_{H_0}$ in a 
small ball centered at the origin and like the identity outside a
slightly bigger ball. So choose a cut-off function 
$k: \R_{\geq 0} \rightarrow [0,1]$ such that, for some $\epsilon > 0$, 
\begin{equation} \label{gl:fn}
    k(t) = \left \{ \begin{array}{l} 
                   1 \ \quad\text{when} \ 0< t \leq \epsilon; \\
                   0 \ \quad\text{when} \ t \geq 2\epsilon
                  \end{array} \right.
\end{equation}
 and define the Hamiltonian
\[ H = k(|u_1|^2 + |u_2|^2) H_0. \]
The Hamiltonian symplectomorphism $\Phi_H$ associated to $H$ satisfies
\[ 
\Phi_H(u_1, u_2 ) = 
\begin{cases} 
 \id_{\C^2},&\quad\text{when}\ |u_1|^2 + |u_2|^2 \geq 2\epsilon;\\
\\
\frac{1}{\sqrt{2}} (u_1 - u_2, u_1 + u_2),&\quad\text{when}\ 
                      |u_1|^2 + |u_2|^2 \leq \epsilon .
\end{cases}
\]
Now let $\Psi$ be the affine symplectomorphism
 \[
\Psi: (v_1,v_2) \mapsto \frac{1}{\sqrt{2}}(v_1-v_2,v_1+v_2-\sqrt{2}).
\]
and finally, define $\Phi = \Psi \circ \Phi_H$. It is clear that 
\[ 
\Phi(u_1, u_2 ) = 
\begin{cases} 
 \Psi,&\quad\text{when}\ |u_1|^2 + |u_2|^2 \geq 2\epsilon;\\
\\
  ( - u_2, u_1 - 1),&\quad\text{when}\ 
                      |u_1|^2 + |u_2|^2 \leq \epsilon .
\end{cases}
\]
Notice that $\Phi$ acts like in (\ref{sympl leg}) on the 
ball of radius $\sqrt{\epsilon}$ around the origin, i.e. in 
a neighborhood of the surface $\{u_1=0\}$, and
like in (\ref{eq. linear sympl}) outside a larger ball. 
We use this $\Phi$ to construct a fibration $f$ using Proposition 
\ref{prop. piecewise smoothness}.
One can then see that $\Phi \circ \Gamma_0$ sends $\Sigma$ to a surface $\Sigma'$ 
such that $\Log (\Sigma^{\prime})\subset\R^2$ is a 3-legged 
amoeba with the end of the horizontal leg pinched down to a straight line. 
The discriminant locus of 
$f$ is then $\Delta=\{ 0 \} \times \Log (\Sigma^{\prime}) \subset \R^3$.
Of course, $f$ fails to be smooth on the slice $\mu^{-1}(0)$.
Using the same method we can twist $\Sigma$ 
suitably and obtain a fibrations having discriminant locus an amoeba with 
three thin legs (cf. Figure \ref{interpol}). For example, to pinch the 
diagonal leg to a thin line, choose a smooth function $H_0$ generating the Hamiltonian 
symplectomorphism
\[ (u_1, u_2) \rightarrow \left( u_1 + \frac{1}{\sqrt{2}}, u_2 + 
                                                 \frac{1}{\sqrt{2}} \right).
                                                                       \]
Cut $H_0$ off with a function $\rho$ which vanishes when 
$|u_2|^2 \leq M/2$, for some big $M$, and is equal to $1$ when
$|u_2|^2 \geq M$. This produces a Hamiltonian $H$. Now one proceeds as 
before. With an almost identical procedure one pinches down the vertical leg. 
The final choice of symplectomorphism $\Phi$ pinching down all three legs simultaneously may look like:
\begin{equation} \label{thin:fi} 
\Phi(u_1, u_2 ) = 
\begin{cases} 
  ( - u_2, u_1 - 1),&\quad\text{when}\ 
                      |u_1|^2 + |u_2|^2 \leq \epsilon; \\
  \\
  (u_1 - 1, u_2 - \sqrt{2}), &\quad\text{when}\ 
                      |u_1|^2 + |u_2 - \sqrt{2}|^2 \leq \epsilon; \\
  \\  
  \frac{1}{\sqrt{2}} ( u_1-u_2, \, u_1+u_2 ), &\quad\text{when}\ 
                       |u_2|^2 \geq M; \\
  \\
  \Psi,&\quad\text{everywhere else}.
\end{cases}
\end{equation}
It is clear that this piecewise smooth example is topologically conjugate to the one in Example~\ref{ex. alt (2,1)}. Here we have made explicit the twistings described there. In \S 7 we will show that this fibration can be modified so that it is actually smooth towards the ends of the three legs. For this we will develop further the smoothing method sketched at the end of Example~\ref{leg}. Also in \S 7, we will show that this fibration can be modified so that it is smooth away from a neighborhood homeomorphic to a 2-disk containing the codimension 1 part of its discriminant.
\end{ex}

The next result states existence of Lagrangian sections of the fibrations in the previous examples.

\begin{prop} \label{lag:section:exist}
 The fibrations in Example~\ref{ex amoebous fibr} and \ref{thin leg} have 
smooth Lagrangian sections which do not intersect the critical surface $\Crit(f)$. 
\end{prop}
\begin{proof}
Consider the symplectomorphism $\Phi$ from Example~\ref{ex amoebous fibr}. The reduced fibration 
at time $t=0$, i.e. the map $G_0 = \Log \circ \Phi \circ \Gamma_0$, has many Lagrangian sections,
since the $\Log$ fibration has many.  In particular we can choose one which does not intersect 
$\Sigma = \Crit(f)$, this follows for example by observing that the following Lagrangian section 
of the $\Log$ fibration 
\begin{equation} \label{log:section}
 (x_1, x_2) \mapsto (i e^{x_1}, e^{x_2}) 
\end{equation}
does not intersect the surface $\Sigma' = \{ v_1 + v_2 + 1 \}$. 
It is easy to see that a section which does not intersect $\Sigma$ can be lifted 
to $\mu^{-1}(0)$. The image of this lift is a coisotropic $2$ dimensional submanifold of $X$. Applying
the coisotropic embedding theorem, we can extend this submanifold to a Lagrangian submanifold along a direction 
which is transversal to $\mu^{-1}(0)$, e.g. along $i \eta$, where $\eta$ is the Hamiltonian vector field of 
the $S^1$ action. This submanifold is then the image of a section of the fibration in 
Example~\ref{ex amoebous fibr}.

In the case of $\Phi$ from Example~\ref{thin leg}, $\Phi(\Sigma)$ is a small perturbation of
$\Sigma'$ as above. One can see that the section in (\ref{log:section}) also avoids $\Phi(\Sigma)$.
Then the argument follows as before.
\end{proof}

We notice that ``smooth section'' in the above statement means a section whose image is a smooth, 
manifold. In fact there is no obvious notion of what a smooth map from the base is, since there 
is no notion of smooth coordinates. 

In view of Proposition \ref{prop. piecewise smoothness}, the fibrations of 
Examples~\ref{ex amoebous fibr} and \ref{thin leg} are all given by piecewise $C^\infty$ maps.
More precisely, away from $\Sigma$, they are the union of two honest 
$C^\infty$ fibrations meeting and coinciding along $\mu^{-1}(0)$. A similar 
phenomenon occurs in  special Lagrangian geometry \cite{Joyce-SYZ}. This kind of piecewise smoothness deserves careful attention and we study it in the next Section.

\section{Stitched fibrations}\label{stitched fibr}
In \cite{CB-M-torino} we proposed to extend the classical theory of action-angle coordinates to a particular type of piecewise smooth fibrations, which we called \emph{stitched fibrations}.  
Here we review how this theory was further developed in \cite{CB-M-stitched} and extend some of those techniques
to fibrations which are not proper.  For details and complete proofs we refer the reader to \cite{CB-M-stitched}. The material in this section is primarily technical but necessary to understand the lack of regularity of the fibrations in \S\ref{sec:pwfibr}. The techniques here are useful, in particular, for the construction of Lagrangian fibrations of negative type \S \ref{negative}.
\begin{defi}\label{defi stitched} 
Let $(X, \omega)$ be a smooth $2n$-dimensional symplectic manifold. Suppose there is a free Hamiltonian $S^1$ action on $X$ with moment map $\mu: X \rightarrow \R$. Let $X^+ = \{ \mu \geq 0 \}$ and $X^- = \{ \mu \leq 0 \}$. Given a smooth $(n-1)$-dimensional manifold $M$, a map $f: X \rightarrow \R \times M$ is said to be a \textit{stitched Lagrangian fibration} if there is a continuous $S^1$ invariant function $G: X \rightarrow M$, such that the following 
holds: 
\begin{itemize}
\item[\textit{(i)}] Let $G^{\pm} = G|_{X^{\pm}}$. Then $G^+$ and $G^-$ are restrictions of $C^\infty$ maps on $X$;
\item[\textit{(ii)}] $f$ can be written as $f =  (\mu, G)$ and $f$ restricted to $X^{\pm}$ is a proper submersion with connected Lagrangian fibres. 
\end{itemize}
We call $Z = \mu^{-1}(0)$ the \textit{seam} and $\Gamma=f(Z) \subseteq \{0\}\times M$ the \textit{wall}. We denote $f^\pm =f|_{X^\pm}$. 
\end{defi}


Notice that we do not require $f$ to be onto $\R \times M$, so we denote $B = f(X)$ and $B^{\pm} = f(X^{\pm})$.
In general, a stitched fibration will only be piecewise $C^\infty$, however all its fibres are smooth Lagrangian tori. Observe also that $f^\pm$ is the restriction of a $C^\infty$ map, it is not a priori required to extend to a smooth Lagrangian fibration beyond $X^\pm$.  Throughout this section we will always assume (unless otherwise 
stated) that the pair $(B, \Gamma)$ is diffeomorphic to the pair $(D^n, D^{n-1})$, where $D^k \subset \R^k$ is an 
open unit ball centered at the origin and $\R^{n-1}$ is embedded in $\R^n$. Later on we will consider more general bases --e.g. non-simply-connected-- 
when we speak about monodromy. 

\medskip
We now review some the examples given in \S\ref{sec:pwfibr}:

\begin{ex} [Stitched focus-focus, revisited] \label{broken focus focus rev}
Consider the piecewise smooth fibration in Example \ref{broken focus focus}. One can easily see that the  restriction of $f$ to $X - f^{-1}(0)$ is a stitched Lagrangian fibration. 
\end{ex}

Analogously, the piecewise smooth fibration in Example \ref{leg} gives rise to a stitched fibration when restricted to the complement of the union of the singular fibres. There is another important example in dimension three:

\begin{ex}[The amoeba, revisited] Consider the fibration in Example \ref{ex amoebous fibr}. When restricted to $X-f^{-1}(\Delta)$, $f$ defines a stitched Lagrangian fibration. The seam is 
$Z = \mu^{-1}(0) - f^{-1}(\Delta)$, notice that in this case $Z$ has three
connected components.
\end{ex}

\medskip
To understand the geometry of stitched fibrations in a neighborhood of a point 
on the wall, it is convenient to allow a more general set of coordinates than 
just the smooth ones.

\begin{defi}\label{defi:admissible}  
A set of coordinates on $B \subseteq \R \times M$, given by a map $\phi: B \rightarrow \R^n$, is said to be 
\textit{admissible} if the components of $\phi=(\phi_1,\ldots ,\phi_n)$ satisfy the following properties:
\begin{itemize}
\item[\textit{(i)}] $\phi_1$ is the restriction to $B$ of the projection map $\R \times M\rightarrow\R$;
\item[\textit{(ii)}]  for $j = 2, \ldots, n$ the restrictions of $\phi_j$ to $B^+$ and $B^-$ are locally 
restrictions of smooth functions on $B$.
\end{itemize}
\end{defi}

Essentially, admissible coordinates are those such that $\phi \circ f$ is 
again stitched. Let $f: X \rightarrow B$ be a stitched Lagrangian fibration and let $\phi$ be a set of admissible 
coordinates. For $j = 2, \ldots , n$,   $f_{j}^{\pm} = \phi_j  \circ  f|_{X^{\pm}}$ is the restriction of a $C^\infty$ function on $X$ to $X^\pm$ and
we can write $f=(\mu, f_2^\pm,\ldots ,f_n^\pm)$.
Let $\eta_1$ and $\eta_{j}^{\pm}$ be the Hamiltonian vector fields of $\mu$ and $f^{\pm}_{j}$ respectively. In order to measure how far $f$ is from being smooth, 
it makes sense to compare $\eta_j^+$ and $\eta_j^-$ in the only place where they 
exist simultaneously, i.e. along $Z$. In fact it is not difficult to show that 
there are $S^1$ invariant functions $a_j$ on $Z$ such that
\begin{equation}\label{discrep eq}
      (\eta^{+}_{j} - \eta^{-}_{j})|_{Z} = a_j \, 
                                         \eta_1|_{Z}. 
\end{equation}
Clearly, when $\phi\circ f$ is smooth $a_2=\cdots =a_n=0$. 

\medskip
It is convenient to interpret the $S^1$ invariant functions $(a_2,\ldots ,a_n)$ in (\ref{discrep eq}) as follows. First observe that the seam of a stitched fibration is an $S^1$-bundle $p:Z\rightarrow \bar Z:=Z\slash S^1$ such that: 
\[
\xymatrix{
Z \ar[dr]_{f|_Z} \ar[rr]^{p} &  & \bar{Z} \ar[dl]^{\bar{f}} \\
&  \Gamma
}
\]
where $\bar{Z}$ has the reduced symplectic form and $\bar f$ is the 
reduced Lagrangian fibration over the wall $\Gamma$.  We also have the vertical $(n-1)$-plane distribution: 
\[
\mathfrak{L}=\ker \bar f_\ast\subset T\bar Z   
\]
tangent to the fibres of $\bar f$. Clearly, a choice of coordinates around $b\in\Gamma$ induces a frame $\bar\eta =(\bar\eta_2,\ldots ,\bar\eta_n)$ of $\mathfrak L$, where $\bar\eta_j=p_\ast\eta_j^+=p_\ast\eta_j^-$.  Define $\ell_1$ to be the section of $\mathfrak L^\ast$ such that: 
\[
\ell_1 (\bar\eta_j)=a_j.
\]
It is not difficult to see (and we prove it in \cite{CB-M-stitched}) that  $\ell_1$ is \textit{fibrewise closed}, i.e. when restricted to the fibres of $\bar f$, $\ell_1$ is a closed 1-form. One can prove that a different choice of coordinates around $b\in\Gamma$ induces a frame $\bar\eta'$ and a section $\ell_1'$ such that $\ell_1'-\ell_1=\delta$, where $\delta$ is fibrewise constant, i.e. the Lie derivative $\mathcal L_{\bar\eta_j}\delta=0$ for all $j=2,\ldots ,n$ (cf. \cite{CB-M-stitched}{Proposition 4.2}). As a corollary, if there is a change of coordinates in the base which makes a stitched fibration smooth, then $\ell_1$ is fibrewise constant. The invariant $\ell_1$ is a first order 
measure of how much $f$ fails to be smooth along $Z$. Of course one also needs to consider ``higher order terms" 
to fully understand the behavior of a stitched fibration near the seam. 

\medskip
In the smooth case, action-angle coordinates defined over $B$ depend on a choice of a 
basis of $H_1(X,\Z)$. In the case of stitched fibrations it is convenient to generalize this idea as follows. 
We choose a pair of bases $\gamma^\pm=(\gamma_1,\gamma_2^\pm, \ldots ,\gamma_n^\pm)$ of $H_1(X,\Z)$ such that 

\begin{itemize}
\item[(a)] $\gamma_1$ is represented by an orbit of the $S^1$ action,
\item[(b)] $\gamma_j^{+} = \gamma_j^{-} + m_j \gamma_1$, for some $m_2, \ldots,m_n\in\Z$.
\end{itemize}

Condition (b) simply means that $p_\ast\gamma^+=p_\ast\gamma^-$ under the map $p_\ast:H_1(X,\Z)\rightarrow H_1(X/S^1,\Z)$. Such a choice of bases will be useful to understand fibrations over non simply connected bases where 
monodromy may occur. The following proposition generalizes the notion of action angle 
coordinates on the base.

\begin{prop}\label{prop:stitched_action} Let $f: X \rightarrow B$ be 
a stitched fibration and let $\gamma^\pm$ be bases of $H_1(X,\Z)$ satisfying the above conditions. 
Then the restrictions of $\gamma^\pm$ to $H_1(X^\pm, \Z)$ induce embeddings,
\[
\Lambda^\pm\hookrightarrow T^\ast_{B^\pm}.
\] 
Let $\alpha^\pm: B^\pm\rightarrow\R^n$ be the corresponding action coordinates 
satisfying $\alpha^\pm(b)=0$ for some $b \in \Gamma$. Then the map
\[
\alpha=
\begin{cases}
\alpha^+ &\textrm{on}\ B^+\\
\alpha^-&\textrm{on}\ B^-
\end{cases}
\] 
is an admissible change of coordinates. If $b_1,\ldots b_n$ denote the action coordinates on $B$ given 
by $\alpha$, then $\{db_1,\ldots db_n\}$ is a basis of $\Lambda^+$ and $\Lambda^-$. Furthermore, the reduced space $\bar Z$ can be identified with $T^\ast\Gamma\slash\langle db_2,\ldots ,db_n\rangle_\Z$ and the reduced fibration $\bar f$ can be identified with the standard projection $\bar\pi$. Moreover $\ell_1$ satisfies
\begin{equation}\label{eq:int_int}
\int_{[db_j]}\ell_1=m_j,\quad j=2,\ldots ,n
\end{equation}
where $[db_j]\in H_1(\bar Z,\Z)$ is the class represented by $db_j$.
\end{prop}
\begin{proof}
The first statements follow from the results in \S \ref{section: aff mfld & lag fib}. For the proof of the last statement we refer the reader to \cite{CB-M-stitched}{\S 4}.
\end{proof}

Recall that to establish the existence of action-angle coordinates, in the classical 
case, one chooses a smooth Lagrangian section. In the stitched case we choose 
a continuous section $\sigma: B \rightarrow X$ such that $\sigma|_{B^{\pm}}$ are 
the restrictions of smooth maps and $\sigma(B)$ is a smooth Lagrangian 
submanifold. Such sections always exist locally, for example the one constructed 
in Proposition~\ref{lag:section:exist} is a section of this type.
We denote a stitched fibration $f: X \rightarrow B$ together with a choice of basis $\gamma$ of $H_1(X, \Z)$ 
and a section $\sigma$ as above by $\mathcal F=(X, B, f, \gamma, \sigma)$.
\begin{defi}\label{def:symp-eq}
Two stitched fibrations $\mathcal F=(X, B, f, \gamma, \sigma)$ and $\mathcal F'=(X', B', f', \gamma', \sigma')$, 
with seams $Z$ and $Z'$ respectively are \textit{symplectically conjugate} if there are neighborhoods $W\subseteq B$ of $\Gamma :=f(Z)$ and $W'\subseteq B'$ of $\Gamma':= f'(Z')$ such that $\mathcal F|_{W}$ and $\mathcal F'|_{W'}$
are $(\psi, \phi)$-conjugate, where $\psi$ is an $S^1$ equivariant $C^\infty$ symplectomorphism sending $Z'$ to $Z$ and $\phi$ is a $C^\infty$ diffeomorphism such that $\psi\circ\sigma'=\sigma\circ\phi$ and $\psi_\ast \gamma' = \gamma$. The set of equivalence classes under this relation will be called \textit{germs of stitched fibrations}.
\end{defi}

Notice that in the above definition we are allowed to shrink to a smaller 
neighborhood of $\Gamma$ but not to a smaller $\Gamma$. So germs are meant 
to be defined around $\Gamma$ and not around a point. 
In \cite{CB-M-stitched} we classified stitched Lagrangian fibrations 
up to symplectic conjugation in terms of certain invariants. We review this classification 
here.

\medskip
First we illustrate a basic construction of stitched fibrations. 

\begin{ex}[Normal forms]
Let $(b_1,\ldots, b_n)$ be the standard coordinates on $\R^n$. Let $(U,\Gamma)$ 
be a pair of subsets of $\R^n$ diffeomorphic to $(D^n,D^{n-1})$ and $\Gamma = U\cap\{ b_1=0\}$. 
Define $U^+=U \cap\{b_1\geq 0\}$ and $U^-=U\cap\{ b_1\leq 0\}$. Consider the lattice $\Lambda= \spn\langle db_1,\ldots ,db_n\rangle_\Z$ and form the 
symplectic manifold $T^\ast U\slash\Lambda$. Denote by $\pi$ the standard 
projection onto $U$. Let $Z =\pi^{-1}(\Gamma)$ and $\bar Z = Z \slash S^1$, where the $S^1$ action is the one generated by $db_1$.
Suppose there is an open neighborhood $V\subseteq T^\ast U\slash\Lambda$ 
of $Z$ and a map $u:V\rightarrow\R^n$ which is a proper, smooth,  
$S^1$-invariant Lagrangian submersion with components $(u_1,\ldots ,u_n)$ such 
that $u|_{Z}=\pi$ and $u_1=b_1$. Now define the following subsets of 
$T^\ast U\slash\Lambda$,
\[
Y^+ := \pi^{-1}(U^+),\quad Y:= Y^+ \cup V,\quad Y^- := Y \cap \pi^{-1}(U^-)
\] 
and define the map $f_u: Y \rightarrow \R^n$ by
\begin{equation} \label{u:st}
       f_u = \begin{cases}
                u \quad\text{on} \  Y^-, \\
               \pi \quad\text{on} \ Y^+.
       \end{cases} 
\end{equation}
Clearly $f_u:Y\rightarrow \R^n$ is a stitched fibration. Denote  
$B_u:=f_u(Y)$. The zero section $\sigma_0$ of $\pi$ is, perhaps after a change of coordinates in the base, a section of $f_u$. Let $\gamma_0$ be the basis of $H_1(Y,\Z)$ induced by $\Lambda$. We call the stitched fibration $\mathcal F_u=(Y,B_u,f_u,\sigma_0,\gamma_0 )$ a \textit{normal form}.
\end{ex}
Now suppose $\mathcal F_u=(Y,B_u,f_u,\sigma_0,\gamma_0 )$ is as above and  let $(b,y) = (b_1, \ldots, b_n, y_1, \ldots, y_n)$ be canonical coordinates on $T^{\ast}B_u$
so that $y$ gives coordinates on the fibre $T_b^{\ast}B_u$. 
Let $W$ be a neighborhood of $\Gamma$ inside $u(V)$.
If $r \in \R$ is a parameter,  for any $b = (0, b_2, \ldots, b_n) \in \Gamma$, 
let $(r,b)$ denote the point $(r,b_2, \ldots, b_n) \in \R^n$. 
Given $(r,b) \in W$, denote by 
$L_{r,b}$ the fibre $u^{-1}((r,b))$. For every fibre $F_b \subset Z$ of 
$\pi$, consider the symplectomorphism
\begin{equation} \label{fb:symp}
 (y_1, \ldots, y_n, \sum_{k=1}^{n} x_k dy_k) 
            \mapsto (x_1, b_2 + x_2, \ldots, b_n + x_n, y_1, \ldots, y_n),
\end{equation}
between a neighborhood of the zero section of $T^{\ast}F_b$ and 
a neighborhood of $F_b$ in $V$. If $W$ is sufficiently small, for every 
$(r,b) \in W$, the Lagrangian 
submanifold $L_{r,b}$ will be the image of the graph of a closed $1$-form 
on $F_b$. Due to the $S^1$ invariance of $u$ and the fact that $u_1=b_1$, this 1-form has to be of the type
\[ r dy_1 + \ell(r,b), \]
where $\ell(r,b)$ is the pull back to $F_b$ of a closed one form on $\bar{F}_b$.  
Denote by $\ell(r)$ the smooth one parameter family of sections of 
$\mathfrak{L}^{\ast}$ such that $\ell(r)|_{\bar{F}_b} = \ell(r,b)$.
The condition $u|_Z=\pi$ implies that $\ell (0,b)=0$. Furthermore, the $N$-th order Taylor series expansion of 
$\ell(r)$ in the parameter $r$ can be written as
\begin{equation} \label{l:tay}
 \ell(r) =\sum_{k=1}^{N} \ell_k \, r^k + o(r^N), 
\end{equation}
where the $\ell_k$'s are fibrewise closed sections of $\mathfrak{L}^{\ast}$. 

\begin{defi}\label{def:seq} 
With the above notation, we define 
\begin{itemize}
\item[(i)] $\mathscr L_{Z}$ the set of sequences 
$\ell = \{\ell_k\}_{k\in\N}$ such that $\ell_k$ is a fibrewise closed section of 
$\mathfrak L^\ast$;
\item[ii)] $\mathscr U_{Z}$ the set of pairs $(V,u)$ where $V\subseteq T^\ast U\slash\Lambda$ is a neighborhood of $Z$ and $u:V\rightarrow\R^n$ is a proper, smooth,  $S^1$-invariant Lagrangian submersion with components $(u_1,\ldots ,u_n)$ such that $u|_{Z}=\pi$ and $u_1=b_1$.
\end{itemize}
\end{defi}

As above, to a given $(V,u) \in \mathscr U_{Z}$ we can associate a unique sequence 
$\ell \in \mathscr L_{Z}$. Conversely, in \cite{CB-M-stitched}{\S 5} we showed that for any given sequence $\ell\in\mathscr L_{Z}$ there is some $(V,u)\in\mathscr U_{Z}$, therefore a normal form, associated to it. Clearly, this $(V,u)$ is not unique.

\medskip
In \cite{CB-M-stitched} we proved that stitched fibrations are normalized 
according to the following:
\begin{prop}\label{prop:normalform} Every stitched fibration
$\mathcal F = (X, B, f, \sigma, \gamma)$ is symplectically conjugate to 
a normal form $\mathcal F_u=(Y, B_u, f_u, \sigma_0, \gamma_0 )$
\end{prop}
\begin{proof}
Let $Z$ be the seam of $\mathcal F$, $\omega_{red}$ the reduced 
symplectic form on $\bar Z$ and $\bar f: \bar Z \rightarrow \Gamma$ the 
reduced fibration.  
Using the coisotropic embedding theorem we can assume w.l.o.g.  that
 $X = \R\times S^1\times\bar Z$ with symplectic form
$\omega = \omega_{red}+ds\wedge dt$, where $(t,s)$ are coordinates on 
$\R\times S^1$ and the projection onto $\R$ is the moment map $\mu$. 
On $X$,  we can define an ``auxiliary" smooth Lagrangian fibration given by 
\[
\tilde\pi (t,s,p)=(t,\bar f(p)).
\]
Fix a basis $\gamma$ of $H_1(X,\Z)\cong H_1(S^1\times\bar Z,\Z)$ and a smooth Lagrangian section of $\tilde\pi$. The action-angle coordinates of $\tilde\pi$ 
with respect to $\gamma$ and $\sigma$ induce a $C^\infty$ symplectomorphism 
\begin{equation}\label{eq ident}
 T^\ast U\slash\Lambda \cong X
\end{equation}
for some open neighborhood $U$ of $0\in\R^n$ with action coordinates
$(b_1,\ldots ,b_n)$. The angle coordinates are $(y_1,\ldots ,y_n)$. 
In these coordinates $Z = \{ b_1 = 0 \}$ and $\Gamma = U \cap \{ b_1 = 0 \}$. 
While $f$ becomes:
\begin{equation}\label{eq:stitched-prime}
       f = \begin{cases}
       		u^+ \quad\text{on} \ X^+;\\
                u^- \quad\text{on} \  X^-,
       \end{cases} 
\end{equation}
where $u^\pm$ correspond to $f^\pm$. It follows that 
$u^+|_{Z}=u^-|_{Z}=\pi|_{Z}$.

One can show that $u^+$ can be extended as a smooth proper 
Lagrangian fibration a little bit beyond $X^+$, i.e. we can find a 
smooth proper Lagrangian fibration $\tilde u^+$ defined on a 
set $X^+ \cup V$, where $V$ is some open neighborhood of $Z$, such 
that $\tilde u^+|_{X^+} = u^+$. For the details 
of this extension see \cite{CB-M-stitched}, Proposition 6.3. 
To put $f$ in normal form, we consider the action-angle coordinates 
associated to $\tilde u^+$ with section $\sigma$ and basis $\gamma$ of $H_1(X,\Z)$ as above. In these coordinates, $X^+ \cup V$ becomes $T^{\ast}U / \Lambda$ 
and $\tilde u^+$ becomes the projection $\pi$. Again in action-angle coordinates of $\tilde u^+$, a Lagrangian extension $\tilde u^-$ of $u^-$, becomes $(W,u) \in \mathscr U_{Z}$ for some $W\subseteq T^\ast U\slash\Lambda$ and some Lagrangian fibration $u$. Then we simply define $Y^+= T^*U^+ / \Lambda$, $Y = Y^+ \cup W$,  $Y^- = Y \cap \pi^{-1} (U^-)$ and 
\begin{equation}
       f_u = \begin{cases}
                u \quad\text{on} \  Y^-, \\
               \pi \quad\text{on} \ Y^+.
       \end{cases} 
\end{equation}
\end{proof}
When $\mathcal F$ is smooth, its normal form is $\mathcal F_\pi$. This is Arnold-Liouville theorem (cf. Corollary \ref{smth:sg:equiv}). 
Given a stitched Lagrangian fibration $\mathcal F=(X,B,f, \sigma, \gamma)$
with normal form $\mathcal F_u=(Y, B_u, f_u, \sigma_0, \gamma_0 )$, we respectively denote by $\znor$ and $\gnor$ the seam and the wall 
of $\mathcal F_u$ and by $\zbnor$ the $S^1$ reduction of $\znor$.

\begin{defi}\label{def:seq-inv} 
Let $\mathcal F=(X,B,f, \sigma, \gamma)$ be a stitched fibration 
with normal form $\mathcal F_u=(Y, B_u, f_u, \sigma_0, \gamma_0 )$. 
Let $\ell\in\mathscr{L}_{\zbnor}$ be the unique sequence determined by $(V,u) \in \mathscr U_{\znor}$ defining $\mathcal F_u$.
We call $\inv (\mathcal F):=(\zbnor,\ell)$ the \textit{invariants} of 
$\mathcal F$. We say that the invariants of $\mathcal F$ vanish if for all $k\in\N$, $\ell_k\equiv 0$ when restricted to the reduced fibres of $\mathcal F_u$. We say that the invariants of $\mathcal F$ are fibrewise 
constant if all the $\ell_k$'s are fibrewise constant.
\end{defi}
We prove in \cite{CB-M-stitched}{Corollary 6.9} that $\inv (\mathcal F)$ is 
independent on the choice of normal form. 

\medskip
We will now see that every specified data $(\zbnor, \ell)$, with $\ell_1$ satisfying an integrality condition can be realized as the invariants of a stitched fibration. Notice that $\zbnor$ is uniquely 
determined by $\Gamma$ as $\zbnor = T^{\ast}\Gamma / \bar \Lambda$, where $\bar \Lambda = \spn \inner{db_2, \ldots}{db_n}_{\Z}$. We have

\begin{thm} \label{broken:constr2} 
Given any pair $(U,\Gamma_\mathrm{nor})$ of subsets of $\R^n$, diffeomorphic 
to $(D^n, D^{n-1})$ and with $\gnor = U \cap \{ b_1 = 0 \}$, 
a sequence $\ell = \{ \ell_k \}_{k \in \N} \in \mathscr L_{\zbnor}$ and integers $m_2, \ldots, m_n$ such that
\begin{equation} \label{int:cond2}
  \int_{[db_j]} \ell_1 = m_j, \quad \text{for all}\ j=2, \ldots, n,
\end{equation}
there exists a smooth symplectic manifold $(X, \omega)$ and a stitched
Lagrangian fibration $f: X \rightarrow U$ satisfying the following
properties:
\begin{itemize}
\item[(i)] the coordinates $(b_1, \ldots, b_n)$ on $U$ are action coordinates of $f$ with $\mu = f^{\ast}b_1$ the moment map of the $S^1$ action;
\item[(ii)] the periods $\{ db_1, \ldots, db_n \}$, restricted to $U^{\pm}$ correspond to bases $\gamma^{\pm} = \{ \gamma_1, \gamma_2^{\pm}, \ldots, \gamma_n^{\pm} \}$ of $H_1(X, \Z)$ satisfying conditions (a) and (b) prior to Proposition \ref{prop:stitched_action};
\item[(iii)]  there is a Lagrangian section $\sigma$ of $f$, such that $(\zbnor, \ell)$ are the invariants of $(X,f, U, \sigma, \gamma^+)$.
\end{itemize}
\end{thm} 
\begin{proof}
We refer the reader to \cite{CB-M-stitched}{Theorem 6.12} for the details. Roughly, one starts with $U^+$ and $U^-$ regarded as disjoint sets. These give two disjoint pieces $X^\pm=T^\ast U^\pm\slash\Lambda^\pm$, where $\Lambda^\pm=\langle db_1,\ldots ,db_n\rangle_\Z$. Let $Z^\pm=\partial X^\pm$. On $X^+$ we have Hamiltonian vector fields $\eta_1=\partial_{b_1}$ and $\eta^+_j=\partial_{b_j}$ for $j=2,\ldots ,n$. We can also define vector fields on $Z^+$:
\[
\eta^-_j=\eta^+_j-a_j\eta_1
\]
where $(a_2,\ldots , a_n)$ are the coefficients of $\ell_1$. One can (topologically) glue $X^+$ and $X^-$ using a map $Q:Z^-\rightarrow Z^+$ defined in terms of the $\R^{n}$ action induced by the flows of $\eta_j^-$. Intuitively, $Q$ identifies the fibres inside each of the two halves $Z^-$ and $Z^+$ after the fibres inside $Z^-$ have been twisted by iteratively flowing in the direction of $\eta_1, \eta^-_2,\ldots ,\eta^-_n$. The integrality condition (\ref{int:cond2}) guarantees that (ii) is satisfied. One can extend $Q$ to give a \emph{smooth} symplectomorphism $\tilde{ Q}$ between open neighborhoods of $Z^\pm$. For this one needs to consider invariants $\ell_k$, for $k>1$. The choice of $\tilde Q$ is determined by $\{\ell_k\}$. This gluing gives a smooth symplectic manifold $(X,\omega)$ and a stitched fibration $f:X\rightarrow U$, which by construction is such that $\inv (\mathcal F)=(\zbnor, \ell)$.
\end{proof}

We also have the following (cf. \cite{CB-M-stitched}{Theorem 6.11}):

\begin{thm}\label{thm: grosso} Let $\mathcal F$ and $\mathcal F'$ be stitched fibrations. Then,
\begin{itemize}
\item[(i)] two stitched fibrations $\mathcal F$ and $\mathcal F'$ are conjugate if and only if $\inv (\mathcal F)=\inv (\mathcal F')$;
\item[(ii)] $\mathcal F$ is smooth if and only if $\inv (\mathcal F)$ vanish;
\item[(iii)] $\mathcal F$ becomes smooth after an admissible change of coordinates on the base if and only if $\inv (\mathcal F)$ are fibrewise constant.

\end{itemize}
\end{thm}
In other words, the set of germs of stitched fibrations is classified by the pairs 
$(\zbnor, \ell)$. We say that a fibration is \textit{fake stitched} if it becomes smooth 
after an admissible change of coordinates on the base. 
One interesting consequence of Theorem~\ref{broken:constr2}, which we will 
exploit later on, is that from a given set of invariants we can form 
another one for example by summing to the sequence $\ell$ another sequence
or by multiplying elements $\ell_k$ by pull backs of smooth functions on 
the base. The new invariants give rise to new stitched fibrations. 

\begin{ex} \label{fake:inv} 
Consider a smooth proper Lagrangian fibration $f: X \rightarrow B$, with $B = \R \times M$ and $f= (\mu, G)$, where $\mu$ is the moment map of a free $S^1$ 
action and $G$ is $S^1$ invariant. Assuming $B$ is contractible and having 
chosen bases 
$\gamma^{\pm}$ of $H_1(X, \Z)$ as in $(a)$ and $(b)$ above, on $B$ 
we can apply the admissible change of coordinates $\alpha$ as in 
Proposition~\ref{prop:stitched_action}. 
Clearly $f' = \alpha \circ f$ is (tautologically) a fake stitched fibration. Given a 
Lagrangian section $\sigma$ of $f'$, it easy to see that the normal form for 
$(X, f'(B), f', \gamma^+, \sigma)$ is of the type $(Y, U, f_u, \gamma_0, \sigma_0)$ where  
$Y = T^*U / \Lambda$ and 
\[ u(y_1, \ldots, y_n, b_1, \ldots, b_n) =(b_1, b_2 - m_2 b_1, \ldots, b_n - m_n b_1), \]
i.e. the projection composed with a linear change of coordinates.
In this case the only non-zero invariant is $\ell_1$ which is given by
\[ \ell_1 = \sum_j \, m_j dy_j. \]
Clearly $\ell_1$ is fibrewise constant.
\end{ex}

\subsection{Monodromy}
We now study stitched fibrations defined over non simply connected bases. In this case, the underlying topological $T^n$ bundle may have monodromy. When $\mathcal F$ is smooth, monodromy can be read from the holonomy of the affine structure on the base. This is no longer true for stitched fibrations in general. This is the case, for instance, of Example \ref{ex amoebous fibr}; in fact, in \cite{CB-M-torino}{Proposition 7} (cf. also Remark 5) we gave explicit evidence of this. We show now that monodromy can alternatively be detected from the behavior of the first order invariant $\ell_1$. We restrict to some specific examples with unipotent monodromy.

\begin{ex}  \label{two:mon}
Let $U \subset \R^2$ be an open annulus in $\R^2$ centered at the origin.
As usual denote $U^+ = U \cap \{ b_1 \geq 0 \}$, $U^- = U \cap \{ b_1 \leq 0 \}$ 
and $\Gamma = U^+ \cap U^-$. This time $\Gamma$ is disconnected.
We let $\Gamma_{u} = \Gamma \cap \{ b_2 \geq 0 \}$ and 
$\Gamma_{d} = \Gamma \cap \{ b_2 \leq 0 \}$ be the upper and lower parts of 
$\Gamma$ respectively. 
Now let $f: X \rightarrow \R^2$ be a stitched Lagrangian fibration such that
$f(X) = U$. Observe that the seam $Z$ has two connected 
components: $Z_u = f^{-1}(\Gamma_u)$ and $Z_d = f^{-1}(\Gamma_d)$. 
Denote by $\bar{Z}_u$ and $\bar{Z}_d$ the respective $S^1$ quotients, i.e.
the connected components of $\bar{Z}$.
Let $b \in \Gamma_u$ and choose as generator of $\pi_1(U, b)$ an anti-clock-wise oriented curve starting at $b$ and going once around $0$. Suppose that with respect to a basis
 $\{ \gamma_1, \gamma_2 \}$ of $H_1(F_b, \Z)$ the monodromy is
\begin{equation} \label{monodr}
  \left( \begin{array}{cc} 1 & -m \\
                     0 & 1 
\end{array} \right),
\end{equation}
for some integer $m\neq 0$. In this case we must have that
$\gamma_1$ is represented by the orbits of the $S^1$ action. 
As usual let $X^{\pm} = f^{-1}(U^{\pm})$. 
Since $U - \Gamma_{d}$ is contractible we can
think of $\{ \gamma_1, \gamma_2 \}$ as a basis of $H_1(f^{-1}(U-\Gamma_{d}), \Z)$.
Consider the diagrams:

\[
\xymatrix{
 &  H_1(X^{+},\Z) \ar[dr]  \\
 H_1(f^{-1}(U-\Gamma_d),\Z) \ar[ur] \ar[rr]^{j_+} & & H_1(f^{-1}(U-\Gamma_u),\Z)
}
\]
or
\[
\xymatrix{
H_1(f^{-1}(U-\Gamma_d),\Z) \ar[dr] \ar[rr]^{j_-} & & H_1(f^{-1}(U-\Gamma_u),\Z) \\
&  H_1(X^{-},\Z) \ar[ur]
}
\]
induced by inclusions and restrictions.  The map $j_+$ identifies 
$\{ \gamma_1, \gamma_2 \}$ with a basis $\{ \gamma_1, \gamma_2^+ \}$ 
of $H_1(f^{-1}(U-\Gamma_{u}), \Z)$, whereas $j_-$ with a basis
$\{ \gamma_1, \gamma_2^- \}$. 
Notice that monodromy is given by $j_+^{-1}\circ j_-$.  
Therefore we must have $\gamma_2^+ = m \gamma_1 + \gamma_2^-$.
Hence $\{ \gamma_1, \gamma_2^+ \}$ and $\{ \gamma_1, \gamma_2^- \}$ satisfy conditions 
(a) and (b) in the previous section.  
Applying Proposition~\ref{prop:stitched_action} to $f$ restricted to  $f^{-1}(U-\Gamma_u)$ 
we can consider the action coordinates map $\alpha$ constructed
by taking action coordinates with respect
to $\{ \gamma_1, \gamma_2^+ \}$ on $U^+$ and with respect to 
$\{ \gamma_1, \gamma_2^- \}$ on $U^-$.
Denote by $(b_1^d, b_2^d)$ such coordinates. Similarly on $U-\Gamma_{d}$ we 
can consider action angle coordinates with respect to the basis
$\{ \gamma_1, \gamma_2 \}$. Denote by $(b_1^u, b_2^u)$ these coordinates.
In particular we have the identifications
\[ \bar{Z}_d = T^{\ast} \Gamma_d \, /  \, \langle db_2^d \rangle_{\Z} \]
and
\[ \bar{Z}_u= T^{\ast} \Gamma_u \, / \, \langle db_2^u \rangle_{\Z}. \]
With respect to this choice of coordinates we can compute
the first order invariants of $f$, $\ell_1^u$ and $\ell_1^d$ on $\bar{Z}_u$ and $\bar{Z}_d$, respectively. Then
(\ref{eq:int_int}) should hold, therefore we obtain
\[ \int_{[db_2^u]} \ell_1^u = 0 \ \ \text{and} \ \ \int_{[db_2^d]} \ell_1^d = m. \]
This tells us that monodromy can be read from a jump in cohomology class
of the first order invariant associated to action coordinates.
\end{ex}

Using the methods of Theorem~\ref*{broken:constr2}
we can also construct stitched Lagrangian fibrations with prescribed monodromy and invariants. 
In fact we have

\begin{thm} \label{stitch:mon_ff} Let $U \subset \R^2$ be an annulus as above with 
coordinates $(b_1, b_2)$. Let
$\bar{Z}_d = T^{\ast} \Gamma_d \, /  \, \langle db_2 \rangle_{\Z}$ 
and
$\bar{Z}_u= T^{\ast} \Gamma_u \, / \, \langle db_2 \rangle_{\Z}$
with projections $\bar{\pi}^d$ and $\bar{\pi}^u$ and bundles
 $\mathfrak{L}_d = \ker\bar{\pi}^d_{\ast}$ and $\mathfrak{L}_u = \ker\bar{\pi}^u_{\ast}$
 respectively. Given an integer $m$ and sequences 
$\ell^d = \{ \ell_k^d \}_{k \in \N} \in \mathscr L_{\bar{Z}_d}$ and 
$\ell^u = \{ \ell_k^u \}_{k \in \N} \in \mathscr L_{\bar{Z}_u}$ such that
\[ \int_{[db_2]} \ell_1^u = 0 \ \ \text{and} \ \ \int_{[db_2]} \ell_1^d = m, \]
there exists a smooth symplectic manifold $(X, \omega)$ and a stitched 
Lagrangian fibration $f: X \rightarrow U$ having monodromy (\ref*{monodr}) with
respect to some basis $\gamma = \{ \gamma_1, \gamma_2 \}$ of $H_{1}(f^{-1}(U- \Gamma_{d}), \Z)$
and satisfying the following properties:
  \begin{itemize}
    \item[(i)] the coordinates $(b_1, b_2)$
          are action coordinates of $f$ with moment map $f^{\ast}b_1$; 
    \item[(ii)] the periods $\{ db_1, db_2 \}$, restricted to $U^{\pm}$ 
          correspond to the basis $\{ \gamma_1, \gamma_2  \}$; 
    \item[(iii)]  there is a Lagrangian section 
          $\sigma$ of $f$, such that $(\bar{Z}_u, \ell^u)$ and $(\bar{Z}_d, \, \ell^d)$ are the invariants of $(f^{-1}(U-\Gamma_d),\, f, \, U-\Gamma_d, \,  \sigma, \, \gamma)$ and
          $(f^{-1}(U- \Gamma_u),\, f, \, U-\Gamma_u, \, \sigma, \, j_+(\gamma))$ respectively.
\end{itemize} 
The fibration $(X,f, U)$ satisfying the above properties is unique up to 
fibre preserving symplectomorphism. 
\end{thm}
\begin{proof} This is just a repetition of the arguments in Theorem~\ref*{broken:constr2} for each component of $\Gamma=\Gamma_d\cup\Gamma_u$. We leave the details as an exercise.
\end{proof}

\begin{rem} \label{skew_ff:ex} 
Notice that the stitched fibrations discussed in Example~\ref{two:mon} are 
more general than the ones constructed in Theorem~\ref{stitch:mon_ff}. We  illustrate this with an example. Let $U^-$ and $U^+$ be two 
``half annuli'' of the same width but of different radii (as depicted in 
Figure~\ref{half anuli}). If $b^{\pm} = (b_1^{\pm}, b_2^{\pm})$ denote coordinates 
on $U^{\pm}$ and we let $\Lambda^{\pm} = \inner{db_1^{\pm}}{db_2^{\pm}}_{\Z}$, then 
we can glue together $X^{+} = T^*U^+ / \Lambda^+$ and $X^{-} = T^*U^- / \Lambda^-$ after choosing suitable invariants and applying the usual method of 
Theorem~\ref{broken:constr2}. We first glue the lower boundaries of $X^+$ and 
$X^-$ and then the upper boundaries, (as indicated by the arrows in Figure~\ref{half anuli}). 
This produces a stitched fibration of the type discussed 
in Example~\ref{two:mon}, in fact we would obtain a total space $X$ which fibres over 
a base obtained as the result of the gluing of the two half annuli, which is clearly 
diffeomorphic to an annulus.  The fibration is not of the type constructed in 
Theorem~\ref{stitch:mon_ff}. There are two main differences between the two constructions.  In the examples from Theorem~\ref{stitch:mon_ff} action coordinates extend continuously to the whole annulus and the symplectic form on the total
space is exact. These two facts do not hold in the example just described, in fact if the 
symplectic form were exact then the action coordinates would extend continuously to the 
whole annulus (to show this one can use an argument similar to the one used in 
Proposition~\ref{prop:pos_simple}). 
\begin{figure}[ht]
\psfrag{U+}{$U^+$}
\psfrag{U-}{$U^-$}
\begin{center}
\epsfig{file=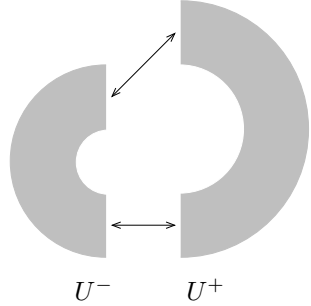, width=4cm,height=4cm}
\caption{Gluing half annuli with different radii.}\label{half anuli}
\end{center}
\end{figure}
\end{rem}

\begin{ex}
An example of a stitched Lagrangian fibration constructed using 
Theorem~\ref{stitch:mon_ff} is the following. We can choose
the elements of the sequence $\ell^u$ to be all zero, while the elements
of the sequence $\ell^d$ to be all zero except $\ell_1^d$ which we 
define to be
\[ \ell_1^d = m \, dy_2. \]
It is clear that the resulting fibration is only fake stitched, in fact
the invariants are fibrewise constant. One can also see 
that, in the case $m=1$ and $U = \R^2 - \{0 \}$,  the
fibration is symplectically conjugate to $(X, \alpha \circ f)$, where
$(X,f)$ is a smooth focus-focus fibration (where the singular 
fibre has been removed) and $\alpha$ is the action coordinates map 
(see the discussion after Example~\ref{smooth:ff} and Example~\ref{fake:inv}). 
In particular 
this fibration induces an affine structure on the base which is simple.
\end{ex}

We now discuss a three dimensional example.

\begin{ex} \label{amoeb:mon} \todo{Check the $m$'s and $g$'s subindices are correct}
In $\R^3$ consider the 3-valent graph 
\[ \Delta = \{(0,0,- t), \ t \geq 0 \} \cup \{ (0,- t,0), \  t \geq 0 \} 
             \cup \{ (0,t,t), \ t\geq 0 \} \]
and let $D$ be a tubular neighborhood of $\Delta$. Take $U = \R^3 - D$
and assume we have a stitched Lagrangian fibration $f: X \rightarrow \R^3$
such that $U=f(X)$ and the seam is $Z= f^{-1}( \{ b_1 = 0 \} \cap U)$. 
Again we let 
$U^+ = U \cap \{ b_1 \geq 0 \}$, $U^- = U \cap \{ b_1 \leq 0 \}$ 
and $\Gamma = U^+ \cap U^-$. Also let $X^{\pm} = f^{-1}(U^{\pm})$. 
This time $\Gamma$ (hence $Z$) has three connected components
\begin{eqnarray*}
  \Gamma_c & = & \{ ( 0, t,s), \ t,s < 0 \} \cap U, \\
  \Gamma_d & = & \{ (0, t,s), \ t> 0, s < t \} \cap U, \\
  \Gamma_e & = & \{ (0, t,s), \ s> 0, t < s \} \cap U.
\end{eqnarray*}
Also denote by $Z_c$, $Z_d$ and $Z_e$ the corresponding connected
components of $Z$ and by $\bar{Z}_c$, $\bar{Z}_d$ and $\bar{Z}_e$
their $S^1$ quotients.

\medskip

Fix $b \in \Gamma_c$ and suppose that there is a basis 
$\{ \gamma_1, \gamma_2, \gamma_3 \}$ of $H_1( F_b, \Z)$
and generators $g_1, g_2, g_3$ of $\pi_1(U, b)$, satisfying $g_1 g_2 g_3 = 1$, with respect to which the monodromy transformations are
\begin{equation} \label{mon:m}
\mathcal{M}_b(g_1) =  T_1  =  \left( \begin{array}{ccc}
                 1 & -m_1 & 0 \\
                 0 & 1  & 0 \\
                 0 & 0  & 1 
              \end{array} \right), \ \ \ 
\mathcal{M}_b(g_2) =  T_2  = \left( \begin{array}{ccc}
                 1 & 0 &  -m_2 \\
                 0 & 1  & 0 \\
                 0 & 0  & 1 
              \end{array} \right)
\end{equation}
and $\mathcal M_b(g_3) = T_3 = T_2^{-1} T_1^{-1}$, for non zero integers $m_1$ and $m_2$.  We have that $\gamma_1$ is represented by the orbits of the
$S^1$ action, since it is the only monodromy invariant cycle. 
Now, since $U - ( \Gamma_d \cup \Gamma_e)$ is
contractible, $\{ \gamma_1, \gamma_2, \gamma_3 \}$ is
a basis of $H_1( f^{-1}(U - ( \Gamma_d \cup \Gamma_e)), \Z)$.
Consider the diagrams:
\[
\xymatrix{
 &  H_1(X^{+},\Z) \ar[dr]  \\
 H_1(f^{-1}(U-(\Gamma_d\cup\Gamma_e)),\Z) \ar[ur] \ar[rr]^{j_+} & & H_1(f^{-1}(U-(\Gamma_c\cup\Gamma_d)),\Z)
}
\]
or
\[
\xymatrix{
H_1(f^{-1}(U-(\Gamma_d\cup\Gamma_e)),\Z) \ar[dr] \ar[rr]^{j_-} & & H_1(f^{-1}(U-(\Gamma_c\cup\Gamma_d)),\Z) \\
&  H_1(X^{-},\Z) \ar[ur]
}
\]
induced by inclusions and restrictions. The map $j_+$ identifies
$\{ \gamma_1, \gamma_2, \gamma_3 \}$ with a basis  of 
$H_1(f^{-1}(U- (\Gamma_c \cup \Gamma_d)), \Z) $, which 
we call $\{ \gamma_1, \gamma_2^+, \gamma_3^+ \}$, while $j_-$ 
identifies it with another basis, which we call 
$\{ \gamma_1, \gamma_2^-, \gamma_3^- \}$. 
Notice that the monodromy map $\mathcal{M}_b(g_2) = j_+^{-1}\circ j_-$. We must have 
\begin{equation} \label{tre:mon} 
 \begin{cases}
  \gamma_2^+ = \gamma_2^-, \\
  \gamma_3^+ = m_2 \gamma_1 + \gamma_3^- .
 \end{cases}
\end{equation}
Applying Proposition~\ref{prop:stitched_action} to $f$ restricted to 
$f^{-1}(U- (\Gamma_c \cup \Gamma_d))$, we can consider the action coordinates 
map $\alpha$ on $U- (\Gamma_c \cup \Gamma_d) $ computed with respect to $\{ \gamma_1, \gamma_2^+ ,\gamma_3^+ \}$ on $U^+$ 
and with respect to $\{ \gamma_1, \gamma_2^-, \gamma_3^- \}$ on $U^-$.
Let us denote these coordinates by $(b_1^e, b_2^e, b_3^e)$. Similarly we can 
consider action coordinates on $U - (\Gamma_d \cup \Gamma_e)$
with respect to the basis $\{ \gamma_1, \gamma_2, \gamma_3 \}$  of 
$H_1(f^{-1}(U- (\Gamma_d \cup \Gamma_e)), \Z) $. We denote them
by $(b_1^c, b_2^c, b_3^c)$.
We have the identifications
\[ \bar{Z}_e = T^{\ast} \Gamma_e \, /  \, \langle db_2^e, db_3^e  \rangle_{\Z} \]
and
\[ \bar{Z}_c= T^{\ast} \Gamma_c \, / \, \langle db_2^c, db_3^c \rangle_{\Z}. \]

With respect to these coordinates we can compute the first order invariants 
$\ell_1^e$ and $\ell_1^c$ on $\bar{Z}_e$ and $\bar{Z}_c$ respectively.
From Proposition~\ref{prop:stitched_action} and identities (\ref*{tre:mon}) 
applied to $\ell_1^c$ and 
$\ell_1^e$ we obtain
\begin{equation*}
  \int_{[db_2^c]} \ell_1^c  
=  \int_{[db_3^c]} \ell_1^c  = 0
\end{equation*}
and
\begin{equation*}
  \int_{[db_2^e]} \ell_1^e =  0 \ \ \text{and} \ \ 
  \int_{[db_3^e]} \ell_1^e  =  m_2.
\end{equation*}
Similarly we construct the first order invariant $\ell_1^d$ on $\bar{Z}_d$. 
It will satisfy
\begin{equation*}
  \int_{[db_2^d]} \ell_1^d  =  m_1 \ \ \text{and} \ \ 
  \int_{[db_3^d]} \ell_1^d =  0.    
\end{equation*}
Again, monodromy is understood in terms of the difference in the
cohomology class of the first order invariant. Example \ref*{ex amoebous fibr} 
is a special case of this situation, where $m_1 = m_2 = 1$. 
\end{ex}

Conversely, we can construct stitched fibrations like the one in previous example by specifying gluing data and applying Theorem~\ref*{broken:constr2}. In fact
we can prove
\begin{thm} \label{stitch:monodr_neg} 
Let $U \subset \R^3$, $\Gamma_c$, $\Gamma_d$ and $\Gamma_e$
be as in Example~\ref*{amoeb:mon} and let $(b_1, b_2, b_3)$ be coordinates
on $U$. Define
$\bar{Z}_c = T^{\ast} \Gamma_c \, /  \, \langle db_2, db_3 \rangle_{\Z}$,
$\bar{Z}_d= T^{\ast} \Gamma_d \, / \, \langle db_2, db_3 \rangle_{\Z}$
and $\bar{Z}_e= T^{\ast} \Gamma_e \, / \, \langle db_2, db_3 \rangle_{\Z}$
with projections $\bar{\pi}^c$, $\bar{\pi}^d$, $\bar{\pi}^e$ and bundles
$\mathfrak{L}_c = \ker\bar{\pi}^c_{\ast}$, $\mathfrak{L}_d = \ker\bar{\pi}^d_{\ast}$,
$\mathfrak{L}_e = \ker\bar{\pi}^e_{\ast}$. 
Suppose we are given integers $m_1$, $m_2$ and sequences
$\ell^c = \{ \ell_k^c \}_{k \in \N} \in \mathscr L_{\bar{Z}_c}$, 
$\ell^d = \{ \ell_k^d \}_{k \in \N} \in \mathscr L_{\bar{Z}_d}$ and 
$\ell^e = \{ \ell_k^e \}_{k \in \N} \in \mathscr L_{\bar{Z}_e}$ satisfying 
\begin{eqnarray} 
  \int_{[db_2]} \ell_1^c  & = & \int_{[db_3]} \ell_1^c = 0, \nonumber\\
  \int_{[db_2]} \ell_1^e & =  &0\quad\ \ \text{and} \quad
  \int_{[db_3]} \ell_1^e  =  m_2, \label{st:neg_int} \\
  \int_{[db_2]} \ell_1^d  & = & m_1 \quad \text{and} \quad
  \int_{[db_3]} \ell_1^d  =  0.    \nonumber 
\end{eqnarray}
Then there exists a smooth symplectic manifold $(X, \omega)$ and a stitched 
Lagrangian fibration $f: X \rightarrow U$ having the same monodromy
of Example~\ref*{amoeb:mon} with respect to some basis 
$\gamma = \{ \gamma_1, \gamma_2, \gamma_3 \}$ of 
$H_1( f^{-1}(U - ( \Gamma_d \cup \Gamma_e)), \Z)$   and satisfying the following 
properties:

\begin{itemize}
\item[(i)] the coordinates $(b_1, b_2, b_3)$ are action coordinates of $f$ with moment map $f^{\ast}b_1$; 
\item[(ii)] the periods $\{ db_1, db_2, db_3 \}$, restricted to $U^{\pm}$ 
          correspond to the basis $\gamma$; 
\item[(iii)]  there is a Lagrangian section $\sigma$ of $f$, such that $(\bar{Z}_c, \ell^c)$, $(\bar{Z}_d, \, \ell^d)$ and $(\bar{Z}_e, \, \ell^e)$ are respectively the invariants of:
\[ (f|_{U-( \Gamma_d \cup \Gamma_e)},  \sigma, \, \gamma),\ 
(f|_{U-( \Gamma_c \cup \Gamma_e)},  \sigma, \, j_+(\gamma))\ \textrm{and}\  
(f|_{U-( \Gamma_c \cup \Gamma_d)},  \sigma, \, j_+(\gamma)).\]
\end{itemize} 
The fibration $(X,f, U)$ satisfying the above properties is unique up to 
fibre preserving symplectomorphism. 
\end{thm}

\begin{rem} \label{skew_neg:re}
Also in this case (cf. Remark~\ref{skew_ff:ex}) we notice that 
fibrations of the type discussed in Example~\ref{amoeb:mon} are more general 
than the ones constructed using Theorem~\ref{stitch:monodr_neg}. 
To show this one can use higher dimensional versions of the fibration 
in Remark~\ref{skew_ff:ex}, with discontinuous action coordinates. We leave
the details to the reader.
\end{rem}

\begin{ex} \label{neg:as_stitched}
A simple example of stitched Lagrangian fibration which can be constructed
using Theorem~\ref{stitch:monodr_neg} is as follows. Define the sequence $\ell^c$ to be identically zero and choose the terms of $\ell^d$ and 
$\ell^e$ to be zero except the first order ones, which we define to be
\[ \ell_1^d = m_1 \, dy_2 \ \ \text{and} \ \ \ell_1^e = m_2 \, dy_3. \]
Clearly $\ell_1^c$, $\ell_1^d$ and $\ell_1^e$ satisfy the integral conditions
of Theorem~\ref{stitch:monodr_neg}, moreover they are fibrewise constant, 
therefore they define fake stitched fibrations. Since the fibration is 
smooth after a change of coordinates on the base, it induces an affine structure
on the base. One can easily see that in the case $m_1=-1$ and $m_2=1$ and \todo{changed the sign of $m_1$, to adjust to new convention} 
$U = \R^3 - \Delta$, this affine structure is simple and affine isomorphic 
to a negative vertex of Example~\ref{neg:aff}. Notice that we could also 
replace $\Delta$ with $\Delta_{\tau}$ and obtain an affine structure which 
is isomorphic to the one in Example~\ref{neg:affvar}.
\end{ex}

\subsection{Non-proper stitched fibrations}
This section is rather technical and the methods introduced will only be used in the proof of 
Lemma~\ref{lem:smth_crit2}, therefore the reader may skip it on first reading. 
Here we study some special cases of piecewise smooth fibrations with non compact fibres. 
The results extend the ones concerning proper maps. For this reason and for 
sake of brevity we shall only give full 
proofs when the arguments do not follow directly from the previous case.

\medskip
Let $X$ be a smooth symplectic $6$-manifold together with a smooth Hamiltonian $S^1$ action with moment map 
$\mu:X\rightarrow\R$. Assume $\mu$ has exactly one critical value $0\in\R$ and a codimension four submanifold 
$\Sigma=\Crit \mu$. Let $M$ be a smooth $2$-dimensional manifold and let $B \subseteq \R \times M$ be a 
contractible open neighborhood of a point $(0,m) \in \R \times M$. Let $\Gamma =B \cap (\{0\}\times M)$. As usual we define $Z = \mu^{-1}(0)$ and $\bar Z$ the $S^1$ quotient of $Z$ and $X^+=\{ \mu\geq 0\}$, $X^-=\{ \mu\leq 0\}$.

We consider fibrations satisfying the following:

\begin{ass}\label{ass: semi-stitched} The map $f:X\rightarrow B$ is a topological $T^3$ fibration with discriminant locus $\Delta \subset \Gamma$ such that $f(\Sigma )=\Delta$ satisfying  
\begin{itemize}
\item[(a)] $(X,\omega,f,B)$ is topologically conjugate to a generic singular fibration.
\item[(b)] There is a continuous $S^1$ invariant map $G:X\rightarrow M$ such that 
\begin{itemize}
\item[(i)] if $G^{\pm} = G|_{X^{\pm}}$ then $G^+$ and $G^-$ are restrictions of $C^\infty$ maps on $X$;
\item[(ii)] $f$ can be written as $f =  (\mu, G)$ and $f$ restricted to $X^{\pm}$ is a proper map with connected Lagrangian fibres. 
\end{itemize}
\item[(c)] There is a connected, $S^1$ invariant, open neighborhood $\mathfrak U\subseteq X$ of $\Sigma$ such that $f(\mathfrak U)=B$ and such that $f_\mathfrak{U}=f|_\mathfrak{U}$ is a $C^\infty$ map with non degenerate singular 
points.

\end{itemize}
 
\end{ass}
We can think of $B$ as $D^2 \times I$ with $\Delta = \{ 0 \} \times I$. Clearly, the restriction of $f$ to $X-f^{-1}(\Delta)$ is a 
stitched fibration in the sense of the previous sections. Example \ref{leg}, as well as the legs of Example \ref{thin leg} satisfy conditions (a) and (b). Furthermore, one can deform such examples near $\Sigma$ to produce fibrations which, in addition, satisfy condition (c) (cf. Lemma \ref{lem:smth_crit}).

\medskip
Let $\mathfrak U'\subset\mathfrak U$ be a smaller open set satisfying condition (c) 
(maybe after shrinking $B$). 
If we remove $\mathfrak U'$ we obtain a topologically trivial compact cylinder fibration
\begin{equation}\label{f_circ}
f|_{X - \mathfrak U'}:X -\mathfrak U'\rightarrow B
\end{equation}
which fails to be smooth along a subset of $Z - (\mathfrak U' \cap Z)$. 
Notice though that the fibration is actually smooth toward the ends of each cylindrical fibre. 
 
\medskip
Let $X^\circ= X -\overline{\mathfrak U'}$ with symplectic structure $\omega^\circ=\omega|_{X^\circ}$. The restriction $f^\circ =f|_{X^\circ}$ defines a piecewise smooth open cylinder fibration
\begin{equation}\label{eq:f^c}
f^\circ:X^\circ\rightarrow B.
\end{equation}
 We denote $F^\circ(b)$ the cylindrical fibre of $f^\circ$ over $b \in B$. 
On the other hand, the smooth part $f_\mathfrak{U}$ of $f$ defines an integrable Hamiltonian system with 
non-degenerate singularities which can be normalized as in Theorem~\ref{thm. normal form}. This normalization defines smooth coordinates $(b_1,b_2,b_3)$ on the base. 

\medskip
Denote by $X^\#=X - \Sigma$ and by $f^\#:X^\#\rightarrow B$ the restriction of $f$ to $X^\#$. Let $(f^\#)^\pm$ be the restriction of $f^\pm$ to $(X^\#)^\pm =X^{\#}\cap X^\pm$ and let $Z^\#=Z - \Sigma$ and $\bar Z^\#$ the corresponding reduced space with reduced symplectic structure $\omega_{red}$ on $\bar Z^\#$. 

\begin{prop}\label{prop:periods_sh} \todo{I fixed the statement. Check everything is ok.}
Let $f: X \rightarrow B$ be a fibration satisfying 
Assumption~\ref{ass: semi-stitched} and let $F_{\bar b} = f^{-1}(\bar b)$ be a smooth fibre. 
There is a  basis $\gamma=(\gamma_1,\gamma_2,\gamma_3)$ of $H_1(F_{\bar b},\Z)$ and coordinates
$(b_1,b_2,b_3)$ on $B$ with respect to which the periods of $f^{\pm}: X^{\pm} \rightarrow B^{\pm}$
can be written
\[
\begin{array}{l}
\lambda_{1}^{\pm} = 2\pi db_1, \\
\lambda_2^\pm=dH^\pm+\lambda_0, \\
\lambda_3^{\pm} = db_3,
\end{array}
\]
where $\lambda_0=\arg (b_1+ib_2)db_1+\log |b_1+ib_2|db_2$ and $H^\pm\in C^\infty (B^\pm)$.
Moreover, there is a fibre preserving symplectomorphism  
\begin{equation}
\Theta^\pm:T^\ast B^\pm\slash\Lambda_{H^\pm}\rightarrow (X^\#)^\pm 
\end{equation}
where $\Lambda_{H^\pm}$ is the integral lattice generated by 
$\lambda_1^{\pm}, \lambda_2^{\pm}, \lambda_3^{\pm}$.

\end{prop}
\begin{proof}
We take as coordinates $(b_1, b_2, b_3)$ on $B$ the ones given by the normalization 
of the singularity in Theorem~\ref{thm. normal form}.
Then the proof goes essentially as in Proposition \ref{prop. generic lattice}. As in the smooth case, one can define $\gamma$ as being represented by an $3$-tuple of sections $b\mapsto (\gamma_1(b), \gamma_2(b),\gamma_3(b))$, each one given by certain composition of Hamiltonian flows. In this case, however, $b\mapsto \gamma_2(b)$ does not vary smoothly but piecewise smoothly, failing to be smooth along $\Gamma$. The contribution of the path $\gamma_2\cap\mathfrak U$ to the periods $\lambda_2^{\pm}$ is $\lambda_0$. On the other hand, the contribution of $\gamma_2\cap X - \mathfrak U$ is $dH^\pm$. In contrast, the other two periods can be computed along paths entirely contained in $\mathfrak U$ which implies that they are smoothly defined on $B$.
\end{proof}

We will from now on denote $\lambda_1^{\pm}$ and $\lambda_{3}^{\pm}$ simply by $\lambda_1$ and 
$\lambda_3$ respectively.

\begin{rem} \label{cor:redper}  Notice that in the above we can assume $H^+|_{\Gamma} = H^-|_{\Gamma}$, therefore 
we can define $\bar \Lambda_H = \Lambda_{H^+} \mod db_1 = \Lambda_{H^-} \mod db_1 $. Via the 
identification in the above Proposition, the space $\bar Z^\#$ corresponds to 
$T^\ast\Gamma \slash \bar\Lambda_H$ and $\bar f^\#:\bar Z^\#\rightarrow \Gamma$ becomes the projection
 $\bar\pi^\#$.
\end{rem}

\medskip
We now introduce a standard model for fibrations satisfying Assumption~\ref{ass: semi-stitched}.

\begin{ex}[Normal form of cylindrical type]\label{cyl:norm}  
Let $(U,\Gamma)$ be a pair of subsets of $\R^2 \times \R$ diffeomorphic to 
$(D^2\times D^{1},D^1\times D^{1})$ with $\Gamma = U \cap \{b_1=0\}$. Let $\Delta=\{ b_1=b_2=0\}$. Given $H\in C^\infty (U)$ denote by $H_\Delta$ the germ of $H$ along $\Delta$. Consider the integral lattice $\Lambda_H$ in $T^\ast U$ generated by:
\begin{equation}\label{Lambda_H}
\begin{array}{l}
\lambda_1=2\pi db_1,\\
\lambda_2=dH+\arg (b_1+ib_2)db_1+\log |b_1+ib_2|db_2,\\
\lambda_3=db_3.
\end{array}
\end{equation}
Let $(y_1,y_2 ,y_3)$ denote the locally defined vertical coordinates on $T^* U$, which it is 
convenient to think of as $\Lambda_H$-periodic coordinates.  
For fixed positive $L\in\R$ consider the following subset of $T^*U$: 
\begin{equation}
C_L =\{ |y_2| < L \}
\end{equation} 
and denote $C_L(b) = T^*_b U \cap C_L$.
If $U$ is a sufficiently small neighborhood of $\Delta$, we can assume that for every $b \in U$, 
$2L < |\log |b|+\partial_{b_2}H|$. Therefore the projection $T^*_bU \rightarrow T^*_bU / \Lambda_H$ maps $C_L(b)$ to a cylinder which closes up in the $y_1$ and $y_3$ direction but not in the 
$y_2$ direction. So let us think of $C_L(b)$ as this cylinder and define 
$J^\circ_{L} =\bigsqcup_{b\in U}C_L (b)$, which is an open subset of $T^*_bU / \Lambda_H$. 
 The projection $\pi$ restricts to an open cylinder fibration:
\[
\pi^\circ:J^\circ_{L} \rightarrow U.
\]
Clearly there is an $S^1$ action on $J^\circ_{L}$ induced by $\lambda_1$, whose moment map 
is $b_1$. Let $Z^\circ_{L}=(\pi^\circ)^{-1}(\Gamma)$ and let $\bar{Z}^\circ_{L}$ be the corresponding $S^1$ 
reduced space. Let $\bar\pi^\circ:\bar{Z}_{L}^\circ \rightarrow \Gamma$ be the reduced fibration. We denote 
the fibre of $\bar\pi^\circ$ by $\bar C_L (b)$. 

For $L' < L$, construct $J^\circ_{L'}$, which is a cylinder fibration with shorter cylinders, and define 
its closure $K_{L'} = \overline{J^\circ_{L'}}$. 
Define the open set $E_{L,L'} = J^\circ_{L} - K_{L'}$, which we can think of as the union 
of the ends of the cylinders.  Suppose now that we have an
open neighborhood $V$ of $Z^\circ_{L}$ and a smooth $S^1$ invariant Lagrangian submersion $u: V \rightarrow 
\R^3$ with cylindrical fibres satisfying: $u|_{Z^\circ_{L}} = \pi^\circ$, $u|_{E_{L,L'}} = \pi^\circ$ and 
$u_1 = b_1$.  Then we can define $Y^+_{L}= (\pi^\circ)^{-1}(U^+)$, $Y_L = Y^+_L \cup V$, $Y^-_L= Y_L \cap (\pi^\circ)^{-1}(U^-)$ and the piecewise smooth function $f^{\circ}_{u}: Y_L \rightarrow B_u\subseteq\R^n$ to 
be the map
\begin{equation}
       f^\circ_{u} = \begin{cases}
               \pi^\circ \quad\text{on} \  Y^+_L, \\
	       u \quad\text{on} \ Y^-_L.
       \end{cases} 
\end{equation}
Clearly, if we think of $Y_L$ as playing the role of $X^\circ$, 
$f^\circ_{u}:Y_L \rightarrow B_u$ is a Lagrangian fibration of type (\ref{eq:f^c}). Notice that the fibres 
of $f^\circ_u$ coincide with the fibres of $\pi^\circ$ inside $E_{L,L'}$, in particular $f^\circ_u$ is smooth 
restricted to $E_{L,L'}$. In some sense, the fibres of $f_u$ are straight towards their ends 
(cf. Figure~\ref{J_0}). 

\medskip
We now compactify by adding the singularities. Let $J_H^\#=T^\ast U\slash\Lambda_H$ and 
let $\pi^\#:J^\#_H\rightarrow U$ be the Lagrangian fibration induced by the standard projection on 
$T^\ast U$. Clearly $J^\circ_{L}$ and therefore $Y_L$ are open subsets of $J_H^\#$.
When $b\in\Delta$, the fibre $C(b)=(\pi^\#)^{-1}(b)$ is an open cylinder, with ends at $+\infty$ and 
$-\infty$ in the $y_2$-direction, otherwise $C(b)$ is a torus. From the results in \cite{RCB1}, 
$J_H^\#$ can be compactified to a symplectic manifold $X$ by adding the singularity at the ends of the 
cylinders $C(b)$ when $b \in \Delta$. The fibration $\pi^{\#}$ extends to a smooth
fibration $f_{H}: X \rightarrow U$ of generic-singular type. The open subset $J^\#_H - K_{L'}$ extends to an open neighborhood 
$E$ of the singular set $\Sigma$. The fibres of $f^\circ_u$ coincide with the fibres of 
$f_{H}$ toward their ends and therefore $f^\circ_u$ may be extended to make it coincide with 
$f_{H}$ on  $E$. More precisely, define $\mathfrak U = f^{-1}_{H}(B_u) \cap E$ and 
$Y = Y_L \cup \mathfrak U$.
Now we can define 
\begin{equation} \label{u:st_o}
       f_{u,H} = \begin{cases}
                f_{H} \quad\text{on} \  \mathfrak U, \\
	       f^\circ_u \quad\text{on} \ Y_L.
       \end{cases} 
\end{equation}
Clearly $f_{u,H}: Y \rightarrow B_u$ is a well defined Lagrangian fibration satisfying 
Assumption~\ref{ass: semi-stitched}. The zero section $\sigma_0$ of $\pi^\circ$ is, perhaps after a change 
of coordinates in the base, a section of $f_u$.
If $F_{\bar b}$ is a smooth fibre of $f_{u,H}$, with $\bar b \in U^+$, let $\gamma_0$ be the 
basis of $H_1(F_{\bar b},\Z)$ determined by $\lambda_1, \lambda_2, \lambda_3$. 
 We call $\mathcal F_{u,H}=(Y, f_{H,u}, \sigma_0,\gamma_0)$ a \textit{normal form of cylindrical type.} 

\medskip
The set $Y_L\subset J^\#_H$ can be visualized in Figure \ref{J_0} as the square with open top and bottom. The straight light-colored lines are the fibres of $\pi^\#$ and the fibres of $f^{\circ}_u:Y_L\rightarrow B_u$
are depicted as dark lines. The upper and lower rectangular regions represent the components of $E_{L,L'}$.
\begin{figure}[!ht]
\begin{center}
\input{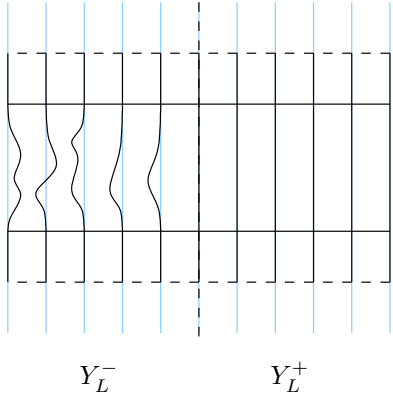}
\caption{Normal form of cylindrical type.}\label{J_0}
\end{center}
\end{figure}
\end{ex}
Given the above construction we denote $Z^{\#}_{H} = (\pi^{\#})^{-1}(\Gamma)$ and by $\bar Z^{\#}_{H}$
its $S^1$ quotient. Notice that if we let $\bar \Lambda_H = \Lambda_H \mod db_1 $, then 
$\bar Z^{\#}_{H} = T^{\ast}\Gamma / \bar \Lambda_H$. If $\bar \pi^{\#}$ is the projection, 
let $\mathfrak L =\ker\bar\pi^{\#}_\ast$. We can assume $u$ is a well defined map in a neighborhood of $Z^{\#}_{H}$ 
which coincides with the projection outside a neighborhood of $Z^{\circ}_L$, therefore we can associate to the 
pair $(V,u)$ a sequence $\ell = \{ \ell_k \}_{ k \in N}$ of fibrewise closed section of $\mathfrak L^*$, just as we did 
in the proper case. We can easily see that the sequence $\ell$ must vanish outside $\bar Z ^{\circ}_L$, in particular 
each $\ell_k$, when restricted to a fibre, has compact support contained in the cylinder $\bar C_L(b)$. 
With respect to the proper case, in this situation we have an additional piece of data, i.e. the smooth 
function $H$.

\medskip
The following is analogous to Definition \ref{def:seq}:

\begin{defi}\label{def:open_inv} With the above notation, 
\begin{itemize}
\item[i)] Let $\mathscr L_{\bar Z^{\#}_{H}}$ the set of sequences of fibrewise closed sections of 
$\mathfrak L^\ast$ which vanish outside $\bar Z^{\circ}_L$ for some positive $L$ such that 
$2L < |\log |b|+\partial_{b_2}H|$ for every $b \in \Gamma$.
\item[ii)] Let $\mathscr U_{\bar Z^{\#}_{H}}$ be the set of pairs $(V,u)$ where, for some positive $L$ and $L'$ 
satisfying $2L'< 2L < |\log |b|+\partial_{b_2}H|$, $V$ is a neighborhood of $Z^\circ_L$ 
and $u:V\rightarrow\R^{n}$ is a smooth,  $S^1$-invariant Lagrangian submersion, with cylindrical 
fibres, with components $(u_1, u_2 ,u_3)$ such that 
$u|_{Z^\circ_L}=\pi^\circ$, $u|_{E_{L,L'}}=\pi^\circ$ and $u_1=b_1$.
\item[iii)] Let $\mathscr H_{\Delta}$ be the set of germs $H_{\Delta}$ of smooth functions $H$ 
defined on neighborhoods of  $\Delta$.
\end{itemize}
\end{defi}

We define the invariants of a normal form of  cylindrical type $\mathcal F_{u, H}$ to be:
\[
\inv (\mathcal F_{u, H}) = (Z^{\#}_{H}, \ell, H_{\Delta}).
\] 

A little explanation is necessary to see in which sense these are invariants.
\begin{rem} \label{inv:cyl} 
Suppose we are given two normal forms of cylindrical type
$\mathcal F_{u, H}$ and $\mathcal F_{u', H'}$. From the results in $\cite{RCB1}$ 
(cf. also Theorem~\ref{thm. main RCB}), a necessary condition for $f_H$ and $f_{H'}$ to be symplectically conjugate is that $H_{\Delta} = H'_{\Delta}$, so suppose this holds. This gives a symplectomorphism, 
which we denote by $\Phi_{H, H'}$, between the total spaces $X$ and $X'$ of the two fibrations which conjugates $(X, f_H, B)$ and $(X', f_{H'}, B')$. By pulling back 
$(V', u')$ via this symplectomorphism and computing the Taylor series, we obtain a sequence of fibrewise closed 
sections of $\mathfrak L^\ast$ which we call $\Phi_{H,H'} \cdot \ell'$. Using the same arguments 
as in the proof of Theorem~\ref{thm: grosso} (cf.\cite{CB-M-stitched}, Theorem 6.11), we can then show 
that $\mathcal F_{u, H}$ and $\mathcal F_{u', H'}$ are symplectically conjugate 
if and only if $\Phi_{H,H'} \cdot \ell' = \ell$. In particular, when $H = H'$, they are symplectically conjugate 
if and only if $\ell = \ell'$.
\end{rem} 

For the classification of fibrations satisfying Assumption \ref{ass: semi-stitched}, it is 
useful to have the following result.
\begin{prop}\label{prop:ext} 
Let $f:X\rightarrow B$ be a Lagrangian fibration satisfying Assumption \ref{ass: semi-stitched}. 
Given a smooth fibre $F_{\bar b}$ of $f$ there is a basis $\gamma$ of $H_1(F_{\bar b},\Z)$ and a 
section $\sigma$ of 
$f$, such that $\mathcal F=(X, f, B, \sigma, \gamma)$ is symplectically conjugate to a normal 
form of cylindrical type $\mathcal F_{u, H}$. 
\end{prop}

\begin{proof} 
One uses the same arguments as in the proof of Proposition \ref{prop:normalform}. Suppose there is an extension of
$f^+:X^+\rightarrow B^+$ to a smooth Lagrangian fibration $\tilde f^+$ defined on a neighborhood $W\subseteq X$ of $Z$
such that $\tilde f^{+}|_{\mathfrak U} = f|_{\mathfrak U}$. 
Then one may compute the period lattice of $\tilde f^+$; this gives a smooth function $H$ 
extending the function $H^+$ in Proposition \ref{prop:periods_sh}. Assuming 
that also $f^-$ has been extended to $\tilde{f}^{-}$ so that $\tilde f^{-}|_{\mathfrak U} = f|_{\mathfrak U}$, one may verify that 
the period map $\Theta^+ :T^\ast U\slash \Lambda_H\rightarrow W^\#$ gives the required equivalence between 
$\mathcal F$ and $\mathcal F_{u, H}$ where $u=\tilde f^-\circ \Theta^+$.

To extend $f^+$, notice that $f_{\mathfrak U}=f|_{\mathfrak U}$ is smooth so, tautologically, $f_{\mathfrak U}$ is an extension of $f^+$ to $\mathfrak U$. It remains to extend $f^+$ away from $\mathfrak U$. Let $\mathfrak U'\subset\mathfrak U$ and define $f^\circ :X^\circ\rightarrow B$  as in (\ref{eq:f^c}). Denote $Z^{\circ} = Z \cap X^\circ$ and by $\bar Z^\circ$ its $S^1$ quotient with 
$\bar f^\circ:\bar Z^\circ\rightarrow\Gamma$ the reduced fibration.
Then $\bar f^\circ$ is a smooth Lagrangian cylinder fibration.

The coisotropic neighborhood theorem allows us to identify a neighborhood of $Z^\circ$ inside $X^\circ$ 
with a neighborhood $V$ of $\{0\} \times S^1 \times\bar Z^\circ$ inside $\R \times S^1 \times \bar Z^\circ$ 
($t$ will denote the $\R$ coordinate). 
Moreover, since $\bar Z^{\#}$ can be identified with $T^{\ast}\Gamma / \bar \Lambda_H$ (see Remark~\ref{cor:redper}), 
$\bar Z^\circ$ can be identified with a subset of $T^{\ast}\Gamma / \bar \Lambda_H$ of the type $\bar Z^\circ _L$ 
for some positive $L$ (see Example~\ref{cyl:norm}).
The pullback of $f^\circ$ under these identifications gives a piecewise smooth Lagrangian fibration on 
$V \subset \R \times S^1 \times \bar Z^\circ_L$
\begin{equation}
       g = \begin{cases}
       		u^+ \quad\text{on} \  V^+;\\
                u^- \quad\text{on} \  V^-
       \end{cases} 
\end{equation}
where $V^+ = V \cap \{ t \geq 0 \}$, $V^- = V \cap \{ t \leq 0 \}$ and $u^\pm$ is the restriction to $ V^\pm$ of a $C^\infty$ map.
The set $Z^\circ \cap \mathfrak U$ where $f^\circ$ is smooth, corresponds (under the above identifications) to the interior of $Z^\circ_L - Z^\circ_{L'}$ which we denote $C_{L, L'}$, 
where $L'  < L$. Notice that the map $g$ above is then smooth along $C_{L,L'}$, in particular the 
Taylor expansions in $t$ of $u^+$ and $u^-$ coincide along $C_{L,L'}$. With the same arguments 
used in the proper case one can show that $u^\pm$ can be smoothly extended to a Lagrangian fibration $\tilde u^\pm$ beyond $V^\pm$ 
(cf. Proposition~\ref{prop:normalform} above, or \cite{CB-M-stitched} Proposition 6.3 for more details). In fact with a little 
more care one can do this so that along $\R \times C_{L, L'}$, where an extension already exists, namely $g$ itself, we have 
$\tilde u^{\pm}|_{\R \times C_{L, L'}} = g|_{\R \times C_{L, L'}}$. The map $\tilde u^+$ gives the required extension $\tilde f^+$ of $f^+$, 
where the last observation guarantees that $\tilde f^{+}|_{\mathfrak U} = f|_{\mathfrak U}$.
\end{proof}

From the above result, it follows that to every Lagrangian fibration $\mathcal F$ satisfying Assumption~\ref{ass: semi-stitched} we can assign the invariants of a normal form for $\mathcal F$, 
i.e. a triple $(Z^{\#}_{H}, \ell, H_{\Delta})$. Notice that two normal forms $\mathcal F_{u, H}$ and 
$\mathcal F_{u', H'}$ for the same fibration $\mathcal F$ must be related in the way described 
in Remark~\ref{inv:cyl}. 
It is worth stating this in the following:
\begin{thm} \label{thm:stitched_open}
Two germs of fibrations $\mathcal F$ and $\mathcal F'$ satisfying Assumption~\ref{ass: semi-stitched} are 
symplectically conjugate if and only if their invariants are related in the way described in Remark~\ref{inv:cyl}. \todo{I am not completely satisfied with the statement of this theorem, I think 
it could be a bit more precise. But let's leave it like this for the moment.}
\end{thm}

We also have:
\begin{prop} \label{cyl:exist_nor}
Given $H_{\Delta} \in \mathscr H_{\Delta}$, there is a function $H$ defined on a neighborhood of 
$\Gamma$ whose germ is $H_\Delta$, such that for every $\ell \in \mathscr L_{\bar Z^{\#}_H}$,
there is a normal form of cylindrical type whose invariants are $(Z^{\#}_H, \ell, H_{\Delta})$.
\end{prop}

The results in this section extend those in \cite{RCB1} to stitched fibrations with generic singularities (satisfying Assumption \ref{ass: semi-stitched}). The arguments here can also be carried through in the stitched focus-focus case, the positive case and their higher dimensional analogues.

\section{Lagrangian negative fibrations} \label{negative} 
The purpose of this section is two-fold. We first use the analysis in \S \ref{stitched fibr} to refine the piecewise smooth fibrations constructed in \S \ref{sec:pwfibr}. Subsequently, we study the affine structures associated to the resulting fibrations. 

\medskip
Recall that we defined a negative vertex to be an integral affine manifold with singularities
modeled on Example \ref{neg:affvar}.

\begin{defi} \label{lag:neg} 
Let $(X, \omega)$ be a $6$-dimensional symplectic manifold and $B \subseteq \R^3$ an open subset. Let $f: X \rightarrow B$ be a piecewise smooth Lagrangian fibration. $\mathcal F=(X, \omega, f, B)$ is called a \textit{Lagrangian negative fibration} if it satisfies the following properties:
\begin{itemize}
\item[(i)] $\mathcal F$ is topologically conjugate to the alternative negative fibration of Example~\ref{ex. alt (2,1)}.
\item[(ii)] there exists a submanifold with boundary $D \subset B$, homeomorphic to a closed disc 
in $\R^2$, such that $\Delta \cap (B - D)$ consists of three one dimensional disjoint segments (the legs of $\Delta$) and $f$ is smooth when restricted to $X - f^{-1}(D)$;
\item[(iii)] let $B_0= B -(D\cup \Delta)$,  $X_0=f^{-1}(B_0)$ and $f_0=f|_{X_0}$. Let $(B_0, \mathscr A)$ be the integral affine manifold induced by the Lagrangian $T^3$ bundle $\mathcal F_0=(X_0, f_0, B_0)$. For some choice of model of negative vertex 
$(\R^3, \Delta_{\tau}, \mathscr A_{\tau})$ as given in Example~\ref{neg:affvar}, there 
exist an open neighborhood $U\subseteq\R^3$ of $0\in\R^3$, a 
submanifold with boundary $D' \subset U$ homeomorphic to a closed disc in $\R^2$, satisfying 
$0 \in D' \subset \{ x_1 = 0\} \subset\R^3$ and an integral affine isomorphism
\[
 (B_0, \mathscr A)\cong(U - (D' \cup \Delta_\tau), \mathscr A_\tau).
\] 
\end{itemize}
\end{defi}

Corollary \ref{aff:symp} directly implies the following:

\begin{prop}\label{prop:neg_glue} 
Let $\mathcal F$ be a Lagrangian negative fibration. With the notation as in 
Definition~\ref{lag:neg}, let $U_0=U - (D' \cup \Delta_\tau)$ and $X(U_0, \mathscr A_{\tau})$ 
be the associated Lagrangian torus bundle. Then, if $\mathcal F$ has a smooth Lagrangian section, 
$\mathcal F_0$ is symplectically conjugate to $X(U_0, \mathscr A_{\tau})$.
\end{prop}

The main result of this section is
\begin{thm}\label{thm:neg_main} There exists a symplectic manifold $(X, \omega)$ and a map $f:X \rightarrow B$ such that $(X,\omega, f, B)$ is a Lagrangian negative fibration. 
\end{thm}

The starting point aiming at the proof of Theorem \ref{thm:neg_main} 
will be the Lagrangian fibration described in 
Example~\ref{thin leg}, which satisfies Definition~\ref{lag:neg}(i). 
The proof will consist essentially of 
three steps. First we modify Example~\ref{thin leg} so to obtain a 
fibration which is smooth towards the ends of the $1$-dimensional 
legs (Smoothing I and II). In the second step (Smoothing III) we use the invariants of stitched Lagrangian 
fibrations to modify the fibration once more so that it satisfies
property (ii). Finally we show that these modifications have been done
in a way that also (iii) holds. 
 
\subsection{Smoothing I}
Let us consider the fibration as in Example~\ref{thin leg} with its discriminant locus $\Delta$. Recall that this fibration is constructed using Proposition~\ref{prop. piecewise smoothness}, 
by taking as symplectomorphism $\Phi$ the one described by (\ref{thin:fi}).
For positive  $M \in \R$ let us define
\begin{equation}\label{far_legs}
\Delta_{h,M} = \Delta \cap \{ b_2 \leq - M \}, \quad 
             \Delta_{v,M} = \Delta \cap \{ b_3 \leq - M \}, \quad 
                      \Delta_{d,M} = \Delta \cap \{ b_2, b_3 \geq  M \}.
\end{equation}

When $M$ is sufficiently big, $\Delta_{h,M}$, $\Delta_{v,M}$ and 
$\Delta_{d,M}$ are $1$-dimensional. In fact, they are the ends of the horizontal,
vertical and diagonal legs of $\Delta$ respectively. Now let $\Sigma_{h,M}$, 
$\Sigma_{v,M}$ and $\Sigma_{d,M}$ be the parts of the critical surface $\Sigma$ which 
are mapped to $\Delta_{h,M}$, $\Delta_{v,M}$ and $\Delta_{d,M}$ respectively.

\medskip
We have the following
\begin{lem} \label{lem:smth_crit} 
The piecewise smooth Lagrangian fibration $\mathcal F=(X,\omega, f, B)$ in 
Example~\ref{thin leg} can be perturbed, without changing its topology,
so that, for sufficiently big $M$, it becomes smooth on small neighborhoods $N_{h,M}$, $N_{v,M}$ and 
$N_{d,M}$ of $\Sigma_{h,M}$, $\Sigma_{v,M}$
and $\Sigma_{d,M}$ respectively. 
\end{lem}
\begin{proof} From the way $f$ is defined in Example~\ref{thin leg}, we can 
assume 
\[ \Sigma_{h,M} = \{ t= 0, \ u_1 = 0, \ |u_2|^2 < \epsilon/4 \}, \]
where $\epsilon$ is as in (\ref{thin:fi}) and $M = \log (\sqrt {\epsilon} / 2 )$. 
For any $\tau > 0$ denote open sets
\[ N^{\tau} = \{(t,u_1, u_2) \ | 
         \  \max(|u_1|, |u_2|^2) < \tau \}. \]
From now on we assume $f$ is restricted to $N^{\epsilon/2}$.
As one can easily see from the construction, the map 
$G_t$ defining $f$, restricted to $N^{\epsilon/2}$ is
\begin{equation} \label{Gt}
    G_t(u_1, u_2) = \left( 
        \log |u_2|, \log \left| \frac{u_1}{\sqrt{|t|+\sqrt{t^2+|u_1|^2}}}-1 \right|
                                                     \right). 
\end{equation}
This is the map that we want to perturb, but just on a smaller neighborhood. 
We do it applying the idea already anticipated at the end of Example~\ref{leg}.
In fact we notice that $G_t$ is invariant with respect to the $S^1$ action 
\[ e^{i\theta}(u_1, u_2) = (u_1, e^{2i\theta}u_2), \]
which is also Hamiltonian with respect to the reduced symplectic form 
$\omega_t$ given in (\ref{om:red}). The moment map is 
\[ (u_1, u_2) \mapsto |u_2|^2. \]
So, if $g$ is a real function depending only on $u_1, t$ and $s=|u_2|^2$,
then $$(u_1, u_2) \mapsto ( \log|u_2|, \ g(u_1, t, |u_2|^2))$$ is a Lagrangian 
fibration with respect to $\omega_t$, provided the level sets of
$u_1 \mapsto g(u_1, t,s)$ are one dimensional submanifolds for every 
$s$ and $t$. For example, consider a real non-negative function 
$\rho$ defined on $\R^3$ such that, for every fixed $(t, s) \in \R^2$, the map
\begin{equation} \label{resc:ro}
u \mapsto \frac{u}{\rho(|u|^2, t,s)}
\end{equation}
is a local homeomorphism of a neighborhood of $u=0$, 
then $g= \log | \frac{u}{\rho} - 1|$ defines a Lagrangian fibration
(at least in a neighborhood of $0$).
In particular  
\[
\rho_0(r,t) = \sqrt{|t|+\sqrt{t^2+ r}},
\] 
with $(r,t) \in \R^2$ gives the map $G_t$ in (\ref{Gt}), but it is not smooth. It is easy to see 
that if $\rho$ is smooth on $\R^3$ and satisfies
\begin{equation} \label{ro1:ro0}
           \rho > \rho_0 
\end{equation}
then the map (\ref{resc:ro}) is an orientation preserving diffeomorphism (at least 
near $u = 0$). 
So let us choose a smooth $\rho_1$, defined on $\R^2$ and satisfying $\rho_1 > \rho_0$, 
and let
\[ g_j = \log \left| \frac{u_1}{\rho_j(|u_1|^2,t)} - 1 \right|, \]
for $j=0,1$. We wish to find a $g$ which interpolates between $g_0$ and $g_1$. 
More precisely, we want $g$ to be equal to $g_0$ outside $N^{3\epsilon/8}$ and 
to $g_1$ on some smaller open neighborhood of $\Sigma_{h,M}$.
Clearly $(u_1, u_2) \in N^{3\epsilon/8}$ if and only if $(|u_1|^2, |u_2|^2)$ is
in the rectangle $$S_0= [- 9 \epsilon^2 / 64 , 9 \epsilon^2 / 64] \times 
 [ - 3 \epsilon /8, 3 \epsilon /8].$$ Now let $S_1$ be a closed neighborhood 
of $0$ in $\R^2$ which is contained in the interior of $S_0$, e.g. a smaller 
rectangle. Taking a $\sigma \in C^{\infty}(\R^2)$, which is $0$ outside $S_0$ and $1$ on $S_1$,
let us define
\[ \rho(r, t, s) = 
            (1 - \sigma(r, s)) \rho_0(r, t) + 
          \sigma(r, s ) \rho_1(r, t),\]
so that $\rho$ is equal to $\rho_0$ outside 
$S_0$ and it is equal to $\rho_1$ on $S_1$. Clearly $\rho > \rho_0$. 
We leave it to the reader to check that choices can be made so that with this 
$\rho$, (\ref{resc:ro}) is indeed a homeomorphism.
Now define
\[ g = \log \left| \frac{u_1}{\rho(|u_1|^2,t, |u_2|^2)} - 1 \right|. \]
Clearly $g$ is equal to $g_0$ outside $N^{3\epsilon/8}$ and to
$g_1$ on 
\[ N_{h,M} = \{ (|u_1|^2, |u_2|^2) \in S_1 \} \]
which, with a suitable choice of $S_1$, is a neighborhood of $\Sigma_{h,M}$. 
Moreover $u \mapsto g(u,t,s)$ has $1$-dimensional level sets. We can therefore replace the 
second component of $G_t$ in (\ref{Gt}) with $g$ and redefine
\[ G_t(u_1, u_2) = ( \log |u_2|, g), \]
which is smooth on $N_{h,M}$. This proves the lemma for $\Sigma_{h,M}$. A schematic picture of this smoothing is described in Figure \ref{fig:smoothing_I_II}. The vertical lines represent fibres of $f$ over the horizontal leg. The base of the fibration is represented by the horizontal line on the bottom of the picture; the bold segment on the right represents the region where the codimension one part of $\Delta$ begins. The shaded region represents the locus where $f$ is not smooth. The dashed region is $N_{h,M}$.
\begin{figure}[!ht] 
\begin{center}
\input{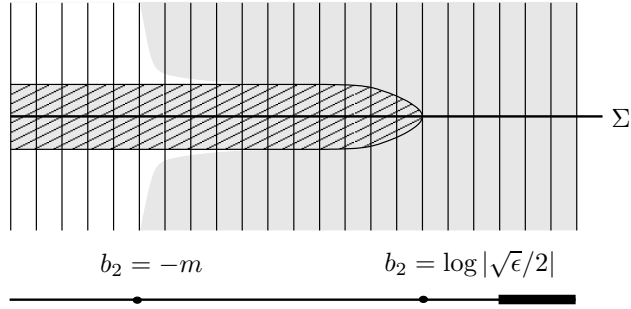}
\caption{Horizontal leg. The dashed region is $N_{h,M}$ as in Lemma \ref{lem:smth_crit}. After Smoothing II there will be a full fibred neighborhood (white region) where the fibration is smooth.}\label{fig:smoothing_I_II}
\end{center}
\end{figure}

The case of the vertical leg is done in the same way.
At first sight it is not so obvious that also the diagonal leg can 
be treated in the same way. So let us give some explanation. 
When $|u_2|^2 \geq M$, the map $G_t$ becomes
\begin{equation} \label{Gtd}
  G_t(u_1, u_2) = \left( 
    \log \left| \frac{u_1}{\rho_0(|u_1|^2,t)}-u_2 \right|, \,
 \log \left| \frac{u_1}{\rho_0(|u_1|^2,t)} + u_2 \right|  \right). 
\end{equation}
The first observation is that this map is invariant under the 
$S^1$-action
\begin{equation} \label{dia:s}
     e^{i\theta}(u_1, u_2) = (e^{i\theta}u_1, e^{i\theta}u_2). 
\end{equation}
After the following change of coordinates on the base
\[ (x_1, x_2) \mapsto \left( \frac{e^{2x_1} + e^{2x_2}}{2}, x_1 - x_2\right) \]
this becomes
\begin{equation} \label{Gtd2}
  G_t(u_1, u_2) = \left( \frac{\sqrt{t^2 + |u_1|^2} + |u_2|^2}{2} -\frac{|t|}{2}, \,
    \log \frac{\left| u_1 /\rho_0 - u_2 \right|}
              {\left| u_1/ \rho_0 + u_2 \right|}  \right). 
\end{equation}
One can check that for every fixed $t \in \R$ the map 
\[ (u_1, u_2) \mapsto \frac{\sqrt{t^2 + |u_1|^2} + |u_2|^2}{2}, \]
is the moment map of the $S^1$-action (\ref{dia:s}), with respect to
the reduced symplectic form $\omega_t$. Moreover, if one replaces 
$u_1 = z_1z_2$, $u_2 = z_3$ and $t = \frac{|z_1|^2 - |z_2|^2}{2}$, 
then the above map becomes 
\[ \nu: 
   (z_1, z_2, z_3) \mapsto  \frac{|z_1|^2 +|z_2|^2}{4} + \frac{|z_3|^2}{2} \]
which is a smooth map on the total space. Let us denote
\[ s = \frac{\sqrt{t^2 + |u_1|^2} + |u_2|^2}{2}. \]

The second component of (\ref{Gtd2}) can be rewritten as
\[ g_0(u_1, u_2) =  \log \left| \frac{ 2u_1 /\rho_0}
              { u_1/ \rho_0 + u_2} - 1  \right|.\]
We can now apply the same strategy we used in the case of the 
horizontal leg. We observe that we could replace
this $g_0$ with any other $S^1$-invariant function $g$. In particular
we could replace $\rho_0$, which is $S^1$-invariant, with another
smooth $S^1$-invariant $\rho_1$.
As before, we then interpolate $\rho_0$ and $\rho_1$ 
with a cut off function $\sigma$ depending on $|u_1|^2$ and $s$.
We avoid writing the details here, as they just follow the same argument
as before. 

In the end we obtain that, in a small neighborhood of $\Sigma_{d,M}$,
$G_t$ can be written as:
\[ G_t = \left( s - \frac{|t|}{2}, \log \left| \frac{ 2u_1 /\rho_1}
              { u_1/ \rho_1 + u_2} - 1  \right| \right), \]
where now the second component is smooth. The first component is 
not quite smooth yet. We saw that $s$ is smooth when lifted to 
the total space, but $|t|$ isn't. The total fibration becomes of 
the type
\[ f(z_1, z_2, z_3) = \left( \mu, \nu - \frac{|\mu|}{2}, g(z_1z_2, z_3, \mu, \nu)\right), \] 
where $g$ is smooth. We see that after a change of coordinates
on the base of the type
\begin{equation} \label{chcoor:dia}
   (b_1, b_2, b_3) \mapsto \left(b_1, b_2 + \frac{|b_1|}{2}, b_3\right) 
\end{equation} 
this fibration becomes
\[ f(z_1, z_2, z_3) = \left( \mu, \nu, g(z_1z_2, z_3, \mu, \nu)\right), \] 
which is smooth. One can find a global change of coordinates on 
the base which acts like (\ref{chcoor:dia}) only in a neighborhood
of the end of the diagonal leg and is the identity elsewhere. 
This ends the proof of the Lemma.
\end{proof}

\begin{rem} Notice that the new perturbed fibration of Lemma~\ref{lem:smth_crit} 
has a Lagrangian section. In fact one can easily see that the 
section of the fibration in Example~\ref{thin leg} survives the smoothing 
above, since it is far from the critical surface $\Sigma$.
\end{rem}

\subsection{Smoothing II}

\begin{figure}[!ht] 
\begin{center}
\input{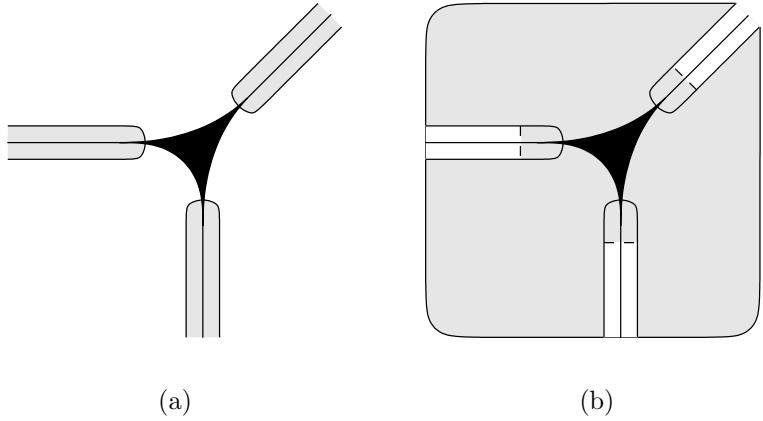}
\caption{Smoothing over the legs.}\label{fig:smoothing_1}
\end{center}
\end{figure}

Lemma \ref{lem:smth_crit} gives us a piecewise smooth fibration $\mathcal F$, topologically conjugate to the one in Example \ref{thin leg} but smooth along $N_{h,M}$, $N_{v,M}$ and $N_{d,M}$. The latter are sets mapping down onto open neighborhoods $B_{h,M}$, $B_{v,M}$ and $B_{d,M}$ of the legs as depicted in Figure \ref{fig:smoothing_1} (a).  Given a positive $m \in \R$, let us denote by $B_{h, m}$, 
$B_{v,m}$ and $B_{d,m}$ neighborhoods of $\Delta_{h,m}$, $\Delta_{v,m}$ and $\Delta_{d,m}$ and for brevity let us define 
$\mathcal F_{h, m} = \mathcal F|_{B_{h,m}}$, $\mathcal F_{v, m} = \mathcal F|_{B_{v,m}}$ and
$\mathcal F_{d, m} = \mathcal F|_{B_{d,m}}$. Clearly when $M$ is as in Lemma~\ref{lem:smth_crit}, $\mathcal F_{h, M}$, $\mathcal F_{v, M}$ and $\mathcal F_{d, M}$ 
satisfy Assumption \ref{ass: semi-stitched}.

\medskip
Our goal now is to use the results on non-proper stitched fibrations in Section~\ref{stitched fibr} to perturb $\mathcal F$
so that for some $m > M$ and neighborhoods $B_{h, m}$, $B_{v,m}$ and $B_{d,m}$, the fibrations 
$\mathcal F_{h, m}$, $\mathcal F_{v, m}$ and 
$\mathcal F_{d, m}$ are smooth. This will produce a fibration whose base is depicted in Figure \ref{fig:smoothing_1} (b). Over the white rectangular regions the fibration is completely smooth but on the shaded region it is still piecewise smooth. The result is the 
following:

\begin{lem}\label{lem:smth_crit2} Let $\mathcal F$ denote the fibration obtained in Lemma \ref{lem:smth_crit}. 
Given a positive real number $m > M$, there exists a perturbation $\tilde{\mathcal F}$ of $\mathcal F$ (perhaps defined over a 
smaller neighborhood of the plane $\{ b_1 = 0 \}$), such that 
\begin{itemize}
\item[(i)] $\tilde{\mathcal F}$ is topologically conjugate to $\mathcal F$;
\item[(ii)] there are open neighborhoods $B_{h, m}$, 
$B_{v,m}$ and $B_{d,m}$ of $\Delta_{h,m}$, $\Delta_{v,m}$ and $\Delta_{d,m}$ respectively so that 
the fibrations $\tilde{\mathcal F}_{h, m}$, $\tilde{\mathcal F}_{v, m}$ and $\tilde{\mathcal F}_{d, m}$ are smooth. 
\end{itemize}
\end{lem}
\begin{proof}
Consider one of the fibrations $\mathcal F_{h,M}$, $\mathcal F_{v,M}$ or $\mathcal F_{d,M}$ as above (whenever necessary, we allow ourselves to restrict to smaller neighborhoods of $\Delta_{h,M}$, $\Delta_{v, M}$ 
or $\Delta_{d, M}$). To keep the notation simple we temporarily drop the subindices and denote it by $\mathcal F$. 

Since $\mathcal F$ satisfies Assumption~\ref{ass: semi-stitched}, it follows from Proposition~\ref{prop:ext} that we can associate 
to $\mathcal F$ 
a normal form of cylindrical type $\mathcal F_{u, H}$ together with its invariants given by a triple $(Z^{\#}_{H}, \ell, H_\Delta )$ which, in view of Theorem~\ref{thm:stitched_open}, uniquely determine $\mathcal F$ as a germ around $\Gamma = B \cap \{ b_1 = 0\}$. By slight abuse of notation we will denote by the same letter $\Gamma$ both $B \cap \{ b_1 = 0\}$ and $B_u \cap \{b_1 = 0 \}$, where $B_u$ is the base of $\mathcal F_{u, H}$. For the duration of this proof $H$ will remain unchanged, so we drop the subindex $H$ and denote $\mathcal F_u:=\mathcal F_{u,H}$ for short.

\begin{figure}[!ht] 
\begin{center}
\input{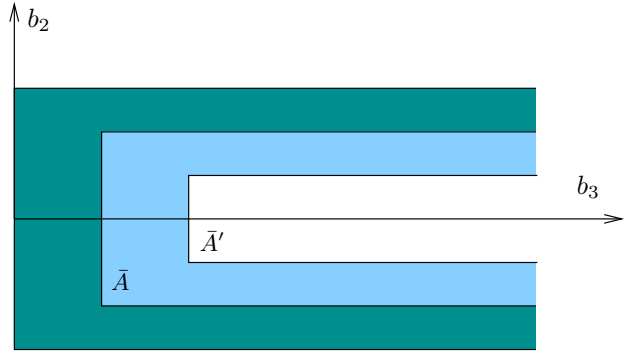}
\caption{$\Gamma$ (or $\Gamma_{h,M}$).}\label{nhbd_leg} 
\end{center}
\end{figure}

\medskip
The proof consists in suitably deforming the sequence $\ell$. 
Let $\bar A\subset\Gamma $ and $\bar A'\subset\bar A$  be (planar) regions as depicted in Figure \ref{nhbd_leg}. 
Given a cut-off function $\rho \in C^\infty (\Gamma)$ such that $\rho$ is 1 on $\Gamma - \bar A$ and $0$ on $\bar A'$, 
define a new (fibrewise closed) sequence $\tilde \ell$ whose elements are $\tilde\ell_k=(\rho \circ \bar\pi^{\#}) \, \ell_k$ for each 
$k\in\N$. We obtain a triple $(Z^{\#}_{H}, \tilde\ell, H_\Delta )$, such that 
$\ell|_{(\bar\pi^{\#})^{-1} (\Gamma - \bar A)} = \tilde\ell|_{(\bar\pi^{\#})^{-1}(\Gamma - \bar A)}$ and 
$\tilde\ell|_{(\bar\pi^{\#})^{-1} (\bar A')} = 0$.

\medskip
In view of Proposition \ref{cyl:exist_nor}, $(Z^{\#}_{H},\tilde \ell, H_\Delta )$ gives rise to a normal form of 
cylindrical type $\mathcal F_{\tilde u}$ defined over a neighborhood of $\Gamma$. By construction and by Theorem~\ref{thm:stitched_open}, $\mathcal F_{u}$ and $\mathcal F_{\tilde u}$ define 
the same germ around $\Gamma - \bar A$, i.e. there are open neighborhoods $U$ and $\tilde U$ of $\Gamma - \bar A$ 
(satisfying $U \cap \{ b_1 = 0 \} = \tilde U \cap \{ b_1 = 0 \} = \Gamma - \bar A$) such that 
$\mathcal F_{u}|_{U}$ and $\mathcal F_{\tilde u}|_{\tilde U}$ are symplectically conjugate. Moreover 
$\mathcal F_{\tilde u}$ is smooth when restricted to any open neighborhood $A'$ of $\bar A'$ such that 
$A' \cap \{ b_1 = 0 \} = \bar A'$. Now recall that $\mathcal F_{u}$ is symplectically conjugate to $\mathcal F$, so 
we have that $\mathcal F_{\tilde u}|_{\tilde U}$ is symplectically conjugate $\mathcal F|_{U}$. 

\medskip
Let us summarize the result using our original notation for the horizontal leg. 
For $\Gamma_{h, M} =  B_{h, M} \cap \{ b_1 = 0 \}$,  we have found sets $\bar A' \subset \bar A \subset \Gamma_{h, M}$ 
(as in Figure~\ref{nhbd_leg}) and a normal form of cylindrical type $\mathcal F_{\tilde u}$, defined over a 
neighborhood of $\Gamma_{h, M}$, smooth over $\bar A'$ and such that 
$\mathcal F_{\tilde u}|_{\tilde U}$ is symplectically conjugate to $\mathcal F_{h, M}|_{U}$,  where $U$ and $\tilde U$ are
neighborhoods of $\Gamma_{h, M} - \bar A$ (satisfying $U \cap \{ b_1 = 0 \} = \tilde U \cap \{ b_1 = 0 \} = \Gamma_{h,M} - \bar A$).

\medskip
If we go back denoting by $\mathcal F$ the fibration of Lemma~\ref{lem:smth_crit}, we can form a new fibration $\tilde{\mathcal F}$ 
in the following way. Let $\mathcal F' = \mathcal F|_{\R^{3} - (\R \times \bar A)}$ and symplectically glue 
$\mathcal F_{\tilde u}$ to $\mathcal F'$ using the conjugation between $\mathcal F_{\tilde u}|_{\tilde U}$ and 
$\mathcal F'|_{U} = \mathcal F_{h, M}|_{U}$. The fibration $\tilde{\mathcal F}$ is the result of this gluing. Notice that 
$\tilde{ \mathcal F}$, due to the properties of $\mathcal F_{\tilde u}$, is such that for some $m > M$ (depending on $\bar A'$) 
and a suitable neighborhood of $B_{h,m}$ of $\Delta_{h,m}$, the restriction $\tilde{\mathcal F}_{h, m}$ is smooth. Notice that $\bar A'$ 
can be chosen so that the latter holds for any $m >M$.

The above method applied to all legs, produces the required result.
\end{proof}

The idea of deforming the sequence $\ell$ by multiplying it by a cut-off function on the base will be used again in the subsection Smoothing III. This is actually the main application of the results on stitched fibrations in this paper.

\begin{rem} \label{section:surv}
We observe that the Lagrangian section of Example~\ref{thin leg} 
survives also this second smoothing.
\end{rem}
\subsection{The normal form}
Consider the Lagrangian fibration $\mathcal F$ produced in Lemma~\ref{lem:smth_crit2}.
If we let $U = \R^3 - \Delta$, then $\mathcal F|_U$ is a stitched 
$T^3$ fibration whose seam consists of three disjoint components.
It is clear that $\mathcal F|_{U}$ is a fibration of the type described in Example~\ref{amoeb:mon}. The goal of this section is to show that $\mathcal F|_{U}$  is in fact symplectically conjugate to a fibration which can be  constructed with 
Theorem~\ref{stitch:monodr_neg}, maybe after restricting the latter to a smaller neighborhood of the vertex 
of $\Delta$ (see Remarks~\ref{skew_ff:ex} and \ref{skew_neg:re}). Essentially, we need to show that the 
action coordinates, a priori defined only on a contractible open set, extend continuously to $\R^3$. 
We need the following
\begin{lem} \label{neg:ex} 
Let $(X, \omega)$ be the total space of the fibration produced in Lemma~\ref{lem:smth_crit2}. 
Then $\omega$ is exact on $X$. 
\end{lem}
\begin{proof}
Recall that the fibration produced in Lemma~\ref{lem:smth_crit2} is a perturbation of the
one in Example~\ref{thin leg}, whose total space is an open set of $\C^3$ with 
standard symplectic form, which is exact. One can see that the successive perturbations of 
this fibration have not modified the cohomology class of $\omega$. 
\end{proof}

To describe the fibration $\mathcal F$ we use the same notation of Example~\ref{amoeb:mon}.
Given $\bar b \in \Gamma_c$, there exists a basis $\gamma = \{ \gamma_1, \gamma_2, \gamma_3 \}$ of 
$H_1(F_{\bar b}, \Z)$ with respect to which monodromy is generated by the matrices in (\ref{mon:m})
with $m_1 = m_2 = 1$. 
We can compute the action coordinates $\alpha: U- (\Gamma_e \cup \Gamma_d) \rightarrow \R^3$ with respect 
to $\gamma$, normalized so that $\alpha(\bar b) = (0,0,0)$  (cf. Proposition \ref{prop:stitched_action}). From Lemma~\ref{neg:ex}, there exists a primitive $\eta$ of $\omega$, 
such that for every $b=(b_1, b_2, b_3) \in  U- (\Gamma_e \cup \Gamma_d)$ we have
\[ \alpha(b) = \left( - \int_{\gamma_1(b)} \eta, \  - \int_{\gamma_2(b)} \eta, \  
                                       - \int_{\gamma_3(b)} \eta \right),\]
where $\gamma_j(b)$ is a cycle in $F_b$ representing $\gamma_j$. Clearly $\alpha$ 
is well defined and continuous on $U - (\Gamma_d \cup \Gamma_e)$. Actually, we have:
\begin{lem} The action coordinates map $\alpha$ extends continuously to
$\R^3$. 
\end{lem}
\begin{proof} We apply a similar argument to the one used in the case of the positive fibre 
(see Proposition~\ref{prop:pos_simple}). Clearly, since $\gamma_1$ is represented by the orbits of the 
$S^1$ action 
\[  - \int_{\gamma_1(b)} \eta = b_1, \]
which is continuous. We now prove that, for $j=2,3$ 
\begin{equation} \label{n:aj}
 \alpha_j(b) = - \int_{\gamma_j(b)} \eta 
\end{equation}
extends continuously to points in $\Gamma_d$ or in $\Gamma_e$. As we did in Proposition~\ref{prop:pos_simple},
we can think of $\alpha_j(b)$ as
\[ \alpha_j(b) = \int_{S} \omega, \]
where $S$ is a surface spanned by the cycles $\gamma_j(b')$ as $b'$ moves along a curve 
joining $\bar b$ and $b$. Suppose $b \in \Gamma_e$ (or $\Gamma_d$), then we need to show that $\alpha_j(b)$ 
is independent of the curve from $\bar b$ to $b$, or equivalently that
\[ \int_{S_1 - S_2} \omega = 0, \]
where $S_1$ and $S_2$ are the surfaces corresponding to two different paths from $\bar b$ to $b$.
The boundary $\partial(S_1 - S_2)$ is determined by monodromy. It is easy to see that 
$\partial(S_1 - S_2)$ is a multiple of $ \gamma_1(b)$, therefore for some integer $k$ we have 
\[ \int_{S_1 - S_2} \omega = - \int_{\partial (S_1 - S_2)} \eta = k \int_{\gamma_1(b)} \eta = 0, \]
where the last equality follows from the fact that $b \in \Gamma_d$ or $\Gamma_e$. To show that $\alpha$ extends continuously also to points of $\Delta$ we can argue that 
(\ref{n:aj}) makes sense also over singular fibres, since both $\eta$ and $\gamma_j(b)$ are
well defined when $b \in \Delta$.
\end{proof} 

We also have:
\begin{lem}
The map $\alpha: \R^3 \rightarrow \R^3$ is a homeomorphism onto its image.
\end{lem}
\begin{proof}
Since $\alpha_1(b) = b_1$, it is enough to show that, if for fixed $t \in \R$ we 
let $U_t = \{ b_1 = t \}$, then $\alpha_t = \alpha|_{U_t}$ is a bijection onto its image. 
If $\lambda_2$ and $\lambda_3$ are the periods of the fibration corresponding to $\gamma_2$ 
and $\gamma_3$, then $\alpha_t$ is computed by taking primitives of $\lambda_2|_{U_t}$ and 
$\lambda_3|_{U_t}$. If we let $X_t$ denote the symplectic reduction of $X$ at $t$ and 
$G_t: X_t \rightarrow \R^2$ the reduced fibration, then it is not difficult to see 
that $\lambda_2|_{U_t}$ and $\lambda_3|_{U_t}$ are in fact periods of $G_t$ (cf. \cite{CB-M-torino}{Lemma 5.9}). Now the conclusion follows by simply observing that 
$G_t$ is a proper Lagrangian submersion, i.e. an integrable system. The argument works 
also when $t=0$. 

An explicit computation of the periods was done in \cite{CB-M-torino}{Proposition 5.10} for the fibration in Example~\ref{thin leg}. There we found that 
\begin{eqnarray}
\lambda_2 & = & \beta_1 \, db_1 - e^{2b_2} db_2, \nonumber \\
 \label{neg.periods} \lambda_3 & = & \beta_2 \, db_1 - e^{2b_3} db_3,
\end{eqnarray}
where $\beta_1$ and $\beta_2$ are functions depending only on $b_1$.
The periods of the perturbed fibration obtained in Lemma~\ref{lem:smth_crit2} will have 
this same expression away from where the perturbation took place (i.e. away from 
the white region in Figure~\ref{fig:smoothing_2}), for example in a neighborhood of 
the codimension 1 part of $\Delta$. It is easy to see from this expression
of the periods that $\alpha$ extends continuously to $\Delta$ and that it is a bijection. 
\end{proof}

\begin{cor}
Let $\mathcal F$ be the fibration constructed in Lemma~\ref{lem:smth_crit2} and let 
$U = \R^3 - \Delta$. The stitched fibration $\mathcal F|_{U}$ is symplectically conjugate to a 
fibration constructed in Theorem~\ref{stitch:monodr_neg}.
\end{cor}

\begin{proof}
The fibrations constructed in Theorem~\ref{stitch:monodr_neg} have smooth Lagrangian 
sections and the action coordinates extend continuously to the whole base. 
Since $\mathcal F|_{U}$ also has a Lagrangian section (cf. Remarks~\ref{section:surv}) and 
the action coordinates extend continuously to the whole base, the statement easily follows from 
the results on stitched fibrations such as the existence of a normal form. The latter is 
found extending the maps $f^+$ and $f^-$ beyond all connected components of the seam and then 
using the Lagrangian section to normalize with the period map. 

\end{proof}

\subsection{Smoothing III}
\begin{figure}[!ht] 
\begin{center}
\input{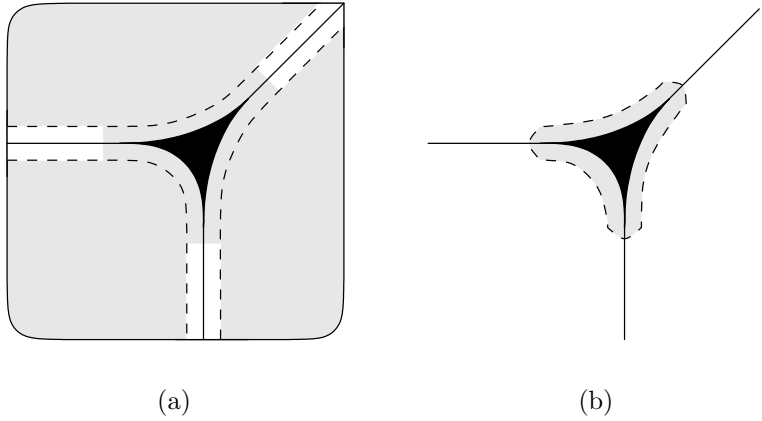}
\caption{Smoothing away from the legs.}\label{fig:smoothing_2}
\end{center}
\end{figure}

Now we show that the fibration in Example \ref{thin leg} can be perturbed to make it smooth on an even larger region. We consider the fibration $\mathcal F$ obtained in Lemma \ref{lem:smth_crit2} whose base is 
depicted in Figure \ref{fig:smoothing_2} (a). Over the white region complete smoothness was achieved. In the 
previous section we saw that over $U = \R^3 - \Delta$ the fibration is (symplectically 
conjugate to) a stitched Lagrangian fibration which can be constructed as  in Theorem~\ref{stitch:monodr_neg}. In this section we want to deform the invariants over 
each connected component of the seam so to achieve smoothness beyond the (planar) gray region in
Figure \ref{fig:smoothing_2} (b).

\begin{lem}\label{lem:smth_crit3} Let $\mathcal F$ be the fibration obtained in Lemma \ref{lem:smth_crit2}. There is a perturbation $\tilde{\mathcal F}$ of $\mathcal F$ such that:
\begin{itemize}
\item[(i)] $\tilde{\mathcal F}$ is topologically conjugate to $\mathcal F$;
\item[(ii)] there exists a submanifold with boundary $D \subset B$, homeomorphic to a closed disc in $\R^2$, with $\Delta \cap (B- D)$ consisting of three disjoint segments, such that $\tilde{\mathcal F}|_{\R^3 - D}$ is a smooth Lagrangian fibration.
\end{itemize}
\end{lem}
\begin{proof}
The proof follows the same lines of Lemma~\ref{lem:smth_crit2}. 
Assume that $\mathcal F|_{\R^3 - \Delta}$ has been constructed with 
Theorem~\ref{stitch:monodr_neg}. In particular the wall $\Gamma$ consists of the union of three 
disjoint sets, denoted $\Gamma_c$, $\Gamma_d$ and $\Gamma_e$. The corresponding components of the 
seam are $Z_c=f^{-1}(\Gamma_c)$, $Z_d=f^{-1}(\Gamma_c)$ and $Z_e=f^{-1}(\Gamma_e)$ with corresponding quotients denoted by $\bar Z_c$, $\bar Z_d$ and $\bar Z_e$. 
The invariants of $\mathcal F|_{\R^3 - \Delta}$ are given by sequences $\ell^c$, $\ell^d$ and 
$\ell^e$. 
In particular the first order invariants satisfy the integral conditions (\ref{st:neg_int}) 
\todo{changed sign of $m_1$ here too} with $m_1=-1$ and $m_2=1$.

\medskip
Over the same wall $\Gamma$ and seam $Z$, we could define another triple of invariants as follows.
Define $(\ell^c)'$ to be the zero sequence, while $(\ell^d)'$ and $(\ell^e)'$ to be sequences 
whose only non-zero terms are the first order ones, which we define to be
\todo{and here}
\[  (\ell_1^d)' =  - dy_2 \ \ \text{and} \ \  (\ell_1^e)' =  dy_3. \]
As we saw in Example~\ref{neg:as_stitched}, these choices of invariants give rise to a fake 
stitched fibration $\mathcal F'$ which is topologically conjugate to $\mathcal F|_{\R^3 - \Delta}$.

\medskip
Using Theorem~\ref{stitch:monodr_neg} we now construct a new stitched fibration with the same 
wall $\Gamma$ and seam $Z$ as $\mathcal F|_{\R^3 - \Delta}$, but whose invariants 
interpolate between those of $\mathcal F'$ and those of $\mathcal F|_{\R^3 - \Delta}$.
Let $A'$ be a small tubular neighborhood of $\Delta$ and denote $\bar A' = A' \cap \{ b_1 = 0 \}$.
Assume that $\bar A'$ is entirely contained in the region in Figure \ref{fig:smoothing_2} (a) 
delimited by the dotted lines. In particular we want the ends of $\bar A'$ to be contained in the white 
region where $\mathcal F$ is smooth. Let $A \subset A'$ be a smaller open neighborhood of $\Delta$ 
and denote $\bar A = A \cap \{ b_1 = 0 \}$. 
Let $\rho \in C^\infty(\Gamma)$ be a cut-off function which is 1 on $\bar A$ and $0$ on $\Gamma - \bar A'$. 
Define $\tilde \ell_k^c = (1 - \rho) (\ell_k^c)' + \rho \, \ell_k^c$ and similarly define 
$\tilde \ell_k^d$ and $\tilde \ell_k^e$.  
It follows from Theorem~\ref{stitch:monodr_neg} that the sequences 
$\tilde\ell_c$, $\tilde\ell_d$ and $\tilde \ell_e$ give rise to a stitched Lagrangian fibration 
$\tilde{\mathcal F}^o$ which is topologically conjugate to $\mathcal F|_{\R^3 - \Delta}$. Moreover
$\tilde{\mathcal F}^o|_{A - \Delta}$ and $\mathcal F |_{A - \Delta}$ are symplectically conjugate
so we can glue $\mathcal F|_{A}$ to $\tilde{\mathcal F}^o|_{A - \Delta}$ along $\mathcal F|_{A - \Delta}$.
This produces a piecewise smooth Lagrangian fibration $\tilde{\mathcal F}$ which is topologically 
conjugate to $\mathcal F$, moreover the chosen invariants guarantee that 
after a change of coordinates on the base $\tilde{\mathcal F}$ satisfies the smoothness condition $(ii)$. 
\end{proof}

The fibration $\tilde{\mathcal F}$ obtained via Lemma \ref{lem:smth_crit3} clearly satisfies properties 
(i) and (ii) of Definition \ref{lag:neg}, but finally we can also give 

\begin{proof}[Proof of Theorem~\ref{thm:neg_main}]
It only remains to show that $\mathcal {\tilde F}$ satisfies property (iii) of Definition~\ref{lag:neg},
but this immediately follows from the construction. In fact, $\tilde{\mathcal F}|_{\R^3 - A'}$ coincides
with the fibration described in Example~\ref{neg:as_stitched} restricted to a suitable neighborhood of 
the vertex. We observed that the latter fibration induces an affine structure on the base which 
is affine isomorphic to a negative vertex of Example~\ref{neg:aff} (or of Example~\ref{neg:affvar}). This concludes the proof.
\end{proof}

\section{The compactification.}\label{section:compact}

%
%

\subsection{The main theorem}

Finally, having completed the construction of the negative fibration, in this last section we prove the 
main result of the article. In order to give a correct statement of the theorem, we need first to make a few
observations.

\medskip
We start with a compact simple integral affine $3$-manifold with singularities $(B, \Delta, \mathscr A)$.
The goal is to symplectically compactify the torus bundle $X(B_0,\mathscr A)$ by gluing to it singular 
fibres. We have already seen in Section~\ref{pos:gener}, Proposition~\ref{prop:aff_glue} how the gluing 
of positive or generic singular fibres is quite straightforward. 
In the case of negative vertices we have seen that our construction gives a fibration whose discriminant 
locus contains components of type $\Delta_a$, i.e. of codimension $1$. For this reason around negative points 
one needs to replace $\Delta$ with a slightly perturbed discriminant locus containing components of type $\Delta_a$. 

\begin{figure}[ht]
\psfrag{A}{$\alpha$}
\begin{center}
\epsfig{file=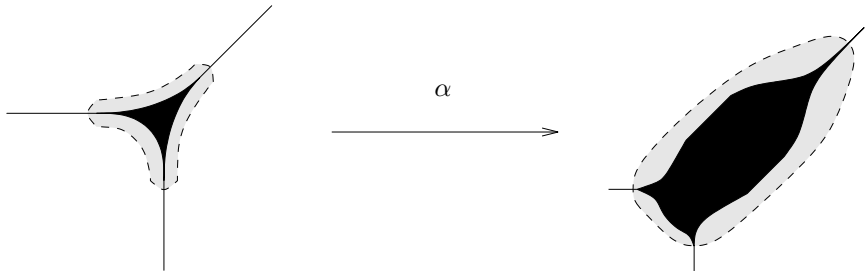, scale=.80}
\caption{The affine image of the amoeba with thin legs}\label{am:aff_shape}
\end{center}
\end{figure}

\medskip
Let us consider the fibration of Example~\ref{thin leg}. The periods of this fibration were computed in 
\cite{CB-M-torino} and they are given by formulas (\ref{neg.periods}). Let us consider the corresponding primitives (action coordinates) restricted to the plane $\{b_1 = 0\}$, which is the plane where the discriminant locus lies. We can easily see that the action coordinates
map $\alpha$ transforms the amoeba with thin legs into a slightly different shape, depicted  in in Figure~\ref{am:aff_shape}. 
This shape does not change much after we have done the smoothings of Lemmas~\ref{lem:smth_crit}, \ref{lem:smth_crit2} 
and \ref{lem:smth_crit3}, what may happen is that the codimension $2$ part --i.e. the legs-- may  become slightly curved. Nevertheless, it is 
not difficult to see that we can prolong the legs of a negative fibration so that they become straight toward their 
ends.  This can be done by gluing suitable generic-singular Lagrangian fibrations using the methods of 
Proposition~\ref{prop:def:leg}. 

\begin{figure}[ht]
\psfrag{-}{$p^-$}
\begin{center}
\epsfig{file=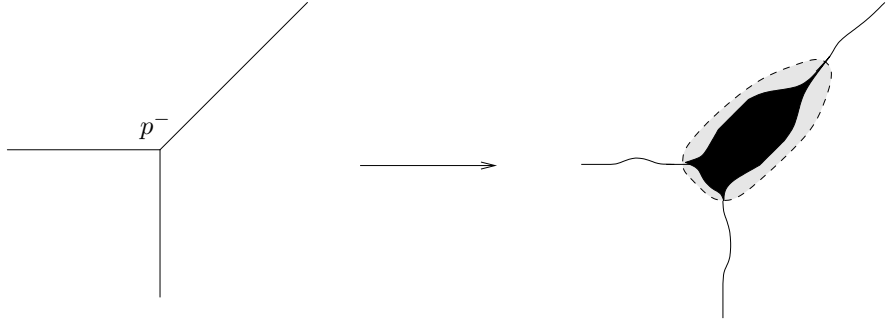, scale=.70}
\caption{A localized thickening of a negative vertex.}\label{local:thicken}
\end{center}
\end{figure}

\begin{defi}\label{loc:thick} Given a simple integral affine $3$-manifold with singularities $(B, \Delta, \mathscr A)$, 
all of whose negative vertices are straight (i.e. locally affine isomorphic to Example~\ref{neg:aff}), a 
\textit{localized thickening} of $\Delta$ is given by the data $(\Delta^{\blacklozenge}, 
\{ D_{p^-} \}_{p^- \in \mathcal N})$
where: 
\begin{itemize}
\item[(i)] $\Delta^{\blacklozenge}$ is the closed subset obtained from $\Delta$ after replacing a neighborhood
            of each negative vertex with a shape of the type depicted in Figure~\ref{local:thicken}. This replacement takes place in the plane corresponding to $\{ x_1 = 0 \}$ of the local model, Example~\ref{neg:aff}.
\item[(ii)] $\mathcal N$ is the set of negative vertices and for each $p^-\in \mathcal N$, $D_{p^-}$ is a submanifold of $B$, homeomorphic  
            to a disk and containing the codimension $1$ component of $\Delta^{\blacklozenge}$ around the 
            negative vertex $p^-$. 
            Moreover, $D_{p^-}$ is contained 
            in the plane $\{ x_1 = 0 \}$. We depict $D_{p^-}$ as the gray area in Figure~\ref{local:thicken}.
\end{itemize}
\end{defi}

The requirement that all negative vertices are straight is only to avoid unnecessary complications. 
Given a localized thickening of $\Delta$, define 
\[
B_{\blacklozenge} = B - \left( \Delta \cup \bigcup_{p^- \in \mathcal N} D_{p^-} \right).
\]
Clearly, the integral affine structure $\mathscr A$ on $B-\Delta$ restricts to an integral affine structure on $B_\blacklozenge$  which we denote by 
$\mathscr A_{\blacklozenge}$, therefore we can form the torus bundle $X( B_{\blacklozenge}, 
\mathscr A_{\blacklozenge})$. 

\medskip
Now we can state and prove the theorem:

\begin{thm} \label{the_symplectic_comp}
Given a compact simple integral affine $3$-manifold with singularities $(B, \Delta, \mathscr A)$,
all of whose negative vertices are straight (i.e. locally isomorphic to Example~\ref{neg:aff}),
there is a localized thickening $(\Delta_{\blacklozenge}, \{ D_{p^-} \}_{p^- \in \mathcal N})$ of $\Delta$ 
and a smooth, 
compact symplectic $6$-manifold $(X, \omega)$ together with a piecewise smooth Lagrangian fibration 
$f: X \rightarrow B$ such that 
\begin{itemize}
\item[(i)] $f$ is smooth except along $\bigcup_{p^- \in \mathcal N} \,  f^{-1}(D_{p^{-}})$;
\item[(ii)] the discriminant locus of $f$ is $\Delta_{\blacklozenge}$;
\item[(iii)] there is a commuting diagram
\begin{equation*} 
\begin{CD}
X(B_{\blacklozenge}, \mathscr A_{\blacklozenge})  @>\Psi>> X\\
@Vf_0VV  @VVfV\\
B_{\blacklozenge} @>\iota>> B
\end{CD}
\end{equation*}
where $\psi$ is a symplectomorphism and $\iota$ the inclusion;
\item[(iv)] over a neighborhood of a positive vertex of $\Delta_{\blacklozenge}$ the fibration is 
             positive, over a neighborhood of a point on an edge the fibration is generic-singular, over
             a neighborhood of $D_{p^-}$ the fibration is Lagrangian negative. 
\end{itemize}
\end{thm}

\begin{proof}
The proof is quite simple. First we glue positive fibrations over sufficiently small neighborhoods of 
positive vertices of $\Delta$ using Proposition~\ref{prop:aff_glue}. 
Now given a negative vertex $p^- \in \mathcal N$, we have that a neighborhood of $p^-$ is 
affine isomorphic to a neighborhood $U$ of zero in the local model Example~\ref{neg:aff}. 
Consider a negative Lagrangian fibration $\mathcal F^- = (X^-, \omega^-, f^-, B^-)$ 
(cf. Definition~\ref{lag:neg}), which we have constructed in Theorem~\ref{thm:neg_main}. 
The discriminant locus $\Delta^-$ of $f^-$ has the shape of an amoeba with thin legs and
there is a disc $D$ containing the codimension $1$ part of $\Delta^-$ such that $f^-$ is smooth 
except at points of $(f^-)^{-1}(D)$ (cf. part (i) and (ii) of Definition~\ref{lag:neg}). 
Moreover we may assume that $B^- - (\Delta^- \cup D)$ is affine isomorphic to 
$(U' - (D' \cup \Delta_{\tau}), \mathscr A_{\tau})$, where $U'$ is a neighborhood of $0$ in 
the affine manifold with singularities of Example~\ref{neg:affvar} and $D' \subset \{ x_1 = 0 \}$
contains $0$ and is homeomorphic to a disc (cf. point (iii) of Definition~\ref{lag:neg}). 
It may happen that $U'$ is too big for us to glue the Lagrangian negative fibration as it is.
However, if we replace $\omega^-$ with $\epsilon \, \omega^-$ for a sufficiently small $\epsilon > 0$, 
this has the effect of scaling the affine coordinates on the base by a factor of $\epsilon$
 (i.e. of making the amoeba as small as we please). 
Therefore we may assume that $U' \subset U$. Moreover, we may also assume that the legs of $\Delta^-$ (in affine coordinates) are straight towards their ends, i.e. they coincide with the legs of 
$\Delta$ outside an open subset $U''$ such that 
$D' \subset \bar U '' \subset U'$.  The localized thickening $\Delta_{\blacklozenge}$ of $\Delta$ around $p^-$ 
consists in replacing $U' \cup \Delta$ with $\Delta^-$ and defining $D_{p^-} = D'$. The affine structure 
$\mathscr A_{\blacklozenge}$ is inherited from $\mathscr A$. This can be done at every negative 
vertex $p^{-}$. Tautologically, we have that  
$X(U' - (D_{p^-} \cup \Delta_{\blacklozenge}), \mathscr A_{\blacklozenge})$ is symplectically conjugate 
to $(f^-)^{-1}(B^{-} - (\Delta^{-} \cup D))$ and therefore we can glue $X^-$ to 
$X(B_{\blacklozenge}, \mathscr A_{\blacklozenge})$.

\medskip
Finally, now that singular fibres have been glued on top of all vertices, it only remains to glue 
generic-singular fibres along the edges. This can be easily done by applying directly 
Proposition~\ref{prop:def:leg}, notice in fact that Lagrangian negative fibrations are smooth and 
generic-singular towards the ends of the legs. 
\end{proof}

We remark that the manifolds we obtain with this theorem are diffeomorphic to Gross' semi-stable 
compactifications of Theorem~\ref{thm. semi-stable}. Also, as a corollary of this construction we have
\begin{cor} A smooth quintic $X$ in $\PP^4$ has a symplectic form $\omega$ with a piecewise smooth
Lagrangian fibration $f: X \rightarrow S^3$.
\end{cor}
\begin{proof}
If we apply Theorem~\ref{the_symplectic_comp} to Example~\ref{quint:aff} we obtain 
a symplectic manifold $X$ with a piecewise smooth Lagrangian fibration $f: X \rightarrow S^3$. 
By Gross' Theorem~\ref{quint:mirr}, $X$ is homeomorphic to a non-singular quintic.
\end{proof}

We do not know whether the symplectic manifold $(X, \omega)$ obtained in this corollary is actually symplectomorphic to a quintic with a K\"ahler form, although we conjecture it is.

%

\bibliographystyle{plain}

\vspace{2cm}

\begin{flushleft}

Ricardo Casta\~no-Bernard \\
Mathematics Department\\
Kansas State University, \\
138 Cardwell Hall\\
Manhattan KS 66502 USA\\
e-mail: \texttt{rcastano@math.ksu.edu}\\ \
\\ \ 
\  
\\ 
Diego~Matessi\\
Dipartimento
 di Scienze e Tecnologie Avanzate\\
Universit\`{a} del Piemonte Orientale\\
Via Bellini 25/G\\
 I-15100 Alessandria, Italy\\
e-mail: \texttt{matessi@unipmn.it}\\

\end{flushleft}

\end{document}